\documentclass[12pt]{article}
\usepackage{amsmath, amsthm, amssymb}

\theoremstyle{plain}
  \newtheorem{thm}{Theorem}[section]
  \newtheorem{lem}[thm]{Lemma}
  \newtheorem{cor}[thm]{Corollary}
  \newtheorem{prop}[thm]{Proposition}

\theoremstyle{definition}
  \newtheorem{defi}[thm]{Definition}
  \newtheorem{ex}[thm]{Example}
  \newtheorem{rem}[thm]{Remark}

\numberwithin{equation}{section}

\newcounter{assump}

\newcounter{assumpt}

\newcommand{\tr}{\operatorname{tr}}

\newcommand{\eps}{\epsilon}
\newcommand{\esup}{\operatornamewithlimits{ess.sup}}
\newcommand{\argmin}{\operatornamewithlimits{arg \, min}}
\newcommand{\argmax}{\operatornamewithlimits{arg \, max}}

\renewcommand{\labelenumi}{(\roman{enumi})}

\title{Max-plus Stochastic Control and Risk-sensitivity}
\author{Wendell H.\,Fleming\footnote{Division of Applied Mathematics and Lefschetz Center
	 for Dynamical Systems, Brown University, Providence, RI 02912, USA. whf@cfm.brown.edu.},	 
	Hidehiro Kaise\footnote{Graduate School of Information Science, Nagoya University, Furo-cho,
  Chikusa-ku, Nagoya 464-8601, Japan. kaise@is.nagoya-u.ac.jp}
  and Shuenn-Jyi Sheu\footnote{Institute of Mathematics, Academia Sinica, Nankang, Taipei
 11529, Taiwan, Republic of China. sheusj@math.sinica.edu.tw.}
  }
\date{January 9, 2009}

\begin{document}
\maketitle

\begin{abstract}
In the Maslov idempotent probability calculus, expectations of random variables are defined so as to be linear with respect to max-plus addition and scalar multiplication. This paper considers control problems in which the objective is to minimize the max-plus expectation of some max-plus additive running cost. Such problems arise naturally as limits of some types of risk sensitive stochastic control problems. The value function is a viscosity solution to a quasivariational inequality (QVI) of dynamic programming. Equivalence of this QVI to a nonlinear parabolic PDE with discontinuous Hamiltonian is used to prove a comparison theorem for viscosity sub- and super-solutions. An example from math finance is given, and an application in nonlinear $H$-infinity control is sketched.


\end{abstract}

\section{Introduction}
A wide variety of asymptotic problems, including large deviations for stochastic processes, can be considered in the framework of the Maslov idempotent probability calculus. In this theory, probabilities are assigned which are additive with respect to ``max-plus'' addition and expectations are defined so as to be linear with respect to max-plus addition and scalar multiplication. For this reason we use the term ``max-plus probability'' instead of ``idempotent probability.'' There is extensive literature on max-plus probability and max-plus stochastic processes. See 
\cite{AQV94}, \cite{BCOQ92}, \cite{DMD99}, \cite{F02}, \cite{MS92} and references cited there. The max-plus framework is also important for certain problems in discrete mathematics and in computer science applications. See 
\cite{AGW05}, \cite{BCOQ92} for example.

To provide background for the results of this paper, let us begin by mentioning two results from the Freidlin-Wentzell theory of large deviations for small random perturbations of dynamical systems 
\cite{FrW84}. Consider a finite time interval $[t,T]$ and $X(s)$ which satisfies the stochastic differential equation (SDE)
\begin{equation}
	\left\{
	\begin{aligned}
		dX(s) &= f(X(s)) ds +\theta^{-1/2} \sigma (X(s)) dW(s), \ t \leq s \leq T, \\
		X(t) &= x \in \mathbb{R}^n.	
	\end{aligned}
	\right.
\label{diff}
\end{equation}
Here $\theta > 0$ is a ``large'' parameter and $W(s)$ is a $d$-dimensional Brownian motion. Sometimes we write $X(s) = X_\theta(s)$ to emphasize dependence on $\theta$. Asymptotic large deviations results as $\theta \rightarrow \infty$ are typically described through a deterministic optimization problem. In the limit problem, $X_\theta(s)$ in \eqref{diff} is replaced by $x(s)$, which satisfies an ordinary differential equation (ODE)
\begin{equation}
	\left\{
	\begin{aligned}
		\frac{dx}{ds}(s) &= f(x(s)) +\sigma (x(s)) v(s), \ t \leq s \leq T, \\
		x(t) &= x.
	\end{aligned}	
	\right.
\label{dyn}
\end{equation}
The unknown function $v(\cdot)$ is a ``control,'' with $v(\cdot) \in L^2 ([t,T]; \mathbb{R}^d).$ In the language of nonlinear $H$-infinity control theory, $v(s)$ is called a ``disturbance'' \cite{HJ99}. The disturbance control is chosen to maximize an expression of the form $\Phi (x (\cdot)) - \frac{1}{2} \| v(\cdot) \|_2$ with $\|  \cdot \|_2$ the $L^2$-norm. The maximum is the max-plus expectation $E^+[\Phi(x(\cdot))]$ as defined in Section 2.

We may consider, in particular, the following two kinds of choices for $\Phi$.

\vspace*{.1in}

\noindent {\bf Case 1}
$$ \Phi_1(x(\cdot)) = \int_t^T l (x(s))ds + G(x(T)).
$$
The term $G(x(T))$ is called a terminal cost. When $G=0$, we say that $\Phi_1$ is a ``max-plus multiplicative'' running cost. Under suitable assumptions on $f$, $\sigma$ and $l$, the Freidlin-Wentzell theory implies that
\begin{equation}
	\lim_{\theta\rightarrow \infty} \theta^{-1} \log E[\exp \{\theta \Phi_1 (X_\theta(\cdot))\}]
	= E^+ [\Phi_1 (x(\cdot)) ].
\end{equation}

\vspace*{.1in}

\noindent {\bf Case 2}
$$
\Phi_2 (x(\cdot)) = \int^{\oplus}_{[t,T]} l (x(s))ds = \max_{s\in [t,T]} l (x(s)).
$$
This is the case of ``max-plus additive'' running cost. In Case 2, the large deviations result (see \cite{F04}) is
\begin{equation}\
\lim_{\theta \rightarrow \infty}\theta^{-1} \log
 E \left[ \int_t^T\exp \{\theta l (X_\theta(s))\}ds \right]
 = E^+ [\Phi_2 (x(\cdot))].
 \label{LD}
\end{equation}

In this paper, we are concerned with problems in which the time evolution of the state $x(s)$ depends on a control $u(s)$. Thus, equation \eqref{dyn} above will be replaced by the ODE \eqref{system}. The case of max-plus multiplicative running costs has already been studied using methods of risk sensitive stochastic control theory and differential games. See \cite[Chap.\,6]{FS06}, \cite{F02}. We consider the problem of choosing the control $u(s)$ to minimize
\begin{equation}
	E^+\left[ \int^\oplus_{[t,T]} l (x(s), u(s))ds\right]
	= E^+\left[ \esup_{t\leq s\leq T} l (x(s),u(s)) \right].
\label{m-p_add}
\end{equation}
We call this a ``max-plus stochastic control problem with max-plus additive running cost.''

The control $u(s)$ at time $s$ is to be chosen using information about disturbances $v(r)$ for times $r<s$. A precise definition of the class $\Gamma(t,T)$ of admissible strategies is given in Section 3. $\Gamma(t,T)$ consists of Elliott-Kalton strategies with the additional properties (S1), (S2). Another possible class consists of those Elliott-Kalton strategies which are ``strictly  progressive,'' as defined in \cite[Section 11.9]{FS06} and \cite{KS05}. See also Section 4.3.

The max-plus stochastic control problem is studied using the method of dynamic programming. The associated value function $V(t,x)$ is defined by \eqref{value}. It satisfies a dynamic programming principle (Theorem \ref{thm_DPP}), which is proved by standard arguments using properties of max-plus conditional expectations. The dynamic programming equation for $V$ takes the form of a quasivariational inequality (QVI) \eqref{QVI}. It is shown later (Theorem \ref{ch_value}) that $V$ is the unique bounded, Lipschitz continuous viscosity solution to \eqref{QVI} with boundary  data \eqref{QVI_T} at time $T$. 

In  Section 4, the QVI \eqref{QVI} is shown to be equivalent to the nonlinear parabolic PDE \eqref{u_Isaacs}, which involves a discontinuous Hamiltonian. In the special case when there are no disturbances in the model $(\sigma = 0)$, the PDE \eqref{u_Isaacs} is of a kind considered by Barron-Ishii \cite{BI89}. The treatment of viscosity subsolutions and supersolutions in Section 4, and the proof of the comparison Theorem \ref{comp}, make use of similar ideas in \cite{BI89}.

Section 4.2 considers a nonlinear two time parameter semigroup associated with the max-plus control problem. This semigroup is expressed in terms of a family of operators $F_{t,r}$, which are related to max-plus linear operators. The semigroup property is a consequence of the dynamic programming principle. The operator $\frac{\partial V}{\partial t} + {\cal H}$ in \eqref{u_Isaacs} has an interpretation as the generator of this semigroup. See Theorem \ref{cal_gen}.

In the max-plus control problem, there are actually two controls. One is the control $u(s)$, chosen to minimize \eqref{m-p_add}. The other (maximizing) control is $v(s)$, which enters through the max-plus expectation $E^+$. Although our problem has a differential game interpretation, we make no use of results about differential games. See remarks in Section 4.3.

In Section 6, we consider a stochastic control version of the model, in which the ODE \eqref{system} for the state dynamics is replaced by the SDE \eqref{stoch_sys} and the goal is to minimize the expectation of the risk-sensitive criterion in \eqref{r-s_cr}. Let $\Psi_\theta(t,x)$ be the value function for this problem. Theorem \ref{r-s_value} states that
\begin{equation}
	\lim_{\theta \rightarrow \infty} \theta^{-1} \log \Psi_\theta (t,x) = V(t,x).
\label{r-s_asym}
\end{equation}
The corresponding large deviations result is \eqref{LD}. Theorem \ref{r-s_value} is proved using a version of the Barles-Perthame method for viscosity solutions.

As an example to illustrate the limit in \eqref{r-s_asym}, we consider in Section 7 the classical Merton optimal investment-consumption problem in mathematical finance. For the Merton problem, $\Psi_\theta(t,x)$, $V(t,x)$ and the corresponding optimal controls can be found by explicit calculations.

In Section 8 we consider some infinite time horizon problems. The interest is in inequalities of the form \eqref{dis_ineq} which hold on every finite time interval $[0,T]$. Such inequalities have an interpretation in nonlinear $H$-infinity control theory, with max-plus additive running cost. The discussion follows mainly \cite[Section 8]{F04}.

\section{Max-plus stochastic control}
\subsection{Preliminaries on max-plus probability}

We start with reviewing some notions and facts from max-plus probability theory
which will be needed later.
The readers should refer to \cite{AQV94}, \cite{BCOQ92}, \cite{DMD99}, \cite{F04},
\cite{HJ99}, \cite{MS92} for more details. 
 
Let us consider extended reals $\mathbb{R}^{-}=\mathbb{R} \cup \{ -\infty \}$.
For $a$, $b \in \mathbb{R}^-$,  we define new addition and multiplication
by
\[
	a \oplus b = \max \{ a, b \}, \ a \otimes b =a+b.
\]
$\mathbb{R}^-$ with these new operations is  called max-plus algebra
and it satisfies all the axioms of rings except for the existence of additive inverse.
In addition, the idempotency $a \oplus a =a$ holds.

To introduce the notions of max-plus probability, we focus on a particular case
discussed in the present paper.
Let $t \in [0,T]$ and $\Omega=\Omega_{t,T}=L^2([t,T];\mathbb{R}^d)$ 
be a sample space. We use the notation $L^{2}[t,T]$ for $L^2([t,T]; \mathbb{R}^d)$
if the dimension is clear from the context.
On this sample space, we consider  max-plus probability density $Q : \Omega \to \mathbb{R}^{-}$
\[
	Q(v)=-\frac{1}{2}\int_t^T |v(s)|^2 ds, \ v \in \Omega.
\]
Then,  the max-plus probability $P^+$ for $A \subset \Omega$ is defined by
\[
	P^+ (A)= \sup_{ v \in A} Q(v)
		=\sup_{v \in A} \left\{ -\frac{1}{2}\int_t^T |v(s)|^2 ds \right\}.
\]
The supremum on empty set is understood to be $-\infty$.
We call $(\Omega, P^+)$ (or $(\Omega, Q)$)  a max-plus probability space.

Let $Z : \Omega \to \mathbb{R}^-$ be a random variable.
Here we do not require any measurability conditions.
The max-plus expectation of $Z$ is 
\[
	E^{+}[Z] = \sup_{v \in \Omega} \{ Z(v) \otimes Q(v) \}
		=\sup_{v \in L^2[t,T]} \left\{ Z(v) -\frac{1}{2}\int_t^T |v(s)|^2 ds \right\}.
\]
$Z$ is max-plus integrable if $ E^{+}[Z] < \infty$.
It is easily seen that the max-plus expectation is linear under max-plus algebra:
for max-plus integrable random variables $Z$, $Y$ and $a \in \mathbb{R}^-$,
\[
	E^+ [Z \oplus Y] = E^+ [Z] \oplus E^+ [Y], \
	E^+ [a \otimes Z] = a \otimes E^+ [Z].
\]
If $Z \leq Y$, then
\[
	E^+ [Z] \leq E^+ [Y].
\]

To consider max-plus conditional expectations under this special probability space,
let $0 \leq t <r < T$ and denote $\Omega_1$, $\Omega_2$ by
\[
	\Omega_1=\Omega_{t,r}=L^2[t,r], \
	\Omega_2=\Omega_{r,T}=L^2[r,T].
\]
On the product space $\Omega_1 \times \Omega_2$, we define max-plus probability density
$Q_1 \otimes Q_2$ by
\[
	(Q_1 \otimes Q_2) (v_1, v_2)=Q_1(v_1) \otimes Q_2(v_2), \ (v_1, v_2) \in \Omega_1 \times \Omega_2,
\]
where
\[
	Q_1(v_1)=-\frac{1}{2}\int_t^r |v_1(s)|^2 ds, \
	Q_2(v_2)=-\frac{1}{2}\int_r^T|v_2(s)|^2ds.
\]
Note that if we set $v_1=v|_{[t,r]}$, $v_2=v|_{[r,T]}$ for $v \in \Omega$,
\[
	Q(v)=Q_1 (v_1) \otimes Q_2 (v_2).
\]
Thus,  ($\Omega,Q)$ can be identified with
$(\Omega_1 \times \Omega_2, Q_1 \otimes Q_2)$.

Let $Z(v)=Z(v_1, v_2)$ be a random variable on $\Omega$.
For $v_1 \in \Omega_1$, the max-plus conditional expectation of $Z$ under $v_1$ is given by
\[
	E^+ [Z | v_1] = \sup_{v_2 \in \Omega_2} [Z(v_1, v_2)\otimes Q_2(v_2) ].
\]
If $Z$ has the form $Z(v_1,v_2)=Z_1(v_1) \oplus Z_2(v_1, v_2)$ for
some $Z_1: \Omega_1 \to \mathbb{R}^-$ and $Z_2: \Omega \to \mathbb{R}^-$,
it is seen that
\begin{equation}
	E^+ [ Z ] = E^+ [Z_1\oplus E^+ [Z_2(v_1, \cdot ) | v_1]].
\label{tower}
\end{equation}

\subsection {Problem formulation and verification theorem}

We consider max-plus control problems on a finite time interval $[t,T]$.
The final time $T$ is fixed throughout the paper, and the initial time satisfies $0 \leq t <T$.
Let us consider the system governed by 
\begin{equation} \label{system}
\left\{
\begin{aligned}
	\frac{dx}{ds}(s) &= f(x(s),u(s))+\sigma(x(s),u(s))v(s), \ t \leq s \leq T, \\
	x(t) &= x \in \mathbb{R}^n,
\end{aligned}
\right.
\end{equation}
where 
$U \subset \mathbb{R}^m$, $f:\mathbb{R}^n \times U \to \mathbb{R}^n$,
$\sigma: \mathbb{R}^n \times U \to M(n,d)$, 
$M(n,d)$ is the set of $n \times d$ matrices.
$x(s)$ is the state of the system,
$u \in L^\infty([t,T]; U)$ is a control and
$v \in L^2[t,T]$ is a disturbance (or uncertainty) in the system.
Equation \eqref{system} is an ordinary differential equation with two parameters.
However,  in terms of max-plus diffusion processes,
we can regard \eqref{system} as a controlled max-plus stochastic differential equation
under $(\Omega, Q)$
(cf. \cite{F04}). 

When we choose a control at a certain time, it is natural to require that
the decision has to be made by past information of the disturbance.
It can be realized
by the notion of Elliott-Kalton strategy in the theory of differential games
(cf. \cite{EK72}).
Let $\alpha: L^2[t,T] \to L^\infty([t,T];U)$. 
$\alpha$ is
an Elliott-Kalton strategy from $L^2[t,T]$ into $L^\infty([t,T];U)$
if $\alpha$ satisfies the following:
Let $v$, $\tilde{v} \in L^2[t,T]$ and $t \leq s \leq T$.
\begin{equation}
	\text{If }v=\tilde{v} \text{ a.e.\,on }[t,s],
	\text{ then }\alpha[v] =\alpha[\tilde{v}] \text{ a.e.\,on }[t,s].
\label{EK-cond}
\end{equation}
We denote by $\Gamma_{EK}(t,T)$ the set of Elliott-Kalton strategies
from $L^2[t,T]$ into $L^\infty([t,T];U)$.
If we choose $\alpha \in \Gamma_{EK}(t,T)$,  \eqref{system} becomes
\begin{equation} \label{mp-SDE}
\left\{
\begin{aligned}
	\frac{dx}{ds}(s) &= f(x(s),\alpha[v](s))
		+\sigma(x(s),\alpha[v](s))v(s), \ t \leq s \leq T, \\
	x(t) &= x \in \mathbb{R}^n
\end{aligned}
\right.
\end{equation}
for each disturbance $v \in L^2 [t,T]$.

The goal in max-plus stochastic control is to minimize the max-plus expectation
of some criterion $\mathcal{J}$
on a suitable subclass of Elliott-Kalton strategies.
For the controlled system \eqref{mp-SDE}, there can be three natural criteria:
\begin{enumerate}
\item Terminal cost: For $\Phi : \mathbb{R}^n \to \mathbb{R}$,
\[
	\mathcal{J}=\Phi(x(T)).
\]
\item  Max-plus multiplicative running cost:  For $l: \mathbb{R}^n \times U \to \mathbb{R}$,
\[
	\mathcal{J}=\int_t^T l(x(s),\alpha[v](s))ds.
\]
\item Max-plus additive running cost: For $l: \mathbb{R}^n \times U \to \mathbb{R}$,
\[
	\mathcal{J}=\int_{[t,T]}^{\oplus} l(x(s),\alpha[v](s))ds
	\equiv \esup_{s \in [t,T]} l(x(s),\alpha[v](s)).
\]
\end{enumerate}
(i) and (ii) are considered in \cite[Section 11.7]{FS06} and \cite{KS05}.
In the present paper, we shall discuss the criterion of type (iii).

Unless otherwise stated, we assume that the following conditions hold:
\renewcommand{\labelenumi}{(A\arabic{enumi})}
\begin{enumerate}
\item $U \subset \mathbb{R}^m$ is compact.

\item $f$ and $\sigma$ are of $C^1$ on $\mathbb{R}^n \times U$.
	$f$, $\sigma$ and their derivatives $f_x$, $\sigma_x$ in $x$ are bounded on
	$\mathbb{R}^n \times U$.

\item $l$ is of $C^1$ on $\mathbb{R}^n \times U$.
$l$, $l_x$ and $l_u$ are bounded on $\mathbb{R}^n \times U$.
\end{enumerate}
\renewcommand{\labelenumi}{(\roman{enumi})}
Under (A2), there exists a unique solution of \eqref{system} for
any control $u \in L^\infty([t,T]; U)$ and disturbance $v \in L^2 [t,T]$.

We formulate our max-plus control problem more specifically.
Let a subclass $\Gamma (t,T) \subset \Gamma_{EK}(t,T)$ be given.
For $\alpha \in \Gamma (t,T)$,
consider the max-plus expectation of the max-plus additive running cost criterion:
\begin{equation}
\begin{aligned}
	J(t,x; \alpha)
		&=E^{+}_{tx} \left[ \int_{[t,T]}^{\oplus} l(x(s),\alpha[v](s))ds \right] \\
		&=\sup_{v \in L^2[t,T]}
		\left\{ \int_{[t,T]}^\oplus l(x(s),\alpha[v](s))ds -\frac{1}{2}\int_t^T |v(s)|^2ds \right\},
\end{aligned}
\label{exp_mp-add}
\end{equation}
where $x(s)$ is the solution of \eqref{mp-SDE}.  
We indicate the dependence on the initial condition of the system  by the subscript $tx$
of $E^{+}_{tx}$. 
Our concern in max-plus stochastic control is to minimize $J(t,x; \alpha)$ on $\Gamma(t,T)$.
Thus, the value function associated with strategy class $\Gamma(t,T)$ is defined by
\begin{equation}
	V(t,x)
	=\inf_{\alpha \in \Gamma(t,T)} J(t,x;\alpha)
	=\inf_{\alpha \in \Gamma(t,T)} E^{+}_{tx}\left[ \int_{[t,T]}^\oplus l(x(s),\alpha[v](s))ds \right]. 
\label{g-value}
\end{equation}

For given class $\Gamma(t,T) \subset \Gamma_{EK}(t,T)$,
it is a fundamental problem to characterize $V(t,x)$ as a (unique) solution
of the associated dynamic programming equation (DPE).
In order to guess the DPE for $V(t,x)$, consider $J(t,x;\alpha)$ with
constant strategy $\alpha[v](s)\equiv u$ for $u \in U$, \textit{i.e}, let us define $V^u(t,x)$ by
\begin{equation}
	V^u(t,x) = E^{+}_{tx} \left[ \int_{[t,T]}^\oplus l(x(s),u) ds \right].
\label{value_u}
\end{equation}
In \cite{F04}, it is proved that under (A1)--(A3), $V^u(t,x)$ is the unique bounded Lipschitz
continuous viscosity solution of
\begin{equation}
\begin{gathered}
	\max \left\{ \frac{\partial V^u}{\partial t}+H^u(x,\nabla V^u(t,x)), l(x,u)-V^u(t,x) \right\}=0,
	\ (t,x) \in (0,T) \times \mathbb{R}^n, \\
	V^u(T,x)= l(x,u), \ x \in \mathbb{R}^n,
\end{gathered}
\label{QV_u}
\end{equation}
where for $x$, $p \in \mathbb{R}^n$ and $u \in U$,
\begin{equation}
\begin{aligned}
	H^u(x,p) &=\max_{v \in \mathbb{R}^d}
	\left\{ (f(x,u)+\sigma (x,u)v)\cdot p -\frac{1}{2}|v|^2 \right\} \\
	&=\frac{1}{2}a(x,u)p\cdot p+f(x,u) \cdot p,
	\ 	a(x,u)=\sigma(x,u) \sigma (x,u)^T. 
\end{aligned}
\label{H^u} 
\end{equation}
Note that
$\sigma(x,u)^T p$ attains the maximum in $H^u(x,p)$, \textit{i.e.},
\begin{equation}
	\sigma (x,u)^T p \in 
	\argmax_{v \in \mathbb{R}^d}
		\left\{ (f(x,u)+\sigma (x,u)v)\cdot p -\frac{1}{2}|v|^2 \right\}.
\label{argmax_H^u}
\end{equation}
We expect that a value function $V(t,x)$ for some strategy class may satisfy the following
quasivariational inequality (QVI):
\begin{gather}
	\min_{u \in U} \max \left\{ \frac{\partial V}{\partial t}+H^u(x,\nabla V(t,x)), l(x,u)-V(t,x) \right\}=0,
	\ (t,x) \in (0,T) \times \mathbb{R}^n, \label{QVI} \\
	V(T,x)= \min_{u \in U}l(x,u), \ x \in \mathbb{R}^n. \label{QVI_T}
\end{gather}
We will see that \eqref{QVI} with \eqref{QVI_T} is the correct DPE
for our max-plus control problem in \eqref{g-value} if $\Gamma (t,T)$ is properly chosen.
In Section 4, we will show that \eqref{QVI} is equivalent to a nonlinear parabolic PDE
with discontinuous Hamiltonian. When $\Gamma(t,T)$ is chosen to satisfy (S1), (S2) in Section 3, then the value function is the unique bounded Lipschitz viscosity solution of \eqref{QVI}, \eqref{QVI_T}. 
See Theorem \ref{ch_value}.


We conclude this section by considering strategies determined by Markov control policies
$\underline{u}$.
Let $\underline{u}(s,y)$ be a Lipschitz continuous function of
$(s,y) \in [t,T] \times \mathbb{R}^n$ into $U$.
The corresponding strategy $\alpha^{\underline{u}}$ is defined by
\begin{equation}
	\left\{
	\begin{aligned}
		\frac{dx}{ds}(s) &= f(x(s), \underline{u}(s,x(s)))
				+\sigma (x(s), \underline{u}(s,x(s))) v(s), \ t \leq s \leq T, \\
		x(t) &= x
	\end{aligned}
	\right.
\label{sys_fdbk}
\end{equation}
and $\alpha^{\underline{u}}[v](s)=\underline{u}(s,x(s))$.
$\underline{u}$ is called a Markov control policy.

\begin{thm} \label{ver_thm}
Let $W(t,x)$ be a $C^1$-solution of \eqref{QVI} and \eqref{QVI_T}
such that $\nabla W(t,x)$ satisfies a uniform Lipschitz condition on $x$. Then: \\
\textup{(a)} $W(t,x) \leq J(t,x;\alpha^{\underline{u}})$ for every Lipschitz Markov control policy
$\underline{u}$; \\
\textup{(b)} If there exists a Lipschitz Markov control policy $\underline{u}^\ast$ such that
for any $(s,y) \in [t,T] \times \mathbb{R}^n$
\begin{equation}
	\underline{u}^\ast(s,y)
	\in \argmin_{u \in U}
	\max \left\{ \frac{\partial W}{\partial s}(s,y)+H^u(y,\nabla W(s,y)), l(y,u)-W(s,y) \right\},
\label{arg_min}
\end{equation}
then $W(t,x)=J(t,x; \alpha^{\underline{u}^\ast})$.
\end{thm}
\noindent
\textit{Proof.} 
We shall first show (a). By the fundamental theorem of calculus,
\begin{multline}
	W(s,x(s)) -\frac{1}{2}\int_t^s |v(r)|^2 dr\\
	= W(t,x) +\int_t^s  \bigg\{ \frac{\partial W}{\partial r}(r,x(r))
			+(f(x(r), \alpha^{\underline{u}}[v](r))
			+\sigma (x(r),  \alpha^{\underline{u}}[v](r))v(r)) \cdot \nabla W(r,x(r))  \\
				 -\frac{1}{2} |v(r)|^2 \bigg\} dr.
\label{fund_v}
\end{multline}
where $x(r)$ is the solution of \eqref{sys_fdbk}.

Let us consider the following closed system on $[t,T]$:
\begin{equation}
	\left\{
	\begin{aligned}
		\frac{d\hat{x}}{ds}(s) &= f(\hat{x}(s), \underline{u}(s,\hat{x}(s)))
				+\sigma (\hat{x}(s), \underline{u}(s,\hat{x}(s)))
						\sigma (\hat{x}(s), \underline{u}(s,\hat{x}(s)))^T \nabla W(s,\hat{x}(s)),  \\
		\hat{x}(t) &= x.
	\end{aligned}
	\right.
\label{sys_closed}
\end{equation}
Define $\hat{v}(\cdot)$ by
\[
	\hat{v}(r) = \sigma (\hat{x}(r), \underline{u}(r,\hat{x}(r)))^T \nabla W(r,\hat{x}(r)),
	\ t \leq r \leq T.
\]
From  \eqref{fund_v} with $v(\cdot)=\hat{v}(\cdot)$,
\begin{multline}
	W(s,\hat{x}(s)) -\frac{1}{2}\int_t^s |\hat{v}(r)|^2 dr\\
	= W(t,x) +\int_t^s  \bigg\{ \frac{\partial W}{\partial r}(r,\hat{x}(r))
			+(f(\hat{x}(r), \alpha^{\underline{u}}[\hat{v}](r))
			+\sigma (\hat{x}(r),  \alpha^{\underline{u}}[\hat{v}](r))\hat{v}(r)) 
				\cdot \nabla W(r,\hat{x}(r))  \\
				 -\frac{1}{2} |\hat{v}(r)|^2 \bigg\} dr.
\label{fund_v_hat}
\end{multline}
By noting \eqref{argmax_H^u}, we have
\begin{multline}
	W(s,\hat{x}(s)) -\frac{1}{2}\int_t^s |\hat{v}(r)|^2 dr \\
	= W(t,x) +\int_t^s  \left\{ \frac{\partial W}{\partial r}(r,\hat{x}(r))
		+H^{\alpha^{\underline{u}}[\hat{v}](r)}\left(\hat{x}(r), \nabla W(r,\hat{x}(r)) \right)	\right\}dr.
\label{fund_H^u}
\end{multline}

We shall consider two cases:
Suppose that
\[
	\frac{\partial W}{\partial r}(r,\hat{x}(r)) 
	+H^{\alpha^{\underline{u}}[\hat{v}](r)}(\hat{x}(r), \nabla W(r,\hat{x}(r))) \geq 0, 
	\ t <  r <T.
\]
Then, by \eqref{fund_H^u} with $s=T$,
\[
	W(T, \hat{x}(T)) -\frac{1}{2}\int_t^T |\hat{v}(r)|^2 dr
	\geq W(t,x).
\]
From \eqref{QVI_T},
\begin{align*}
	W(T,\hat{x}(T)) = \min_{u \in U}l(\hat{x}(T), u)
	&\leq l(\hat{x}(T), \alpha^{\underline{u}}[\hat{v}](T)) \\
	&\leq \sup_{t \leq s \leq T} l(\hat{x}(s),\alpha^{\underline{u}}[\hat{v}](s))  
	=\int_{[t,T]}^{\oplus} l(\hat{x}(s), \alpha^{\underline{u}}[\hat{v}](s)) ds.
\end{align*}
Here note that $\esup$ on $[t,T]$ coincides with $\sup$ on $[t,T]$
because $ l(\hat{x}(s),\alpha^{\underline{u}}[\hat{v}](s))$
is continuous on $[t,T]$.
Thus we have
\[
	W(t,x) \leq
	\int_{[t,T]}^{\oplus} l(\hat{x}(s), \alpha^{\underline{u}}[\hat{v}](s)) ds
	-\frac{1}{2}\int_t^T |\hat{v}(s)|^2 ds
	\leq J(t,x; \alpha^{\underline{u}}).
\]

We next consider the case where there exists $r_0 \in (t,T)$ such that
\[
	\frac{\partial W}{\partial r}(r_0,\hat{x}(r_0)) 
	+H^{\alpha^{\underline{u}}[\hat{v}](r_0)}(\hat{x}(r_0), \nabla W(r_0,\hat{x}(r_0))) < 0.
\]
Define $\tau$ by
\[
	\tau \equiv \inf \left\{ r \in [t, T] \ ; \
	\frac{\partial W}{\partial r}(r,\hat{x}(r)) 
	+H^{\alpha^{\underline{u}}[\hat{v}](r)}(\hat{x}(r), \nabla W(r,\hat{x}(r))) < 0 \right\}
\]
We can show that the following claims hold:
\begin{gather}
	 \frac{\partial W}{\partial r}(r,\hat{x}(r)) 
	+H^{\alpha^{\underline{u}}[\hat{v}](r)}(\hat{x}(r), \nabla W(r,\hat{x}(r))) 
	\geq 0, \ t \leq  r < \tau, 
	\label{before_tau} \\
	 l(\hat{x}(\tau), \alpha[\hat{v}](\tau)) \geq W(\tau, \hat{x}(\tau)). \label{at_tau}
\end{gather}
\eqref{before_tau} is obvious from the definition of $\tau$.
For the proof of \eqref{at_tau}, since $W(t,x)$ is a solution of \eqref{QVI},
\[
	\max\left\{ \frac{\partial W}{\partial r}(r,x)+H^u(x,\nabla W(r,x)), l(x,u)-W(r,x) \right\} \geq 0,
	\ \forall u \in U.
\]
In particular,
\begin{multline}
	\max\left\{ \frac{\partial W}{\partial r}(r,\hat{x}(r))
	+H^{\alpha^{\underline{u}}[\hat{v}](r)}(\hat{x}(r),\nabla W(r,\hat{x}(r))), 
		l(\hat{x}(r),\alpha^{\underline{u}}[\hat{v}](r))-W(r,\hat{x}(r)) \right\} \geq 0, \\
	\ t < r <T.
\label{QVI_u}
\end{multline}
Take a sequence $\{ r_n \}$ such that $r_n \downarrow \tau \ (n \uparrow \infty)$ and
\[
	\frac{\partial W}{\partial r}(r_n,\hat{x}(r_n))
	+H^{\alpha^{\underline{u}}[\hat{v}](r_n)}(\hat{x}(r_n),\nabla W(r_n,\hat{x}(r_n))) < 0, 
	\ \forall n.
\]
Thus, from \eqref{QVI_u},
\[
	l(\hat{x}(r_n), \alpha^{\underline{u}}[\hat{v}](r_n)) \geq W(r_n, \hat{x}(r_n)), \ \forall n.
\]
Since $s \mapsto l(\hat{x}(s), \alpha^{\underline{u}}[\hat{v}](s))
= l(\hat{x}(s), \underline{u}(s,\hat{x}(s)))$ is continuous, 
we have  by taking $n \to \infty$
\[
	l(\hat{x}(\tau),\alpha^{\underline{u}}[\hat{v}](\tau)) \geq W(\tau, \hat{x}(\tau)).
\]

In \eqref{fund_H^u}, if we take $s=\tau$,
\begin{multline*}
	W(\tau,\hat{x}(\tau))-\frac{1}{2}\int_t^\tau |\hat{v}(r)|^2 dr \\
	= W(t,x) +\int_t^\tau  \left\{ \frac{\partial W}{\partial r}(r,\hat{x}(r)) 
	+H^{\alpha^{\underline{u}}[\hat{v}](r)}(\hat{x}(r), \nabla W(r,\hat{x}(r))) \right\} dr.
\end{multline*}
By \eqref{before_tau} and \eqref{at_tau},
\[
	l(x(\tau), \alpha^{\underline{u}}[\hat{v}](\tau))
		-\frac{1}{2}\int_t^\tau |\hat{v}(r)|^2 dr \\
	\geq W(t,x) .
\]
If we make the dependence on $\alpha^{\underline{u}}[\hat{v}]$ and $\hat{v}$ clear, 
this inequality can be written by
\begin{equation}
	l(x^{\alpha^{\underline{u}}[\hat{v}], \hat{v}}(\tau), \alpha^{\underline{u}}[\hat{v}](\tau))
	-\frac{1}{2}\int_t^\tau |\hat{v}(r)|^2 dr \\
	\geq W(t,x)  
\label{case2_fund}
\end{equation}
where $x^{\alpha^{\underline{u}}[\hat{v}],\hat{v}}$ is the solution of \eqref{mp-SDE}
corresponding to $\alpha^{\underline{u}}[\hat{v}]$, $\hat{v}$.

Define $\tilde{v}: [t,T] \to \mathbb{R}^n$ by
\[
	\tilde{v}(s) =
	\begin{cases}
		\hat{v}(s), & t \leq s \leq \tau, \\
		0, & \tau < s \leq T.	
	\end{cases}
\]
Then, from the definition of $\alpha^{\underline{u}}$, we can see that
\[
	x^{\alpha^{\underline{u}}[\tilde{v}], \tilde{v}}(r)
	= x^{\alpha^{\underline{u}}[\hat{v}], \hat{v}}(r),
	\ t \leq  r \leq \tau.	
\]
Therefore, from \eqref{case2_fund},
\begin{align*}
	 W(t,x) &\leq
	l(x^{\alpha^{\underline{u}}[\tilde{v}], \tilde{v}}(\tau), \alpha^{\underline{u}}[\tilde{v}](\tau))
	-\frac{1}{2}\int_t^T |\tilde{v}(r)|^2 dr	\\
	&\leq \int_{[t,T]}^{\oplus}
	 l(x^{{\alpha}^{\underline{u}}[\tilde{v}], \tilde{v}}(s), \alpha^{\underline{u}}[\tilde{v}](s)) ds
	-\frac{1}{2}\int_t^T |\tilde{v}(r)|^2 dr	
	\leq {J}(t,x; \alpha^{\underline{u}}).
\end{align*}

We shall now prove (b).
By \eqref{arg_min}, we have for $(s,y) \in (t,T) \times \mathbb{R}^n$
\[
	\max \left\{ \frac{\partial W}{\partial s}(s,y) +H^{\underline{u}^\ast(s,y)} (y, \nabla W(s,y)),
	l(y,\underline{u}^\ast(s,y))-W(s,y) \right\} = 0.
\]
This implies that
\begin{equation}
	\frac{\partial W}{\partial s}(s,y) +H^{\underline{u}^\ast(s,y)} (y, \nabla W(s,y)) \leq 0
	\text{ and }
	l(y,\underline{u}^\ast(s,y))-W(s,y) \leq 0.
\label{arg_min_ineq}
\end{equation}
For any $v \in L^2[t,T]$, consider \eqref{sys_fdbk} with $\underline{u}(s,y)=\underline{u}^\ast(s,y)$.
By \eqref{fund_v} with $\underline{u}(s,y)=\underline{u}^\ast (s,y)$,
\begin{align*}
	&W(s,x(s)) -\frac{1}{2}\int_t^s |v(r)|^2 dr\\
	&= W(t,x) +\int_t^s  \bigg\{ \frac{\partial W}{\partial r}(r,x(r))
			+(f(x(r), \alpha^{\underline{u}^\ast}[v](r))
			+\sigma (x(r),  \alpha^{\underline{u}^\ast}[v](r))v(r)) \cdot \nabla W(r,x(r))  \\
	& \qquad \qquad \qquad \qquad \qquad \qquad \qquad
			-\frac{1}{2} |v(r)|^2 \bigg\} dr \\
	&\leq W(t,x)+\int_t^s  \left\{  \frac{\partial W}{\partial r}(r,x(r))
			+H^{\alpha^{\underline{u}^\ast}[v](r)}(x(r), \nabla W(r,x(r)))  \right\} dr.	
\end{align*}
Thus, from \eqref{arg_min_ineq}, 
\[
	l(x(s), \alpha^{\underline{u}^\ast}[v](s)) -\frac{1}{2}\int_t^s |v(r)|^2 dr
	\leq  W(t,x), \ t \leq s \leq T.
\]
Hence we have
\[
	\int_{[t,T]}^\oplus l(x(s), \alpha^{\underline{u}^\ast}[v](s))ds -\frac{1}{2}\int_t^T |v(r)|^2 dr
	\leq  W(t,x).
\]
Since $v \in L^2[t,T]$ is taken arbitrarily, 
$$
	J(t,x;\alpha^{\underline{u}^\ast})
	=\sup_{v \in L^2[t,T]}
	\left\{ \int_{[t,T]}^\oplus l(x(s),\alpha^{\underline{u}^\ast}[v](s))ds
		-\frac{1}{2}\int_t^T |v(s)|^2 ds \right\} \leq W(t,x).
\eqno{\qed}
$$

\begin{rem}
(i) Theorem \ref{ver_thm} holds under weaker assumptions than (A1)--(A3).
An inspection of the proof shows that Theorem \ref{ver_thm} is true if we omit (A1)
and also the assumption that $f$ is bounded in (A2).
Instead of (A3), $l$ can be any continuous function.

\noindent
(ii) Theorem \ref{ver_thm} (a) can be extended to a general class including 
Markov control policies.
If we consider any strategy $\alpha$ satisfying (S1), (S2) (Section 3)
and $\alpha[v](T-)=\alpha[v](T)$, we can have $W(t,x) \leq J(t,x; \alpha)$
under more regularity assumptions on $W(t,x)$.
This can be shown with modifications of the proof of Theorem \ref{ver_thm}
by using Lemma \ref{lem_approx_str} for a particular Elliott-Kalton strategy
$\beta: L^\infty([t,T];U) \to L^2[t,T]$;
\[
	\beta[u](s)=\sigma(x(s),u(s))^T \nabla W(s,x(s)), \ t \leq s \leq T,
	\ u \in L^\infty([t,T];U),
\]
where $x(s)$ is the solution of \eqref{system} with $v(s)=\beta[u](s)$.
If we use $\alpha[v_\eps]$ and $v_\eps$ (see Lemma \ref{lem_approx_str})
instead of $\underline{u}$ and $\hat{v}$, respectively,
the arguments in the proof are still valid by working with additional discussions
on approximations.
Such idea is also used in the proof of \eqref{l_est_gen}.

\end{rem}

\section{Value function and dynamic programming principle}
To characterize a value function as a solution of the DPE, 
we usually need  (i) the dynamic programming principle (DPP) for the value function,
(ii) the correct form of the infinitesimal generator of the semigroup associated with the DPP.
Under (i) and (ii), one could show that the value function is a solution
of the evolution equation in some sense
(cf. \cite[Chapter 2]{FS06}). 
In this section, we shall introduce a strategy class which will be related to \eqref{QVI}
and show the DPP.

As a strategy class for our max-plus control problem,
we consider $\Gamma (t,T) \subset \Gamma_{EK}(t,T)$ defined by
the set of $\alpha : L^{2}[t,T]  \to L^\infty([t,T];U)$
satisfying the following conditions:
\renewcommand{\labelenumi}{(S\arabic{enumi})}
\begin{enumerate}
\item For any $v \in L^2[t,T]$, $s \mapsto \alpha[v](s)$ is right-continuous
			with left limits on $[t,T]$.
\item Let $v$, $\tilde{v} \in L^2[t,T]$ and $t \leq s \leq T$.
If $v = \tilde{v}$ a.e.\,on $[t,s]$, then $\alpha[v]=\alpha[\tilde{v}]$ on $[t,s]$.
\end{enumerate}
\renewcommand{\labelenumi}{(\roman{enumi})}
Note that we require $\alpha[v](r)=\alpha[\tilde{v}](r)$ for \textit{all} $r \in [t,s]$ in (S2)
(compare with \eqref{EK-cond} for Elliott-Kalton strategy).
We always take $\Gamma(t,T)$ satisfying (S1) and (S2) for our strategy class
in the rest of the arguments.

We point out some properties on instantaneous time delay of $\Gamma(t,T)$.
The proofs are immediate from the definition.
\begin{lem} \label{inst_delay}
For $\alpha \in \Gamma(t,T)$, \textup{(i)} and \textup{(ii)} hold:

\noindent
\textup{(i)} $\alpha[v](t)$ does not depend on $v \in L^2[t,T]$.

\noindent
\textup{(ii)} Let $v$, $\tilde{v} \in L^2[t,T]$ and $t <s \leq T$.
If $v=\tilde{v}$ a.e.\,on $[t,s)$, then $\alpha[v](r)=\alpha[\tilde{v}](r)$ for all $ r \in [t,s]$.
\end{lem}
\begin{rem}
If $\underline{u}$ is any Lipschitz Markov control policy (Section 2)
and we consider $\Gamma(t,T)$ satisfying (S1) and (S2), then $\alpha^{\underline{u}} \in \Gamma (t,T)$. On the other hand, if $\alpha [v] (s) = \phi [v(s)]$ for some nonconstant Borel measurable function $\phi\colon\mathbb{R}^d \rightarrow U$, then $\alpha \in \Gamma_{EK} (t,T)$ but $\alpha \not\in \Gamma (t,T)$. 
For any $\alpha \in \Gamma(t,T)$ and $t<s<T$, $\alpha[v]$ has the left hand limit $\alpha[v](s-)$
which depends only on $v(r)$ for $r<s$. Moreover, $\alpha[v](s-)=\alpha[v](s)$ except for countably many $s$.
This expresses (in a mathematically imprecise way) the intuitive idea that condition (S1) allows $\alpha[v](s)$
to depend on past values $v(r)$ for $r<s$, but not on the current value $v(s)$.
\end{rem}

For $t<T$, let $V(t,x)$ be the value function associated with $\Gamma(t,T)$, \textit{i.e.},
\begin{equation}
\begin{aligned}
	V(t,x)
	&= \inf_{\alpha \in \Gamma(t,T)} E^{+}_{tx}
		\left[ \int_{[t,T]}^\oplus l(x(s),\alpha[v](s))ds \right] \\
	&=\inf_{\alpha \in \Gamma(t,T)} \sup_{v \in L^2[t,T]}
	\left\{ \int_{[t,T]}^\oplus l(x(s),\alpha[v](s))ds -\frac{1}{2}\int_t^T |v(s)|^2 ds \right\},
\end{aligned}
\label{value}
\end{equation}
where $x(s)$ is the solution of \eqref{mp-SDE}.
When $t=T$, we let $\Gamma (T,T)=U$ and define $V(T,x)$ by
\[
	V(T,x)= \min_{u \in U} l(x,u).
\]

As the following theorem shows, we have the DPP for $V(t,x)$.
\begin{thm}\label{thm_DPP}
For any $(t,x) \in [0,T] \times \mathbb{R}^n$ and $t \leq r \leq T$,
\begin{equation}
	V(t,x) = \inf_{\alpha \in \Gamma(t,r)}
		E^{+}_{tx} \left[ \int_{[t,r]}^\oplus l(x(s),\alpha[v](s)) ds \oplus V(r,x(r)) \right].
\label{DPP}
\end{equation}
\end{thm}
\noindent
\textit{Proof.} \ 
In the case where $t=r=T$ or $t=r<T$, \eqref{DPP} is immediate.

We next consider the case where $t<r<T$. Let $W(t,x)$ be the right-hand side (RHS) of \eqref{DPP}.
As in Section 2, we identify $v$ with $(v_1, v_2)$, where $v_1=v|_{[t,r]}$ and $v_2=v|_{[r,T]}$.
For any $\eps>0$, we take $\alpha^1 \in \Gamma(t,r)$ such that
\begin{equation}
	W(t,x) +\eps > E^{+}_{tx} \left[ \int_{[t,r]}^\oplus l(x(s),\alpha^1[v_1](s)) ds \oplus V(r,x(r)) \right],
\label{W-approx}
\end{equation}
where $x(s)$ $( t \leq s \leq r)$ is the solution of \eqref{system} on $[t,r]$ for $u=\alpha^1[v_1]$
and $v=v_1 \in L^2[t,r]$.
For each $\xi \in \mathbb{R}^n$, choose $\alpha^2_\xi \in \Gamma(r,T)$ such that
\begin{equation}
	V(r,\xi)+\eps >E_{r\xi}^+ \left[ \int_{[r,T]}^\oplus l(x^2(s), \alpha^2_\xi [v_2](s)) ds \right],
\label{V_r-approx}
\end{equation}
where $x^2(s)$ $(r \leq s \leq T)$ is the solution of \eqref{system}
on $[r,T]$ for $u =\alpha^2_\xi[v_2] $ and $v= v_2 \in L^2[r,T]$
with initial condition $x^2(r)=\xi$.

For $v \in L^2[t,T]$, define $\alpha[v] \in L^\infty([t,T];U)$ by
\begin{equation}
	\alpha[v](s)=
	\begin{cases}
		\alpha^1 [v_1](s),& t \leq s <r, \\
		\alpha^2_{x(r)} [v_2](s),& r \leq s \leq T.
	\end{cases}
\label{mk_str1}
\end{equation}
$x(r)$ is the solution of \eqref{system} given by $u=\alpha^1[v_1]$ and $v=v_1$.
Note that $\alpha_{x(r)}[v_2](r)$ does not depend on $v_2$
because of Lemma \ref{inst_delay} (i).
Then, 
it is not difficult to check that $\alpha$ satisfies (S1) and (S2), that is, 
$\alpha \in \Gamma (t,T)$.

By using \eqref{tower},
\begin{align*}
	&E^{+}_{tx}\left[ \int_{[t,T]}^\oplus l(x(s),\alpha[v](s)) ds \right]  \\
	&=E^{+}_{tx} \left[  \int_{[t,r]}^\oplus l(x(s),\alpha^1[v_1](s))ds \oplus 
			E^{+}_{tx} \left[ \left. \int_{[r,T]}^\oplus l(x(s),\alpha[v](s)) ds  \right| v_1 \right]  \right] \\
	&=E^{+}_{tx} \left[  \int_{[t,r]}^\oplus l(x(s),\alpha^1[v_1](s))ds \oplus 
			E^{+}_{r x(r)} \left[ \int_{[r,T]}^\oplus l(x^2(s),\alpha^2_{x(r)}[v_2] (s))ds \right]\right]
\end{align*}
From \eqref{V_r-approx}, the last term can be estimated by
\begin{align*}
	&E^{+}_{tx} \left[  \int_{[t,r]}^\oplus l(x(s),\alpha^1[v_1](s))ds \oplus 
			E^{+}_{r x(r)} \left[ \int_{[r,T]}^\oplus l(x^2(s),\alpha^2_{x(r)}[v_2] (s))ds \right]\right] \\
	&\leq E^{+}_{tx} \left[  \int_{[t,r]}^\oplus l(x(s),\alpha^1[v_1](s))ds \oplus (V(r,x(r)) +\eps) \right] \\
	&\leq E^{+}_{tx} \left[  \int_{[t,r]}^\oplus l(x(s),\alpha^1[v_1](s))ds \oplus V(r,x(r)) \right] +\eps.
\end{align*}
Thus we have from \eqref{W-approx}
\[
	E^{+}_{tx}\left[ \int_{[t,T]}^\oplus l(x(s),\alpha[v](s)) ds \right] 
	\leq W(t,x) + 2\eps.
\]
Since $\alpha \in \Gamma(t,T)$,  we obtain
\[
	V(t,x) \leq W(t,x)+2\eps.
\]
Sending $\eps$ to $0$,  we have
\[
	V(t,x) \leq W(t,x).
\]

To prove $W(t,x) \leq V(t,x)$, take any $\alpha \in \Gamma(t,T)$.
We define $\alpha_1:L^2[t,r] \to L^\infty([t,r];U)$ by
\begin{equation}
	\alpha_1[v_1](s) = \alpha[v_1 \cdot v_2^0](s), \ t \leq s < r,
			\ v_1 \in L^2[t,r],
\label{mk_str2}
\end{equation}
where $v_2^0 \in L^2(r,T]$ is a (dummy) disturbance on $(r,T]$ and 
$v_1 \cdot v_2^0$ is the concatenation of $v_1$ and $v_2^0$:
\[
	v_1 \cdot v_2^0 (s)=
	\begin{cases}
		v_1(s),& t \leq s \leq r, \\
		v_2^0(s), & r < s \leq T.
	\end{cases}
\]
Note that $\alpha_1$ does not depend on the choice of $v^0_2$ by (S2).
Since $\alpha \in \Gamma(t,T)$, $\alpha_1 \in \Gamma(t,r)$.

For a given $v_1 \in  L^2 [t,r]$,  
we define $\alpha_2^{v_1}: L^2[r,T] \to L^\infty([r,T];U)$ as follows:
\begin{equation}
	\alpha^{v_1}_2[v_2](s)
	=\alpha[v_1 \cdot v_2](s), \ r \leq s \leq T, \ v_2 \in L^2[r,T].
\label{mk_str3}
\end{equation}
We can see $\alpha^{v_1}_2 \in \Gamma (r,T)$
because $\alpha \in \Gamma(t,T)$.
Thus, by the definition of $W(t,x)$ and $V(r,x(r))$,
\begin{align}
	W(t,x)
	&\leq E^{+}_{tx}\left[ \int_{[t,r]}^\oplus l(x(s), \alpha_1[v_1](s)) ds \oplus V(r,x(r)) \right]  \notag\\
	&\leq E^{+}_{tx}\left[ \int_{[t,r]}^\oplus l(x(s), \alpha_1[v_1](s)) ds \oplus
			E_{r x(r)}^+\left[ \int_{[r, T]}^\oplus l(x_2(s), \alpha_2^{v_1}[v_2](s)) ds \right]\right] 
	\label{W-est-ab}
\end{align}
By using \eqref{tower}, 
\begin{align*}
	\text{RHS  of \eqref{W-est-ab}}
	&=E^{+}_{tx}\left[ \int_{[t,r]}^\oplus l(x(s), \alpha_1[v_1](s)) ds \oplus
			E^{+}_{tx} \left[  \left. \int_{[r, T]}^\oplus l(x(s), \alpha_2^{v_1}[v_2](s)) ds \right| v_1 
			\right] \right] \\
	&=E^{+}_{tx}\left[ \int_{[t,T]}^\oplus l(x(s),\alpha[v](s)) ds \right]. 
\end{align*}
Thus we have
\[
	W(t,x) \leq E^{+}_{tx}\left[ \int_{[t,T]}^\oplus l(x(s),\alpha[v](s)) ds \right].
\]
Since $\alpha \in \Gamma(t,T)$ is taken arbitrarily,
$$
	W(t,x) \leq V(t,x).
$$

Finally, let $t<r=T$. For any $\alpha \in \Gamma (t,T)$,
\[
	V(T,x(T)) \leq l(x(T),\alpha[v](T-)) \leq \sup_{t<s<T}  l (x(s),\alpha[v](s))
	=\int_{[t,T]}^\oplus l(x(s),\alpha[v](s))ds.
\]
In this case, \eqref{DPP} is immediate from \eqref{g-value}. $\qed$

By using (A1)--(A3) and the DPP for the value function,
we can obtain the following regularity result on the value function.
\begin{prop} \label{prop_reg}
$V(t,x)$ is bounded Lipschitz continuous on $[0,T] \times \mathbb{R}^n$.
\end{prop}
\noindent
\textit{Proof.} \
It is obvious that $V(t,x)$ is bounded.
Since $f$, $\sigma$, $l$ are time-independent, if $t<T$ then $V(t,x)$ can be rewritten as follows:
\begin{align}
	V(t,x)
	&=\inf_{\alpha \in \Gamma(0,T-t)} E^{+}_{0x}
		\left[ \int_{[0,T-t]}^\oplus l(x(s),\alpha[v](s)) ds \right] \notag \\
	&=\inf_{\alpha \in \Gamma(0,T-t)} \sup_{v \in L^2[0,T-t]}
		\left\{ \int_{[0,T-t]}^\oplus l(x(s), \alpha[v](s)) ds -\frac{1}{2}\int_{0}^{T-t} |v(s)|^2ds\right\},
	\label{shift_value}
\end{align}
where $x(s)$ is the solution of 
\begin{equation}
	\left\{
	\begin{aligned}
		\frac{dx}{ds}(s) &= f(x(s),\alpha[v](s))+\sigma(x(s),\alpha[v](s))v(s), \ s \geq 0, \\
		x(0) &= x.
	\end{aligned}	
	\right.
\label{shift_sys}
\end{equation}
Since $l$ is bounded, it is sufficient to consider the supremum of \eqref{shift_value} on 
\begin{equation}
	\| v \|_{L^2[0,T-t]} \leq K
\label{rest_v}
\end{equation}
 for some constant $K$
where $K$ does not depend on $t$, $x$, $\alpha$.

Let $\tilde{x}(s)$ be the solution of \eqref{shift_sys} with the initial condition $\tilde{x}(0)=y$
and set $\zeta (s)=x(s)-\tilde{x}(s)$.
Then, from (A2), we can find constants $C_1$, $C_2$ such that
\[
	|\zeta (s)| \leq |x-y| +\int_0^s C_1 |\zeta (r)| + C_2 |\zeta (r)| |v(r)| dr,
	\ 0 \leq s \leq T.
\]
By Grownwall's inequality,
\[
	|\zeta(s)| \leq |x-y|\left\{ 1+\int_0^s (C_1+C_2 |v(r)|)
		e^{\int_r^s (C_1+C_2|v(\tau)|)d\tau} dr\right\}, \ 0 \leq s \leq T.
\]
Since we only consider $v$ satisfying \eqref{rest_v},
we can find $C_K>0$ such that
\[
	|\zeta(s)| \leq C_K |x-y|, \ 0 \leq s \leq T.
\]
Since $l_x$ is bounded,
\[
	l(x(s),\alpha[v](s)) \leq l(\tilde{x}(s), \alpha[v](s)) + \| l_x \|_{\infty} |\zeta(s)|.	
\]
Thus we have
\[
	\int_{[0,T-t]}^\oplus 	l(x(s),\alpha[v](s)) ds
	\leq \int_{[0,T-t]}^\oplus 	l(\tilde{x}(s),\alpha[v](s)) ds + L_K |x-y|,
\]
where $L_K=\| l_x \|_{\infty} C_K$. Taking the max-plus expectation,
\[
	E^{+}_{0x} \left[ 	\int_{[0,T-t]}^\oplus 	l(x(s),\alpha[v](s)) ds\right]
	\leq 	E^{+}_{0y} \left[ 	\int_{[0,T-t]}^\oplus 	l(\tilde{x}(s),\alpha[v](s)) ds\right]
	+L_K |x-y|.
\]
Since $\alpha \in \Gamma(0,T-t)$ is taken arbitrarily,
\[
	V(t,x) \leq V(t,y) +L_K |x-y|.
\]
Therefore we have
\begin{equation}
	|V(t,x)-V(t,y)| \leq L_K|x-y|, \ 0 \leq t < T, \ x,y \in \mathbb{R}^n.
\label{Lip_x}
\end{equation}

We now show that $V(\cdot,x)$ is uniformly Lipschitz.
Let $0 \leq t<r \leq T$. From \eqref{shift_value}, it is easy to see that
\[
	V(r,x) \leq V(t,x).
\]
So it suffices to show for some $M>0$ that
\[
	V(t,x) \leq V(r,x) + M(r-t).
\]
By the DPP for \eqref{shift_value} and max-plus linearity of the max-plus expectation,
\begin{align}
	V(t,x)&=\inf_{\alpha \in \Gamma(0,r-t)} E^{+}_{0x}\left[ \int_{[0,r-t]}^\oplus l(x(s),\alpha[v](s))ds
		\oplus V(r, x(r-t)) \right]  \notag \\
	&=\inf_{\alpha \in \Gamma(0,r-t)}
		E^{+}_{0x}\left[ \int_{[0,r-t]}^\oplus l(x(s),\alpha[v](s))ds \right]
		\oplus E^{+}_{0x}  \left[ V(r, x(r-t)) \right].
	\label{shift_DPP2}
\end{align}
By (A2), there exists $C>0$ such that for any $v \in L^2[0,r-t]$,
\[
	|x(s) -x| \leq Cs + C \int_0^{s}|v(r)|dr, \ 0 \leq s \leq r-t.
\]
Since $x \mapsto V(t,x)$ is uniformly Lipschitz by \eqref{Lip_x},
\[
	V(r,x(r-t)) \leq V(r,x) + \tilde{C}_K (r-t)+\tilde{C}_K  \int_0^{r-t}|v(s)| ds,
\]
where $\tilde{C}_K=L_K C$.
Thus, we can estimate the second expectation of \eqref{shift_DPP2} by
\begin{align}
	& E^{+}_{0x}  \left[ V(r, x(r-t)) \right]  \notag \\
	&\leq \sup_{v \in L^2[0,r-t]}
	 	\left\{ V(r,x) + \tilde{C}_K (r-t)+\tilde{C}_K  \int_0^{r-t}|v(s)| ds
		 -\frac{1}{2}\int_0^{r-t}|v(s)|^2 ds \right\} \notag \\
	&\leq V(r,x)+\tilde{C}_K (r-t) +\frac{1}{2}\tilde{C}_K^2 (r-t)
	 = V(r,x) +{M}_K^1 (r-t), \label{est_at_r-t}
\end{align}
where ${M}_K^1= \tilde{C}_K+(1/2)\tilde{C}_K^2$.
Therefore we have from \eqref{shift_DPP2}
\begin{align}
	V(t,x)
	&\leq \inf_{\alpha \in \Gamma(0,r-t)}
		E^{+}_{0x}\left[ \int_{[0,r-t]}^\oplus l(x(s),\alpha[v](s))ds \right]
		\oplus ( V(r,x) +{M}_K^1 (r-t))	  \notag \\
	&=  ( V(r,x) +{M}_K^1 (r-t)) \oplus \inf_{\alpha \in \Gamma(0,r-t)}
		E^{+}_{0x}\left[ \int_{[0,r-t]}^\oplus l(x(s),\alpha[v](s))ds \right].
	\label{r-t_est}
\end{align}

Let us consider the constant strategy $\alpha[v](s)\equiv u_0 $,
\[
	\inf_{\alpha \in \Gamma(0,t-r)}
		E^{+}_{0x}\left[ \int_{[0,r-t]}^\oplus l(x(s),\alpha[v](s))ds \right]
	\leq E^{+}_{0x} \left[ \int_{[0,t-r]}^\oplus l(x^0(s), u_0) ds \right],
\]
where $x^0(s)$ is the solution of \eqref{shift_sys} for $\alpha[v](s)\equiv u_0$.
Since $l_x$ is bounded, by the same argument as  \eqref{est_at_r-t},
there exists $M^2_K>0$ such that
\[
	 E^{+}_{0x} \left[ \int_{[0,t-r]}^\oplus l(x^0(s), u_0) ds \right]
	 \leq l(x,u_0) + M^2_K (r-t).
\]
If we take $u_0 \in \textup{arg} \min_{u \in U} l(x,u)$,
\[
	l(x,u_0) =\min_{u \in U} l(x,u) =V(T,x) \leq V(r,x).
\]
Thus, we have
\[
	\inf_{\alpha \in \Gamma(0,t-r)}
		E^{+}_{0x}\left[ \int_{[0,r-t]}^\oplus l(x(s),\alpha[v](s))ds \right]
	\leq V(r,x)+ M^2_K (r-t).
\]
Hence, if we take $M_K = \max \{ M^1_K, M^2_K\}$, we obtain from \eqref{r-t_est}
$$
	V(t,x) \leq V(r,x)+M_K (r-t).
\eqno{\qed}
$$

\section{Viscosity approach to value function}
\subsection{Nonlinear parabolic equation equivalent to QVI}
As mentioned in Section 2, we want to show that the DPE for \eqref{value} is \eqref{QVI}.
Equation \eqref{QVI} seems quite reasonable because it is derived
from the generator with constant control case.
One difficulty on \eqref{QVI} is very nonlinear, for instance, 
it is not linear on $\partial V/ \partial t$.

On the other hand, by using the idea of \cite{BI89},
it is shown that \eqref{QVI} is equivalent to the following equation
(see Proposition \ref{QVI-u_Isaacs}):
\begin{equation}
	\frac{\partial V}{\partial t}+\mathcal{H}(x,V(t,x), \nabla V(t,x))=0,
	\ (t,x) \in (0,T) \times \mathbb{R}^n,
\label{u_Isaacs}
\end{equation}
where for $x \in \mathbb{R}^n$, $r \in \mathbb{R}$, $p \in \mathbb{R}^n$,
\begin{gather}
	\mathcal{H}(x,r,p)=\min_{u \in A(x,r)} H^u(x,p)
	=\min_{u \in A(x,r)} \max_{v \in \mathbb{R}^d}
		\left\{ (f(x,u)+\sigma(x,u)v) \cdot p -\frac{1}{2}|v|^2 \right\},
	\label{gen}\\
	A(x,r) =\{ u \in U \,;\, l(x,u) \leq r \}. \notag
\end{gather}
We define $\mathcal{H}(x,r,p)=\infty$ if $A(x,r) = \emptyset$.
Equation \eqref{u_Isaacs} looks like a standard nonlinear equation of parabolic type backward in time.
But we need to be careful in \eqref{u_Isaacs} because the Hamiltonian $\mathcal{H}(x,r,p)$
is discontinuous and it can be $\infty$.
As we will see later,
since \eqref{u_Isaacs} will be naturally derived by calculating the generator,  
we  mainly discuss \eqref{u_Isaacs} instead of \eqref{QVI}. 
But one should notice that \eqref{QVI} is useful to prove Theorem \ref{ver_thm}.
See also Section 7.
We first show that \eqref{QVI} and \eqref{u_Isaacs} are equivalent in viscosity sense.

Let us recall the definition of viscosity solutions
with a discontinuous Hamiltonian. 
For the later argument, we shall define discontinuous viscosity solutions.
Let $h(\xi)$ be a function on a subset in  a Euclidean space.
We denote by $h^\ast(\xi) $ and $h_\ast(\xi)$
the upper and the lower semi-continuous envelopes of $h$,
respectively:
\[
	h^\ast(\xi) =\limsup_{ \eta \to \xi}h(\eta), \
	h_{\ast}(\xi)=\liminf_{ \eta \to \xi}h(\eta).
\]
\begin{defi} \label{visc_def}
Let $W(t,x)$ be a locally bounded function on $(0,T) \times \mathbb{R}^n$.
$W(t,x)$ is a \textit{viscosity subsolution}
(\textit{resp.\,viscosity supersolution}) of \eqref{u_Isaacs} if the following holds:
if $(\hat{t},\hat{x}) \in (0,T) \times \mathbb{R}^n$ is a maximum point
(\textit{resp.\,}minimum point)
of $W^\ast-\varphi$  (\textit{resp.\,} $W_\ast-\varphi$)
for a $C^1$-function $\varphi(t,x)$ in $(0,T) \times \mathbb{R}^n$,
then
\begin{gather*}
	\frac{\partial \varphi}{\partial t}(\hat{t},\hat{x})
	+\mathcal{H}^\ast (\hat{x}, W^\ast(\hat{t},\hat{x}), \nabla \varphi (\hat{t},\hat{x})) \geq 0,  \\
	\left( \text{ resp. } 
	\frac{\partial \varphi}{\partial t}(\hat{t},\hat{x})
	+\mathcal{H}_\ast (\hat{x}, W_\ast(\hat{t},\hat{x}), \nabla \varphi (\hat{t},\hat{x}))
		\leq 0 \right).
\end{gather*}
If $W(t,x)$ is a viscosity sub and supersolution, $W(t,x)$ is called a \textit{viscosity solution}.
\end{defi}


We will give some properties on $\mathcal{H}(x,r,p)$ used in \cite{BI89}.
They can be proved in a similar way to \cite{BI89}, so we omit the proofs.
\begin{lem}[cf.\,Lemma 2.3, \cite{BI89}] \label{est_set}
\textup{(i)}  If $r < r'$, then $A(x,r) \subset A(x,r')$.  \\
\textup{(ii)}  Let $r<r'$, $x \in \mathbb{R}^n$. Then there exists $\delta=\delta(|r'-r|)>0$ such that
\[
	A(y,r) \subset A(x,r'), \ \forall y \in B_\delta (x),
\]
where $B_\delta(x)$ is the open ball centered at $x$ with radius $\delta$.
\end{lem}

\begin{lem}[cf.\,Lemma 2.4, \cite{BI89}] \label{H-err}
\textup{(i)} If $r<r'$, then $\mathcal{H}(x,r,p) \geq \mathcal{H}(x,r',p)$. \\
\textup{(ii)} Let $r<r'$, $x \in \mathbb{R}^n$. Then there exists $\delta=\delta(|r'-r|)>0$ such that
\[
	\mathcal{H}(y,r,p) \geq \mathcal{H}(x,r',p) -L(|p|+1)|p||x-y|, \
	\forall y \in B_\delta(x), \ \forall p \in \mathbb{R}^n
\]
for some constant $L>0$.
\end{lem}

\begin{lem}[cf.\,Proposition 2.5, \cite{BI89}] \label{H-semi-lim}
For $x \in \mathbb{R}^n$, $r \in \mathbb{R}$, $p \in \mathbb{R}^n$, \\
\textup{(i)} $\mathcal{H}_\ast(x,r,p)=\mathcal{H}(x,r+0,p)$, \\
\textup{(ii)} $\mathcal{H}^\ast(x,r,p)=\mathcal{H}(x,r-0,p)$.
\end{lem}

In Lemma \ref{H-semi-lim}, we can identify the semi-continuous envelopes as follows.
\begin{lem} \label{mono_lim}
\textup{(i)} $\mathcal{H}_\ast(x,r,p) =\mathcal{H}(x,r,p)$. \\
\textup{(ii)} $\mathcal{H}^\ast(x,r,p)=\inf_{u \in A'(x,r)}H^u(x,p)$
where $A'(x,r)\equiv \{ u \in U \, ; \, l(x,u)<r \}$.
\end{lem}
\noindent
\textit{Proof.} \
We shall prove (i).
Since $s \mapsto \mathcal{H}(x,s,p)$ is non-increasing,
\[
	\mathcal{H}(x,r+1/n,p) \leq \mathcal{H}(x,r,p).
\]
Taking $n \to \infty$, 
\begin{equation}
	\mathcal{H}(x,r+0,p) \leq \mathcal{H}(x,r,p).
\label{ineq_l1}
\end{equation}

To show the equality holds in \eqref{ineq_l1},  suppose that
\begin{equation}
		\mathcal{H}(x,r+0,p) <\mathcal{H}(x,r,p).
\label{l_contr}
\end{equation}
By noting that $s \mapsto \mathcal{H}(x,s,p)$ is non-increasing, we have
\begin{equation}
	\min_{u \in A(x,r+1/n)}H^u(x,p) = \mathcal{H}(x,r+1/n,p)
	\leq \mathcal{H}(x,r+0,p), \ \forall n \in \mathbb{N}.
\label{ineq_l2}
\end{equation}
Let $u_n \in A(x,r+1/n)$ be a minimizer of $\min_{u \in A(x,r+1/n)}H^u(x,p)$, \textit{i.e.},
\begin{equation}
	H^{u_n}(x,p)=\min_{u \in A(x,r+1/n)}H^u(x,p)=\mathcal{H}(x,r+1/n,p).
\label{max_pt_l}
\end{equation}
Note that $A(x,r+1/n)\not= \emptyset$ 
because $\mathcal{H}(x,r+1/n,p)<\infty$ by \eqref{l_contr} and \eqref{ineq_l2}.
Since $U$ is compact, there exists a subsequence $\{ u_{n_j} \}$ and $\bar{u} \in U$ such that
\[
	u_{n_j} \to \bar{u} \ (j \to \infty).	
\]
By \eqref{ineq_l2} and \eqref{max_pt_l},
\[
	H^{u_{n_j}}(x,p) \leq \mathcal{H}(x,r+0,p).
\]
Thus, taking the limit as $j \to \infty$,
\begin{equation}
	H^{\bar{u}}(x,p) \leq \mathcal{H}(x,r+0,p).
\label{u_bar1}
\end{equation}
Since $u_{n_j} \in A(x,r+1/{n_j})$,
\[
	l(x,u_{n_j}) \leq r+1/n_j, \ j =1,2,\cdots.
\]
Taking the limit as $j \to \infty$, we have
\begin{equation}
	l(x,\bar{u}) \leq r, \textit{ i.e. },
	\bar{u} \in A(x,r).
\label{u_bar2}
\end{equation}
Therefore, we have from \eqref{u_bar1} and \eqref{u_bar2}
\[
	\mathcal{H}(x,r,p) =\min_{u \in A(x,r)}H^u(x,p) \leq H^{\bar{u}}(x,p)
	\leq \mathcal{H}(x,r+0,p).
\]
This contradicts to \eqref{l_contr}. Hence we obtain
\[
		\mathcal{H}(x,r+0,p) = \mathcal{H}(x,r,p).
\]
From Lemma \ref{H-semi-lim} (i), we have (i).

We shall next prove (ii). It is enough to show that
\[
	\mathcal{H}(x,r-0,p)=\inf_{u \in A'(x,r)} H^u(x,p).
\]
Since $A(x,r-1/n) \subset A'(x,r)$ $(n=1,2,\cdots)$,
\[
	\mathcal{H}(x,r-1/n,p)=\min_{u \in A(x,r-1/n)}H^u(x,p)
	\geq \inf_{u \in A'(x,r)}H^u(x,p).
\]
Sending $n$ to $\infty$,
\begin{equation}
	\mathcal{H}(x,r-0,p) \geq \inf_{u \in A'(x,r)}H^u(x,p).
\label{ineq_u1}
\end{equation}

Since $s \mapsto \mathcal{H}(x,s,p)$ is non-increasing,
\[
	\mathcal{H}(x,r-1/n,p) \geq \mathcal{H}(x,r-0,p), 
\]
\begin{equation}
\textit{ i.e. } \
	H^{u}(x,p) \geq \mathcal{H}(x,r-0,p), \ \forall u \in A(x,r-1/n).
\label{ineq_u2}
\end{equation}
Note that 
\[
	A'(x,r) = \bigcup_{n=1}^\infty A(x,r-1/n).
\]
So, we have from \eqref{ineq_u2}
\[
	H^u(x,p) \geq \mathcal{H}(x,r-0,p),\ \forall u \in A'(x,r).
\]
Thus, we have
$$
	 \inf_{u \in A'(x,r)}H^u(x,p) \geq \mathcal{H}(x,r-0,p).
\eqno{\qed}
$$

By using Lemma \ref{mono_lim}, we can prove the following result on the relation 
between \eqref{QVI} and \eqref{u_Isaacs}.
\begin{prop} \label{QVI-u_Isaacs}
Let $W(t,x)$ be a locally bounded function on $(0,T) \times \mathbb{R}^n$.
$W(t,x)$ is a viscosity subsolution (\textup{resp.\,}viscosity supersolution) of  \eqref{QVI} 
if and only if $W(t,x)$ is a viscosity subsolution (\textup{resp.\,}viscosity supersolution) 
 of \eqref{u_Isaacs}. 
 Here $W(t,x)$ is a viscosity subsolution (\textup{resp.}\,viscosity supersolution)
 of \eqref{QVI} if the following holds:
 if $(\hat{t},\hat{x}) \in (0,T) \times \mathbb{R}^n$ is a maximum  (\textup{resp.}\,minimum) point
 of $W^\ast(t,x)-\varphi(t,x)$
 (\textup{resp.}\,$W_\ast(t,x)-\varphi(t,x)$) for a $C^1$-function 
 $\varphi(t,x)$ on $(0,T) \times \mathbb{R}^n$,
 then
 \begin{gather*}
	\min_{u \in U} \max
	\left\{ \frac{\partial \varphi}{\partial t}(\hat{t},\hat{x})
	+H^u (\hat{x}, \nabla \varphi(\hat{t},\hat{x})), l(\hat{x},u)-W^\ast(\hat{t},\hat{x}) \right\} \geq 0 \\
	\left( \textup{resp.} \ 
	\min_{u \in U} \max
	\left\{ \frac{\partial \varphi}{\partial t}(\hat{t},\hat{x})
	+H^u (\hat{x}, \nabla \varphi(\hat{t},\hat{x})), l(\hat{x},u)-W_\ast(\hat{t},\hat{x}) \right\} \leq 0
	\right).
 \end{gather*}
\end{prop}
\noindent
\textit{Proof.} \
Note that we may assume $W^\ast(\hat{t}, \hat{x})=\varphi(\hat{t},\hat{x})$
in the definitions of subsolutions without loss of generality.
Then, it is sufficient for the equivalence of subsolutions to show that for a $C^1$-function $\varphi(t,x)$
and $(\hat{t},\hat{x}) \in (0,T)\times \mathbb{R}^n$, the following inequalities are equivalent:
\begin{gather}
	\min_{u \in U} \max \left\{ \frac{\partial \varphi}{\partial t}(\hat{t}, \hat{x})
		+H^u(\hat{x},\nabla \varphi(\hat{t},\hat{x})), l(\hat{x},u)-\varphi(\hat{t},\hat{x})
		\right \} \geq 0,
	\label{QVI-sub} \\
	\frac{\partial \varphi}{\partial t}(\hat{t}, \hat{x})
	+\mathcal{H}^\ast(\hat{x}, \varphi(\hat{t},\hat{x}), \nabla \varphi(\hat{t},\hat{x})) \geq 0.
	\label{EVL-sub}
\end{gather}

Let us assume \eqref{QVI-sub} holds.
Then, we have
\begin{equation}
	\max\left\{  \frac{\partial \varphi}{\partial t}(\hat{t}, \hat{x})
	+H^{u} (\hat{x},\nabla \varphi(\hat{t},\hat{x})),
	l(\hat{x},u) - \varphi(\hat{t},\hat{x}) \right\} \geq 0, \ \forall u \in U.
\label{QVI-sub-all}
\end{equation}
If $u \in A'(\hat{x},\varphi(\hat{t},\hat{x}))$, \textit{i.e.}, 
$l(\hat{x},u)-\varphi(\hat{t},\hat{x})<0$, \eqref{QVI-sub-all} implies that
\[
	\frac{\partial \varphi}{\partial t}(\hat{t}, \hat{x})
	+H^{u} (\hat{x},\nabla \varphi(\hat{t},\hat{x})) \geq 0.
\]
Therefore,
\[
	\frac{\partial \varphi}{\partial t}(\hat{t}, \hat{x})
	+H^{u} (\hat{x},\nabla \varphi(\hat{t},\hat{x})) \geq 0, \ \forall u \in A'(\hat{x},\varphi(\hat{t},\hat{x})).
\]
Thus, we have
\[
	\frac{\partial \varphi}{\partial t}(\hat{t}, \hat{x})
	+\inf_{u  \in  A'(\hat{x},\varphi(\hat{t},\hat{x}))}H^{u} (\hat{x},\nabla \varphi(\hat{t},\hat{x})) \geq 0.
\]
By Lemma \ref{mono_lim} (ii), we have
\[
	\frac{\partial \varphi}{\partial t}(\hat{t}, \hat{x})
	+\mathcal{H}^\ast(\hat{x}, \varphi(\hat{t},\hat{x}), \nabla \varphi(\hat{t},\hat{x})) \geq 0.
\]

We suppose \eqref{EVL-sub} holds.
From \eqref{EVL-sub} and Lemma \ref{mono_lim},
\begin{equation}
	\frac{\partial \varphi}{\partial t}(\hat{t}, \hat{x})
	+\inf_{u  \in  A'(\hat{x},\varphi(\hat{t},\hat{x}))}H^{u} (\hat{x},\nabla \varphi(\hat{t},\hat{x})) \geq 0.
\label{EVL-sub-expl}
\end{equation}
If $A'(\hat{x},\varphi(\hat{t},\hat{x}))=\emptyset$,
\[
	l(\hat{x},u) \geq \varphi(\hat{t},\hat{x}), \forall u \in U.
\]
Thus we have
\[
	\max \left\{\frac{\partial \varphi}{\partial t}(\hat{t}, \hat{x})
	+H^{u} (\hat{x},\nabla \varphi(\hat{t},\hat{x})),
	l(\hat{x},{u}) - \varphi(\hat{t},\hat{x}) \right\}
	\geq l(\hat{x},u) -\varphi(\hat{t},\hat{x}) \geq 0, \forall u \in U,
\]
which implies
\[
	\min_{u \in U} \max \left\{\frac{\partial \varphi}{\partial t}(\hat{t}, \hat{x})
	+H^{u} (\hat{x},\nabla \varphi(\hat{t},\hat{x})),
	l(\hat{x},{u}) - \varphi(\hat{t},\hat{x}) \right\} \geq 0.
\]

If  $A'(\hat{x},\varphi(\hat{t},\hat{x}))\not=\emptyset$,
we have from \eqref{EVL-sub-expl},
\begin{equation}
	\frac{\partial \varphi}{\partial t}(\hat{t}, \hat{x})
	+H^{u} (\hat{x},\nabla \varphi(\hat{t},\hat{x})) \geq 0, 
	\ \forall u  \in  A'(\hat{x},\varphi(\hat{t},\hat{x})).
\label{EVL-sub-all}
\end{equation}
Let $u \in U$ be arbitrarily taken. If $l(\hat{x}, u) < \varphi(\hat{t},\hat{x})$,
\textit{i.e.}, $u \in A'(\hat{x},\varphi(\hat{t},\hat{x}))$,
then \eqref{EVL-sub-all} implies
\[
	\max \left\{ \frac{\partial \varphi}{\partial t}(\hat{t}, \hat{x})
	+H^{u} (\hat{x},\nabla \varphi(\hat{t},\hat{x})),
	l(\hat{x}, u) - \varphi(\hat{t},\hat{x}) \right\}
	=\frac{\partial \varphi}{\partial t}(\hat{t}, \hat{x})
	+H^{u} (\hat{x},\nabla \varphi(\hat{t},\hat{x}))\geq 0.
\]
On the other hand, if $l(\hat{x}, u)  \geq \varphi(\hat{t},\hat{x})$,
\[
	\max\left\{ \frac{\partial \varphi}{\partial t}(\hat{t}, \hat{x})
	+H^{u} (\hat{x},\nabla \varphi(\hat{t},\hat{x})),
	l(\hat{x}, u) - \varphi(\hat{t},\hat{x}) \right\}
	\geq l(\hat{x}, u) - \varphi(\hat{t},\hat{x}) \geq 0.
\]
Thus,
\[
	\max\left\{ \frac{\partial \varphi}{\partial t}(\hat{t}, \hat{x})
	+H^{u} (\hat{x},\nabla \varphi(\hat{t},\hat{x})),
	l(\hat{x}, u) - \varphi(\hat{t},\hat{x}) \right\} \geq 0, \ \forall u \in U.
\]
Therefore, we have
\[
	\min_{u \in U} \max \left\{ \frac{\partial \varphi}{\partial t}(\hat{t}, \hat{x})
	+H^{u} (\hat{x},\nabla \varphi(\hat{t},\hat{x})),
	l(\hat{x}, u) - \varphi(\hat{t},\hat{x}) \right\} \geq 0.
\]

For the definitions of supersolutions, we may also suppose 
$W_\ast(\hat{t},\hat{x})=\varphi(\hat{t},\hat{x})$.
Thus, in the proof of the equivalence of supersolutions,
it suffices to check that the following inequalities are equivalent:
\begin{gather}
	\min_{u \in U} \max \left\{ \frac{\partial \varphi}{\partial t}(\hat{t}, \hat{x})
		+H^u(\hat{x},\nabla \varphi(\hat{t},\hat{x})), l(\hat{x},u)-\varphi(\hat{t},\hat{x})
		\right\} \leq 0,
	\label{QVI-super} \\
	\frac{\partial \varphi}{\partial t}(\hat{t}, \hat{x})
	+\mathcal{H}_\ast(\hat{x}, \varphi(\hat{t},\hat{x}), \nabla \varphi(\hat{t},\hat{x})) \leq 0.
	\label{EVL-super}
\end{gather}

Suppose \eqref{QVI-super} holds.
Since $U$ is compact,
there exists $\bar{u} \in U$ such that
\begin{align*}
	&\max \left\{ \frac{\partial \varphi}{\partial t}(\hat{t}, \hat{x})
		+H^{\bar{u}}(\hat{x},\nabla \varphi(\hat{t},\hat{x})), l(\hat{x},\bar{u})-\varphi(\hat{t},\hat{x}) \right\} \\
	&=\min_{u \in U} \max \left\{ \frac{\partial \varphi}{\partial t}(\hat{t}, \hat{x})
		+H^u(\hat{x},\nabla \varphi(\hat{t},\hat{x})), l(\hat{x},u)-\varphi(\hat{t},\hat{x})
		\right\} \leq 0.
\end{align*}
Thus,
\[
	\frac{\partial \varphi}{\partial t}(\hat{t}, \hat{x})
		+H^{\bar{u}}(\hat{x},\nabla \varphi(\hat{t},\hat{x})) \leq 0 \ \text{ and }
	\  l(\hat{x},\bar{u})-\varphi(\hat{t},\hat{x}) \leq 0.
\]
Since $l(\hat{x},\bar{u}) \leq \varphi(\hat{t},\hat{x})$,
we see  $\bar{u} \in A(\hat{x},\varphi(\hat{t},\hat{x}))$.
Therefore we have
\begin{align*}
	\frac{\partial \varphi}{\partial t}(\hat{t}, \hat{x})
	+\mathcal{H}(\hat{x},\varphi(\hat{t},\hat{x}), \nabla \varphi(\hat{t},\hat{x}))
	&=\frac{\partial \varphi}{\partial t}(\hat{t}, \hat{x})
	+\min_{u \in A(\hat{x},\varphi(\hat{t},\hat{x}))}H^u (\hat{x},\nabla \varphi(\hat{t},\hat{x}))  \\
	&\leq \frac{\partial \varphi}{\partial t}(\hat{t}, \hat{x})
		+H^{\bar{u}}(\hat{x},\nabla \varphi(\hat{t},\hat{x})) \leq 0.
\end{align*}
By Lemma  \ref{mono_lim} (i), we obtain
\[
	\frac{\partial \varphi}{\partial t}(\hat{t}, \hat{x})
	+\mathcal{H}_\ast(\hat{x},\varphi(\hat{t},\hat{x}), \nabla \varphi(\hat{t},\hat{x})) \leq 0.
\]

We assume \eqref{EVL-super} holds.
From \eqref{EVL-super} and Lemma \ref{mono_lim} (i),
\begin{equation}
	\frac{\partial \varphi}{\partial t}(\hat{t}, \hat{x})
	+\min_{u \in A(\hat{x},\varphi(\hat{t},\hat{x}))}H^u (\hat{x},\nabla \varphi(\hat{t},\hat{x})) 
	\leq 0.
\label{EVL-super-expl}
\end{equation}
Note that $A(\hat{x},\varphi(\hat{t},\hat{x})) \not= \emptyset$ because the minimum
of \eqref{EVL-super-expl} is finite.
Taking a minimum point $\bar{u} \in A(\hat{x},\varphi(\hat{t},\hat{x}))$ in \eqref{EVL-super-expl},
\[
	\frac{\partial \varphi}{\partial t}(\hat{t}, \hat{x})
	+H^{\bar{u}} (\hat{x},\nabla \varphi(\hat{t},\hat{x})) 
	\leq 0
\]
Combining  $\bar{u} \in A(\hat{x},\varphi(\hat{t},\hat{x}))$, \textit{i.e.},
$l(\hat{x},\bar{u}) - \varphi(\hat{t},\hat{x}) \leq 0$, 
\[
	\max \left\{ \frac{\partial \varphi}{\partial t}(\hat{t}, \hat{x})
	+H^{\bar{u}} (\hat{x},\nabla \varphi(\hat{t},\hat{x})),
	l(\hat{x},\bar{u}) - \varphi(\hat{t},\hat{x}) \right\} \leq 0.
\]
Thus, we obtain
$$
	\min_{u \in U} \max \left\{ \frac{\partial \varphi}{\partial t}(\hat{t}, \hat{x})
	+H^{u} (\hat{x},\nabla \varphi(\hat{t},\hat{x})),
	l(\hat{x},u) - \varphi(\hat{t},\hat{x}) \right\}  \leq 0.
\eqno{\qed}
$$

\subsection{Generators and evolution equations}
In optimal control problems, 
DPPs  can be often described by the semigroups associated with the value functions.
Once the control problem is related to the semigroup and the generator is identified,
one can  prove that the value function is a viscosity solution
of the evolution equation with the generator (cf. \cite{FS06}).
We shall discuss our max-plus control problem by the semigroup-generator
approach and prove that the value function  \eqref{value} is a unique viscosity solution of
\eqref{u_Isaacs} with \eqref{QVI_T}.

For $0 \leq t <r \leq T$, we define operator $F_{t,r}$ acting
on $\phi:\mathbb{R}^n \to \mathbb{R}^{-}$ by
\[
	F_{t,r}\phi (x) = \inf_{\alpha \in \Gamma(t,r)} E^{+}_{tx}
		\left[ \int_{[t,r]}^\oplus l(x(s),\alpha[v](s)) ds \oplus \phi(x(r))\right],  \ x \in \mathbb{R}^n.
\]
We understand that $F_{t,t}$ is an identity operator.
By using $F_{t,r}$, the value function $V(t,x)$ can be written as
\[
	V(t,x) = F_{t,T} \mathbf{0}(x),
\]
where $\mathbf{0}$ is a constant function taking its value $\mathbf{0}=-\infty$.
Note that $\mathbf{0}=-\infty$ is the additive identity in max-plus algebra.
Then, \eqref{DPP} can be written in terms of $F_{t,r}$ as follows:
\[
	F_{t,T}\mathbf{0} =F_{t,r} F_{r,T}\mathbf{0}, \ 0 \leq t<r \leq T.
\]
The DPP of \eqref{value} given by \eqref{DPP} corresponds to the semigroup property of
 $\{ F_{t,r} \}$.

Since $V(t,x)$ involves an inf-sup as seen in \eqref{value},
the form of the generator of $\{ F_{t,r} \}$ is not obvious.
Moreover, it is expected that the form depends on the choice of strategies classes (cf. \cite{KS05}).
By following the argument used in \cite{KS05}, we can show that 
our strategy class $\Gamma(t,T)$ satisfying (S1) and (S2)
is related to $\mathcal{H}(x,r,p)$ in \eqref{u_Isaacs}.
To state the result on the generators, 
we denote by $C^1_b ((0,T) \times \mathbb{R}^n)$ the set of bounded $C^1$-functions
with bounded first order derivatives.
\begin{thm} \label{cal_gen}
For any $\varphi \in C_b^1((0,T) \times \mathbb{R}^n)$
and $(t,x) \in (0,T) \times \mathbb{R}^n$, we have
\begin{gather}
	\limsup_{\delta \to 0+}\frac{1}{\delta} 
		\{ F_{t,t+\delta}\varphi(t+\delta,\cdot)(x) -\varphi(t,x) \}
	\leq \frac{\partial \varphi}{\partial t}(t,x)+\mathcal{H}^\ast(x, \varphi(t,x), \nabla \varphi(t,x)), 
	\label{u_est_gen}\\
	\frac{\partial \varphi}{\partial t}(t,x)+\mathcal{H}_\ast(x, \varphi(t,x), \nabla \varphi(t,x))
	\leq 	\liminf_{\delta \to 0+}\frac{1}{\delta} 
		\{ F_{t,t+\delta}\varphi(t+\delta,\cdot)(x) -\varphi(t,x)\}.
	\label{l_est_gen}
\end{gather}
\end{thm}

We need a result on approximations of strategies in $\Gamma(t,T)$.
The proof is given in the Appendix.
Let $\Delta_{EK}(t,T)$ be the set of Elliott-Kalton strategies
from $L^\infty([t,T];U)$ into $L^2[t,T]$.
\begin{lem} \label{lem_approx_str}
Let $\alpha \in {\Gamma}(t,T)$. For any $\beta \in \Delta_{EK}(t,T)$ and any $\eps>0$,
there exist $u^\eps \in L^\infty([t,T]; U)$ and $v^\eps \in L^2[t,T]$ such that\[
 \ |\alpha[v^\eps](s) -u^\eps(s)|<\eps, \ \forall s \in [t,T) \ \text{ and } \
	\beta[u^\eps](s) =v^\eps(s)\text{ a.e.}\,s \in [t,T].
\]
\end{lem}

\noindent
\textit{Proof of Theorem \textup{\ref{cal_gen}}.} \
We first prove \eqref{u_est_gen}.
In Lemma \ref{mono_lim} (ii) with $r=\varphi(t,x)$, $p=\nabla \varphi(t,x)$,
we consider $u \in U$ such that
\[
	l(x,u) < \varphi(t,x).
\]
It is enough to prove that
\begin{equation}
	\limsup_{\delta \to 0+}\frac{1}{\delta} 
		\{ F_{t,t+\delta}\varphi(t+\delta,\cdot)(x) -\varphi(t,x) \}
	\leq
	\frac{\partial \varphi}{\partial t}(t,x)
	+H^u (x, \nabla \varphi(t,x)).
\label{u_est_gen_red}
\end{equation}
We take $\rho_0>0$ such that
\begin{equation}
	l(x,u)<\varphi(t,x)-\rho_0.
\label{err_l-phi}
\end{equation}
If we take a constant strategy $\alpha[v](s) \equiv u$ for $F_{t,t+\delta}\varphi(t+\delta,\cdot)(x)$,
\begin{align}
	& F_{t,t+\delta}\varphi(t+\delta,\cdot)(x) \notag \\
	& \leq E^{+}_{tx}\left[ \int_{[t,t+\delta]}^\oplus l(x(s), u)ds \oplus\varphi(t+\delta, x(t+\delta)) \right]
	  \notag \\
	& = \sup_{v \in L^2[t,t+\delta]}
	\left\{ \int_{[t,t+\delta]}^\oplus l(x(s), u)ds\oplus \varphi(t+\delta,  x(t+\delta)) 
		 -\frac{1}{2}\int_{t}^{t+\delta} |v(s)|^2 ds \right\}.
	\label{op_const}
\end{align}
Since $l$ and $\varphi$ are bounded,
the supremum in \eqref{op_const} can be restricted to those $v$ satisfying
\begin{equation}
	\int_{t}^{t+\delta}|v(s)|^2 ds \leq M,
\label{rest_dist_const}
\end{equation}
for some $M>0$. Here $M$ does not depend on $u$ and $\delta$.
Then, we can show that there exists $c(M)>0$ such that for any $v$ satisfying \eqref{rest_dist_const},
\begin{gather}
	|x(s) -x| < c(M)\delta^{1/2}, \ t \leq s\leq t+\delta,  \label{unif_conti_ini}\\
	|\varphi(t+\delta,x(t+\delta)) - \varphi(t,x)| < c(M)\delta^{1/2}, \notag \\
	| l(x(s),u)-l(x,u)| < c(M) \delta^{1/2}, \ t \leq s \leq t +\delta. \notag
\end{gather}
For small $\delta>0$, we have
\begin{gather*}
	|\varphi(t+\delta,x(t+\delta)) -\varphi(t,x)|<\frac{1}{2}\rho_0, \\
	|l(x(s),u)-l(x,u)| < \frac{1}{2}\rho_0, \ t \leq s \leq t+\delta.
\end{gather*}
Then, by using \eqref{err_l-phi}, we have for $s \in [t,t+\delta]$
\[
	l(x(s),u) <l(x,u)+\frac{1}{2}\rho_0 <\varphi(t,x)-\frac{1}{2}\rho_0
	<\varphi(t+\delta,x(t+\delta)),
\]
which implies
\[
	\int_{[t,t+\delta]}^\oplus l(x(s),u)ds < \varphi(t+\delta,x(t+\delta))
\]
for small $\delta>0$.
Thus, from \eqref{op_const} with \eqref{rest_dist_const},
\[
	F_{t,t+\delta}\varphi(t+\delta,\cdot)(x) -\varphi(t,x)
	\leq
	\sup_{ \| v \|_{L^2[t,t+\delta]}^2 \leq M}
	\left\{ \varphi(t+\delta,x(t+\delta)) - \varphi(t,x)- \frac{1}{2}\int_t^{t+\delta}|v(s)|^2 ds \right\}
\]
for small  $\delta>0$.
By using the fundamental theorem of calculus and \eqref{unif_conti_ini},
\begin{align*}
	&\varphi(t+\delta,x(t+\delta)) - \varphi(t,x)- \frac{1}{2}\int_t^{t+\delta}|v(s)|^2 ds \\
	&=\int_t^{t+\delta} \frac{\partial \varphi}{\partial s}(s,x(s))
		+(f(x(s),u)+\sigma(x(s),u)v(s)) \cdot \nabla \varphi(s,x(s)) -\frac{1}{2}|v(s)|^2 ds \\
	&=\int_t^{t+\delta}
		\frac{\partial \varphi}{\partial t}(t,x)
		+(f(x,u)+\sigma(x,u)v(s)) \cdot \nabla \varphi(t,x) -\frac{1}{2}|v(s)|^2 ds +o(\delta) \\
	&\leq \left( \frac{\partial \varphi}{\partial t}(t,x)
		+\sup_{v \in \mathbb{R}^d}
		\left\{ (f(x,u)+\sigma(x,u)v) \cdot \nabla \varphi(t,x) -\frac{1}{2}|v|^2  \right\} \right) \delta
		+o(\delta)
\end{align*}
where $o(\delta)$ is uniform on $v$ satisfying \eqref{rest_dist_const}.
Therefore we obtain
\[
	F_{t,t+\delta}\varphi(t+\delta,\cdot)(x) -\varphi(t,x)
	\leq \left( \frac{\partial \varphi}{\partial t}(t,x) +H^u(x,\nabla \varphi(t,x)) \right) \delta +o(\delta).
\]
This implies \eqref{u_est_gen_red}. Hence we have proved \eqref{u_est_gen}.

To prove \eqref{l_est_gen}, we use the idea of the proof of \cite[Proposition 2.6]{KS05}.
If $u \in L^\infty([t,t+\delta];U)$ and $v \in L^2[t,t+\delta]$ are given,
let us define $F^{u,v}_{t,t+\delta}\phi(x)$ for $\phi: \mathbb{R}^n \to \mathbb{R}^{-}$ by
\[
	F^{u,v}_{t,t+\delta}\phi(x)
	=\int_{[t,t+\delta]}^\oplus l(x(s),u(s))ds\oplus \phi(x(t+\delta)) -\frac{1}{2}\int_{t}^{t+\delta}|v(s)|^2 ds,
\]
where $x(s)$ is the solution of \eqref{system}.

We consider a particular $\hat{\beta} \in \Delta_{EK}(t,t+\delta)$:
\[
	\hat{\beta}[u](s) = \sigma (\hat{x}(s), u(s))^{T} \nabla \varphi( s, \hat{x}(s)), \ t \leq s \leq t+\delta,
\]
where $\hat{x}(s)$ is the solution of \eqref{system} with $v(s)=\hat{\beta}[u](s)$

Fix $\rho>0$. We shall show that for small $\delta>0$
\begin{align}
	&F^{u,\hat{\beta}[u]}_{t,t+\delta}\varphi(t+\delta, \cdot)(x) - \varphi(t,x)  \notag \\
	&\geq \min \left\{ \left(
		\frac{\partial \varphi}{\partial t}(t,x)+\mathcal{H}(x, \varphi(t,x)+\rho, \nabla \varphi(t,x)) \right)\delta
		+o (\delta), 
	\frac{\rho}{2} -M\delta \right\},
\label{l_est_delta}
\end{align}
where $o(\delta)$ and $M$ are uniform on $u \in L^\infty([t,t+\delta];U)$.

For the first case, suppose $u(s) \in A(x, \varphi(t,x)+\rho)$ for a.e.\,$s \in [t,t+\delta]$.
From the definition of $\mathcal{H}(x,\varphi(t,x)+\rho,\nabla \varphi(t,x))$,
\begin{align*}
	&\mathcal{H}(x, \varphi(t,x)+\rho, \nabla \varphi(t,x))  \\
	&\leq \sup_{v \in \mathbb{R}^d} \left\{ (f(x,u(s)) +\sigma(x,u(s))v) \cdot \nabla \varphi(t,x)
		-\frac{1}{2}|v|^2 \right\} \\
	&=  \left(f(x,u(s)) +\sigma(x,u(s))\sigma(x,u(s))^T\nabla \varphi(t,x) \right) \cdot \nabla \varphi(t,x) \\
	&\qquad \qquad \qquad \qquad \qquad 
		-\frac{1}{2}|\sigma(x,u(s))^T\nabla \varphi(t,x)|^2, \ \text{a.e.\,}s \in [t,t+\delta].
\end{align*}
Integrating the above inequality on $[t,t+\delta]$, we have
\begin{align*}
	&\mathcal{H}(x, \varphi(t,x)+\rho, \nabla \varphi(t,x)) \delta \\
	&\leq \int_t^{t+\delta} \bigg\{
	 (f(x,u(s)) +\sigma(x,u(s))\sigma(x,u(s))^T\nabla \varphi(t,x)) \cdot \nabla \varphi(t,x) \\
	&\qquad \qquad \qquad \qquad
	-\frac{1}{2}|\sigma(x,u(s))^T\nabla \varphi(t,x)|^2	\bigg\} ds.
\end{align*}
Since $f$, $\sigma$, $\nabla \varphi$ are bounded and continuous, 
there exists $K>0$, which does not depend on $u(\cdot)$, such that
\begin{equation}
	|\hat{x}(s)-x| \leq K (s-t), \ t \leq s \leq t+\delta.
\label{err_sol_ini}
\end{equation}
Thus, we have
\begin{align}
	&\mathcal{H}(x, \varphi(t,x)+\rho, \varphi(t,x)) \delta \notag \\
	&\leq \int_t^{t+\delta}  \bigg\{ (f(\hat{x}(s),u(s)) +\sigma(\hat{x}(s),u(s))\hat{\beta}[u](s))
		\cdot \nabla \varphi(s,\hat{x}(s))  \notag \\
	& \qquad \qquad \qquad \quad \qquad \qquad  \quad \qquad \qquad 
	-\frac{1}{2}|\hat{\beta}[u](s)|^2	 \bigg\} ds + o(\delta).
	\label{u_est_delta}
\end{align}
where $o(\delta)$ is uniform on $u(\cdot)$.

By \eqref{system} for $\hat{x}(s)$, we  have
\begin{align*}
	&\varphi(t+\delta,\hat{x}(t+\delta)) -\varphi(t,x) -\frac{1}{2}\int_{t}^{t+\delta}|\hat{\beta}[u](s)|^2 ds \\
	&=\int_t^{t+\delta} \frac{\partial \varphi}{\partial s}(s,\hat{x}(s))
		+(f(\hat{x}(s),u(s)) +\sigma(\hat{x}(s),u(s))\hat{\beta}[u](s))
		\cdot \nabla \varphi(s,\hat{x}(s))
		-\frac{1}{2}|\hat{\beta}[u](s)|^2	ds \\
	&= \frac{\partial \varphi}{\partial t}(t,x) \delta
		+\int_t^{t+\delta} \bigg\{
		(f(\hat{x}(s),u(s)) +\sigma(\hat{x}(s),u(s))\hat{\beta}[u](s))
		\cdot \nabla \varphi(s,\hat{x}(s)) \\
	& \qquad \qquad \qquad \qquad \qquad \qquad \qquad \qquad
		-\frac{1}{2}|\hat{\beta}[u](s)|^2	\bigg\} ds + o(\delta),
\end{align*}
where $o(\delta)$ is uniform on $u(\cdot)$.
Thus, we have from \eqref{u_est_delta}
\begin{multline*}
	\varphi(t+\delta,\hat{x}(t+\delta)) -\varphi(t,x) -\frac{1}{2}\int_{t}^{t+\delta}|\hat{\beta}[u](s)|^2 ds \\
	\geq
		\left( \frac{\partial \varphi}{\partial t}(t,x)
		+\mathcal{H}(x, \varphi(t,x)+\rho, \nabla \varphi(t,x))\right) \delta
		+o(\delta).
\end{multline*}
Therefore, we obtain
\begin{align}
	&F^{u,\hat{\beta}[u]}_{t,t+\delta}\varphi(t+\delta,\cdot)(x)-\varphi(t,x)  \notag \\
	&=\int_{[t,t+\delta]}^\oplus l(\hat{x}(s), u(s))ds \oplus \varphi(t+\delta,\hat{x}(t+\delta))
	 -\frac{1}{2}\int_{t}^{t+\delta}|\hat{\beta}[u](s)|^2 ds -\varphi(t,x) \notag \\
	&\geq \varphi(t+\delta,\hat{x}(t+\delta)) -\varphi(t,x) 
	-\frac{1}{2}\int_{t}^{t+\delta}|\hat{\beta}[u](s)|^2 ds \notag \\
	&\geq
		\left( \frac{\partial \varphi}{\partial t}(t,x)
		+\mathcal{H}(x, \varphi(t,x)+\rho, \nabla \varphi(t,x))\right) \delta
		+o(\delta).
	\label{l_est_delta1}
\end{align}

For the second case, we suppose there exists a subset $I_\delta \subset [t,t+\delta]$
with positive Lebesgue measure such that
$u(s) \not\in A(x,\varphi(t,x)+\rho)$ for $s \in I_\delta$, \textit{i.e.},
\[
	l(x,u(s)) >\varphi(t,x)+\rho, \ \forall s \in I_\delta.
\]
Since $l$ and $\varphi$ are Lipschitz, 
\[
	l(y,u(s))>\varphi(r,z)+\rho -L(|y-x|+|z-x|+\delta) \text{ for } s \in I_\delta, \
	t \leq r \leq t+\delta.
\]
Then, by \eqref{err_sol_ini}, we have
\[
	l(\hat{x}(s), u(s))>\varphi(r,\hat{x}(r))+\frac{\rho}{2}, \ s \in I_\delta, \ t \leq r \leq t+\delta
\]
for small $\delta>0$.
Therefore we have
\[
	\int_{[t,t+\delta]}^\oplus 	l(\hat{x}(s), u(s))ds > \varphi(r,\hat{x}(r))+\frac{\rho}{2},
	\ t \leq r \leq t+\delta.
\]
By using this,
\begin{align}
	F_{t,t+\delta}^{u,\hat{\beta}[u]}\varphi(t+\delta,\cdot)(x) -\varphi(t,x)
 	&=\int_{[t,t+\delta]}^\oplus 	l(\hat{x}(s), u(s))ds \oplus \varphi(t+\delta,\hat{x}(t+\delta)) \notag \\
	& \qquad \qquad \qquad 
	-\frac{1}{2}\int_{t}^{t+\delta}|\hat{\beta}[u](s)|^2ds -\varphi(t,x)  \notag \\
	&\geq \int_{[t,t+\delta]}^\oplus 	l(\hat{x}(s), u(s))ds
		-\frac{1}{2}\int_t^{t+\delta}|\hat{\beta}[u](s)|^2ds
		-\varphi(t,x) \notag \\
	&\geq \varphi(t,x) +\frac{\rho}{2}
	-\frac{1}{2}\int_t^{t+\delta}|\hat{\beta}[u](s)|^2ds -\varphi(t,x) \notag\\
	&\geq \frac{\rho}{2}-M \delta.
	\label{l_est_delta2}
\end{align}
for some constant $M>0$. Here we used
$\hat{\beta}[u](s)=\sigma (\hat{x}(s),u(s))^T \nabla \varphi(s,\hat{x}(s))$ is bounded.
Hence, by \eqref{l_est_delta1} and \eqref{l_est_delta2}, we obtain \eqref{l_est_delta}.

Now we shall prove \eqref{l_est_gen}. We first fix $\eps>0$.
Let $\alpha \in \Gamma(t,t+\delta)$ be taken arbitrarily.
It implies from Lemma \ref{lem_approx_str} that
there exist $u^{\eps, \delta} \in L^\infty([t,t+\delta]; U)$ and 
$v^{\eps,\delta} \in L^2[t,t+\delta]$ such that
\begin{equation}
	|\alpha[v^{\eps,\delta}](s) -u^{\eps,\delta}(s)| < \eps \delta, \ \forall s \in [t,t+\delta)
	\text{ and }
	\hat{\beta}[u^{\eps,\delta}](s)=v^{\eps,\delta}(s) \ \text{a.e.\,}s \in [t,t+\delta].
\label{err_eps_del}
\end{equation}
Let $x^{\eps,\delta}(s)$ (\textit{resp.\,}$ \hat{x}^{\eps,\delta}(s)$)
be the solution of \eqref{system} for
$u(s)=\alpha[v^{\eps,\delta}](s)$ and $v(s)=v^{\eps,\delta}(s)$
(\textit{resp.\,}$u(s)=u^{\eps,\delta}(s)$ and $v(s)=\hat{\beta}[u^{\eps,\delta}](s)$).
By using \eqref{err_eps_del}, we can show that
\[
	|x^{\eps,\delta}(s)-\hat{x}^{\eps,\delta}(s)| \leq C\eps \delta^2,
	\ t \leq s \leq t+\delta
\]
for some $C>0$.
Then, it is not difficult to see that
\begin{multline*}
	\int_{[t,t+\delta]}^\oplus
		l(x^{\eps,\delta}(s),\alpha[v^{\eps,\delta}](s))ds \oplus \varphi(t+\delta,x^{\eps,\delta}(t+\delta)) \\
	\geq \int_{[t,t+\delta]}^\oplus
		l(\hat{x}^{\eps,\delta}(s),u^{\eps,\delta}(s))ds
		\oplus \varphi(t+\delta,\hat{x}^{\eps,\delta}(t+\delta))-\tilde{C}\eps \delta
		-\tilde{C} \eps \delta^2
\end{multline*}
for some constant $\tilde{C}$.
Since $v^{\eps,\delta}=\hat{\beta}[u^{\eps,\delta}]$, 
\[
	F^{\alpha[v^{\eps,\delta}], v^{\eps,\delta}}_{t,t+\delta}\varphi(t+\delta,\cdot)(x)-\varphi(t,x)
	\geq F^{u^{\eps,\delta}, \hat{\beta}[u^{\eps,\delta}]}_{t,t+\delta}\varphi(t+\delta,\cdot)(x)
	-\varphi(t,x)
	-\tilde{C}\eps \delta - \tilde{C}\eps\delta^2.
\]
By \eqref{l_est_delta}, 
\begin{multline*}
	F^{\alpha[v^{\eps,\delta}], v^{\eps,\delta}}_{t,t+\delta}\varphi(t+\delta,\cdot)(x)-\varphi(t,x) \\
	\geq
	\min \left\{ \left(
		\frac{\partial \varphi}{\partial t}(t,x)+\mathcal{H}(x, \varphi(t,x)+\rho, \nabla \varphi(t,x)) \right)\delta
		+o (\delta),
	\frac{\rho}{2}-M\delta \right\} -\tilde{C}\eps \delta -\tilde{C} \eps \delta^2.
\end{multline*}
Note that $o(\delta)$ is uniform on $\alpha$ because \eqref{l_est_delta} is uniform on $u$.
Therefore we obtain
\begin{multline*}
	F_{t,t+\delta}\varphi(t+\delta,\cdot)(x)-\varphi(t,x) \\
	\geq \min \left\{ \left(
		\frac{\partial \varphi}{\partial t}(t,x)+\mathcal{H}(x, \varphi(t,x)+\rho, \nabla \varphi(t,x)) \right)\delta
		+o (\delta),
	\frac{\rho}{2}-M\delta \right\} -\tilde{C}\eps \delta -\tilde{C} \eps \delta^2.
\end{multline*}
Dividing the both side by $\delta>0$ and taking $\liminf$ as $\delta \to 0+$,
\[
	\liminf_{\delta \to 0+}
	\frac{1}{\delta} \left\{ F_{t,t+\delta}\varphi(t+\delta,\cdot)(x)-\varphi(t,x)\right\}
	\geq \frac{\partial \varphi}{\partial t}(t,x)+\mathcal{H}(x, \varphi(t,x)+\rho, \nabla \varphi(t,x))
		-\tilde{C}\eps.
\]
Then, by sending $\eps \to 0$, we have
\[
	\liminf_{\delta \to 0+}
	\frac{1}{\delta} \left\{ F_{t,t+\delta}\varphi(t+\delta,\cdot)(x)-\varphi(t,x)\right\} 
	\geq  \frac{\partial \varphi}{\partial t}(t,x)+\mathcal{H}(x, \varphi(t,x)+\rho, \nabla \varphi(t,x)).
\]
By Lemma \ref{H-semi-lim} (i), if we take $\rho \to 0+$, we finally have \eqref{l_est_gen}. \
$\qed$

Under the DPP and the estimates on the generator,
we can easily characterize $V(t,x)$ as a unique viscosity solution.
The uniqueness is implied by the comparison theorem in the next section.
\begin{thm} \label{ch_value}
$V(t,x)$ is the unique bounded Lipschitz continuous viscosity solution
of \eqref{u_Isaacs} (equivalently \eqref{QVI}) with the terminal condition \eqref{QVI_T}.
\end{thm}
\noindent
\textit{Proof.} \
First of all, we point out that we may suppose the smooth test function $\varphi(t,x)$
in Definition \ref{visc_def} is taken from $C_b^1((0,T) \times \mathbb{R}^n)$ 
because $V(t,x)$ is bounded and the definition of viscosity sub/supersolutions
only uses a local information of $W(t,x)-\varphi(t,x)$.

We only show that $V(t,x)$ is a subsolution of \eqref{u_Isaacs} since the supersolution part
is proved in a similar way.
Let $(\hat{t},\hat{x}) \in (0,T) \times \mathbb{R}^n$ be a maximum point
of $V(t,x)-\varphi(t,x)$ for  $\varphi \in C_b^1((0,T) \times \mathbb{R}^n)$
and $V(\hat{t},\hat{x})=\varphi(\hat{t},\hat{x})$.
If we use \eqref{DPP} at $(t,x)=(\hat{t},\hat{x})$, then
\[
	V(\hat{t},\hat{x})
	=\inf_{\alpha \in \Gamma(\hat{t},\hat{t}+\delta)}
	E_{\hat{t}\hat{x}}^{+}\left[ \int_{[\hat{t},\hat{t}+\delta]}^\oplus l(x(s),\alpha[v](s))ds
	\oplus V(\hat{t}+\delta,x(\hat{t}+\delta)) \right].
\]
Since $V(\hat{t},\hat{x})=\varphi(\hat{t},\hat{x})$ and $V(t,x) \leq \varphi(t,x)$
for any $(t,x) \in (0,T) \times \mathbb{R}^n$, 
\[
	\varphi(\hat{t},\hat{x})
	\leq\inf_{\alpha \in \Gamma(\hat{t},\hat{t}+\delta)}
	 E^{+}_{\hat{t}\hat{x}}\left[ \int_{[\hat{t},\hat{t}+\delta]}^\oplus l(x(s),\alpha[v](s))ds
	\oplus \varphi(\hat{t}+\delta,x(\hat{t}+\delta)) \right].
\]
By using the notation $F_{t,t+\delta}$, this inequality can be written as 
\[
	\varphi(\hat{t},\hat{x}) \leq F_{\hat{t},\hat{t}+\delta}\varphi(\hat{t}+\delta,\cdot)(\hat{x}).
\]
Thus, we have
\[
	0 \leq \limsup_{\delta \to 0+}\frac{1}{\delta}
	\{  F_{\hat{t},\hat{t}+\delta}\varphi(\hat{t}+\delta,\cdot)(\hat{x})-\varphi(\hat{t},\hat{x}) \}.
\]
From Theorem \ref{cal_gen},
\[
	0 \leq \frac{\partial \varphi}{\partial t}(\hat{t},\hat{x})
	+\mathcal{H}^{\ast}(\hat{x},\varphi(\hat{t},\hat{x}),\nabla \varphi(\hat{t},\hat{x})).
\]
Therefore $V(t,x)$ is a viscosity subsolution of \eqref{u_Isaacs}.

The uniqueness is a consequence of Theorem \ref{comp} in the next section. $\qed$ 

\subsection{Remarks on other strategy classes}
To discuss the max-plus control problems, 
one might want to choose a different class of strategies.
Actually, there exist classes for which the generator-DPP arguments work.
We shall make remarks on other possibilities for strategy classes without details.

First of all, one might think that Elliott-Kalton strategy class should be a natural candidate for the problem.
Let $V_{EK}(t,x)$ be the value function of \eqref{value}
defined by $\Gamma_{EK}(t,T)$ instead of $\Gamma(t,T)$.
It is not difficult to see $V_{EK}(t,x)$ satisfies the DPP in a similar way to
the case of $\Gamma(t,T)$: for $ t \leq r \leq T$, $x \in \mathbb{R}^n$,
\[
	V_{EK}(t,x)= \inf_{\alpha \in \Gamma_{EK}(t,r)}
	E^{+}_{tx} \left[ \int_{[t,r]}^\oplus l(x(s), \alpha[v](s))ds \oplus V_{EK}(r,x(r)) \right].
\]
For $t<r$,  denote by $F_{t,r}^{EK}$ the operator associated with $\Gamma_{EK}(t,r)$:
\[
	F_{t,r}^{EK}\phi(x)= \inf_{\alpha \in \Gamma_{EK}(t,r)}
	E^{+}_{tx} \left[ \int_{[t,r]}^\oplus l(x(s), \alpha[v](s))ds \oplus \phi(x(r)) \right].
\]
If we calculate the generator,
we can prove that for $\varphi \in C_b^1((0,T) \times \mathbb{R}^n)$ and
$(t,x) \in (0,T) \times \mathbb{R}^n$,
\begin{gather}
	\limsup_{ \delta \to 0+} \frac{1}{\delta}
	\left\{ F_{t,t+\delta}^{EK} \varphi(t+\delta,\cdot)(x) - \varphi (t,x) \right\}
	\leq \frac{\partial \varphi}{\partial t}(t,x)
		+\mathcal{K}^\ast(x,\varphi(t,x), \nabla \varphi(t,x)), 
	\label{EK-gen_u} \\
	\frac{\partial \varphi}{\partial t}(t,x)
	+\mathcal{K}_\ast(x,\varphi(t,x), \nabla \varphi(t,x))
	 \leq 
	 \liminf_{ \delta \to 0+} \frac{1}{\delta}
	\left\{ F_{t,t+\delta}^{EK} \varphi(t+\delta,\cdot)(x) - \varphi (t,x) \right\},
	\label{EK-gen_l}
\end{gather}
where 
\begin{equation}
	\mathcal{K}(x,r,p)
	=\max_{v \in \mathbb{R}^d} \min_{u \in A(x,r)}
		\left\{ (f(x,u) +\sigma (x,u)v) \cdot p - \frac{1}{2}|v|^2  \right\}.
\label{l_Hamlitonian}
\end{equation}
See Appendix \ref{EK_gen} for the proof.
Thus, $V_{EK}(t,x)$ can be shown to be a bounded Lipschitz continuous viscosity solution of
\begin{equation}
	\frac{\partial V_{EK}}{\partial t}(t,x)
	+\mathcal{K}(x,V_{EK}(t,x), \nabla V_{EK}(t,x))=0, \ (t,x) \in (0,T) \times \mathbb{R}^n,
	\label{l_Isaacs}
\end{equation}
with the terminal condition \eqref{QVI_T}.
Moreover, $V_{EK}(t,x)$ can be shown to be unique in such function class
since we can have a comparison theorem for \eqref{l_Isaacs}
by using similar ideas to Theorem \ref{comp}.
Here we note that $\mathcal{K}(x,r,p) \leq \mathcal{H}(x,r,p)$ and they do not coincide in general.
Since our DPE is \eqref{QVI}, equivalently \eqref{u_Isaacs},  
using $\Gamma_{EK}(t,T)$ leads to a wrong equation.

There can be another class associated with \eqref{u_Isaacs}.
$\alpha \in \Gamma_{EK}(t,T)$ is called a strictly progressive strategy
if for any $\beta \in \Delta_{EK} (t,T)$,  there exist $u \in L^\infty([t,T];U)$ and
$v \in L^2[t,T]$ such that
\[
	\alpha[v]=u, \ \beta[u]=v \text{ a.e.\,on } [t,T],
\]
where $\Delta_{EK}(t,T)$ is the set of Elliott-Kalton strategies from $L^\infty([t,T];U)$ into $L^2[t,T]$
(see \cite{F02}, \cite{FS06} and \cite{KS05}). 
If we denote by $\Gamma_{SP}(t,T)$ the set of strictly progressive strategies on $[t,T]$
and $V_{SP}(t,x)$ is the value function given by $\Gamma_{SP}(t,T)$,
it is also seen that $V_{SP}(t,x)$ satisfies the DPP:
\[
	V_{SP}(t,x)=\inf_{\alpha \in \Gamma_{SP}(t,r)} 
	E_{tx}^+ \left[ \int_{[t,r]}^\oplus l(x(s), \alpha[v](s))ds \oplus V_{SP}(r,x(r)) \right]
\]
for $t \leq r \leq T $ and $x \in \mathbb{R}^n$
(cf. \cite{KS05} for max-plus multiplicative running cost).
Let the operator corresponding to $\Gamma_{SP}(t,r)$ be
\[
	F_{t,r}^{SP}\phi(x) 
	=\inf_{\alpha \in \Gamma_{SP}(t,r)} 
	E_{tx}^+ \left[ \int_{[t,r]}^\oplus l(x(s), \alpha[v](s))ds \oplus \phi(x(r)) \right].
\]
We can show that the generator of $F_{t,r}^{SP}$ is associated with $\mathcal{H}$
by the exactly same way as $\Gamma(t,r)$ without using Lemma \ref{lem_approx_str}:
\begin{gather*}
	\limsup_{ \delta \to 0+} \frac{1}{\delta}
	\left\{ F_{t,t+\delta}^{SP} \varphi(t+\delta,\cdot)(x) - \varphi (t,x) \right\}
	\leq \frac{\partial \varphi}{\partial t}(t,x)
		+\mathcal{H}^\ast(x,\varphi(t,x), \nabla \varphi(t,x)), \\
	\frac{\partial \varphi}{\partial t}(t,x)
	+\mathcal{H}_\ast(x,\varphi(t,x), \nabla \varphi(t,x))
	 \leq 
	 \liminf_{ \delta \to 0+} \frac{1}{\delta}
	\left\{ F_{t,t+\delta}^{SP} \varphi(t+\delta,\cdot)(x) - \varphi (t,x) \right\}
\end{gather*}
for $\varphi \in C_b^1((0,T) \times \mathbb{R}^n)$ and $(t,x) \in (0,T) \times \mathbb{R}^n$.
Therefore, $V_{SP}(t,x)$ is the unique bounded Lipschitz continuous viscosity solution
of \eqref{u_Isaacs} with \eqref{QVI_T}.
Hence, $V_{SP}(t,x)=V(t,x)$.
But we note that it is not clear that Markov control policies are strictly progressive in general.
Thus, $\Gamma(t,T)$ may be a better class if we consider the connection to 
the verification theorem.
See also \cite[Remark XI.9.1]{FS06} for the discussion on Markov control policies
of a differential game problem.
Lemma \ref{lem_approx_str} states that a property like strict progressivity holds
approximately, if $\alpha \in \Gamma(t,T)$.
This approximation property is useful in deriving Theorem \ref{cal_gen},
instead of directly using strict progressivity.

We have not made use of the following differential game interpretation of
the max-plus control problem. 
In this interpretation, $u(s)$ is a minimizing control and $v(s)$ is a maximizing control. 
The differential game dynamics are \eqref{system} and the game payoff is
\begin{equation}
	P(t,x; u,v) = \int^\oplus_{[t,T]} l (x(s), u(s))ds - \frac{1}{2} \int^T_t |v(s)|^2 ds. 
\label{game-pf}
\end{equation}
According to \eqref{value}
\begin{equation}
	V(t,x) = \inf_{\alpha \in \Gamma(t,T)} \sup_{v\in L^{2}[t,T]} P (t,x; \alpha [v], v).
\label{i-s_value}
\end{equation}
In the Elliott-Kalton definition of lower differential game value, $\Gamma(t,T)$ is replaced by
$\Gamma_{EK}(t,T)$.
Hence we should not expect $V(t,x)$ to agree in all cases with the Elliott-Kalton lower value.

The right side of \eqref{game-pf} involves both a max-plus integral and an ordinary integral.
For this reason we could not use standard differential game arguments
to obtain a dynamic programming principle.
Instead, properties of max-plus conditional expectations were used
to prove the dynamic programming principle in Theorem \ref{thm_DPP}.

A different kind of differential game payoff is
\begin{equation}
	P_1(t,x;u,v) = \int^T_t l(x(s), u(s)) - \frac{1}{2} \int_t^T |v(s)|^2 ds. 
\end{equation}
This is the case of ``max-plus multiplicative running cost" mentioned in the Introduction.
For this payoff, standard differential game methods can be used. Let
\begin{equation}
V_1(t,x) = \inf_{\alpha \in \Gamma_{SP}(t,T)} \sup_{v\in L^{2}[t,T]} P_1(t,x; \alpha [v],v).
\end{equation}
Then $V_1(t,x)$ agrees with the upper Elliott-Kalton value for the differential game (not the lower value). See \cite[ p.\ 393]{FS06}.

\section{Comparison theorem}
In the proof of Theorem \ref{ch_value}, we used the uniqueness
of bounded Lipschitz continuous viscosity solutions of \eqref{u_Isaacs} with \eqref{QVI_T}.
This is a consequence of a comparison theorem stated in the present section.
So far, we  only  consider continuous viscosity solutions.
We prove a comparison theorem for semi-continuous cases
because we need it in the argument for risk-sensitive limits.



As we noted, the Hamiltonian $\mathcal{H}(x,r,p)$ is not sufficiently regular to apply the general comparison
results of viscosity solutions. Following the ideas used in the proof of \cite[Theorem 4.2]{BI89}
with Lemmas \ref{est_set}--\ref{H-semi-lim},
we will get a weak comparison theorem for \eqref{u_Isaacs}.
The statement of the weak  comparison theorem is similar to \cite[Theorems 8.1 and 8.2]{FS06}
by assuming  a Lipschitz viscosity solution.

\begin{thm} \label{comp}
Let $W(t,x)$ (\textup{resp.\,}$\overline{W}(t,x)$) be a bounded function  on $(0,T] \times \mathbb{R}^n$
and a viscosity subsolution
(\textup{resp.\,}viscosity supersolution) of \eqref{u_Isaacs}.
Assume that $x \mapsto W^{\ast}(T,x)$ and $x \mapsto \overline{W}_\ast (T,x)$ are continuous and 
there exists a bounded continuous viscosity solution $U \in C((0,T] \times \mathbb{R}^n)$
of \eqref{u_Isaacs} such that 
$x \mapsto U(t,x)$ is Lipschitz continuous uniformly on $(0,T)$
and $W^\ast(T,x) \leq U(T,x)  \leq \overline{W}_\ast(T,x)$ for $x \in \mathbb{R}^n$.
Then $W^\ast(t,x) \leq U(t,x)$ and $U(t,x) \leq \overline{W}_\ast(t,x)$ on $(0,T] \times \mathbb{R}^n$.
Hence, $W^\ast(t,x) \leq \overline{W}_\ast(t,x)$ on $(0,T] \times \mathbb{R}^n$.
\end{thm}
\noindent
\textit{Proof.} \
We shall only show $U(t,x) \leq \overline{W}_\ast(t,x)$.
 $W^\ast(t,x) \leq U(t,x)$ can be proved in a similar way
if we change the role of $U(t,x)$ (\textit{resp.\,} $\overline{W}_\ast(t,x)$)
to $W^\ast(t,x)$ (\textit{resp.\,}$U(t,x)$).
Let $0<\theta<1$ and $\beta>0$. 
We choose a smooth nonnegative function $g :\mathbb{R}^n \to \mathbb{R}$ such that
\begin{gather}
	\overline{W}_\ast(t,x) \leq g(x), \ \forall (t,x) \in (0,T] \times \mathbb{R}^n, \label{g_dom}\\
	\| \nabla g \|_{\infty} < \infty, \ g(x) \to \infty \ (|x| \to \infty). \label{bd_grad_g}
\end{gather}
Define $\overline{W}^{\theta, \beta} : (0,T] \times \mathbb{R}^n \to \mathbb{R}$ by
\[
	\overline{W}^{\theta,\beta}(t,x) =
	(1-\theta)\overline{W}_\ast(t,x)+\theta g(x)+L\theta (T-t)+\frac{\beta}{t},
\]
where $L>0$ is a constant.

We can show that 
there exists a constant $L>0$ depending on $\| f \|_\infty$, $\| \sigma \|_\infty$, $\| \nabla g \|_\infty$
such that
\begin{gather}
	\frac{\partial \overline{W}^{\theta,\beta}}{\partial t}+
		\mathcal{H}_{\ast}(x,\overline{W}^{\theta,\beta}(t,x), \nabla \overline{W}^{\theta,\beta}(t,x))
		\leq -\frac{\beta}{t^2}
		\text{ in } (0,T) \times \mathbb{R}^n \text{ (viscosity sense)}, \label{eq_approx_super} \\
	\overline{W}_\ast \leq \overline{W}^{\theta,\beta} \text{ in } (0,T] \times \mathbb{R}^n. \label{dom_approx_super}
\end{gather}
These can be seen if $\overline{W}(t,x)$ is a classical supersolution in the following way.
The viscosity case can be also verified by using the same kind of arguments.

By $L>0$, $\beta>0$ and \eqref{g_dom}, 
\begin{align*}
	\overline{W}^{\theta,\beta}(t,x)
	&= (1-\theta)\overline{W}(t,x) +\theta g(x) +L\theta (T-t)+\frac{\beta}{t} \\
	&\geq (1-\theta)\overline{W}(t,x)+\theta g(x) \\
	&\geq (1-\theta)\overline{W}(t,x)+\theta \overline{W}(t,x) 
	=\overline{W}(t,x).
\end{align*}
Thus we have \eqref{dom_approx_super}.

We shall show \eqref{eq_approx_super}. Since $p \mapsto H^u (x,p)$ is convex,
\[
	H^u(x,\nabla \overline{W}^{\theta,\beta}(t,x))
	\leq (1-\theta)H^u (x, \nabla \overline{W}(t,x)) + \theta H^u (x, \nabla g(x)).
\]
By noting that $f$, $\sigma$ and $\nabla g$ are bounded, there exists $C>0$
depending only on $\|f \|_\infty$, $\| \sigma \|_\infty$, $\| \nabla g \|_\infty$ such that
\[
	H^u(x,\nabla \overline{W}^{\theta,\beta}(t,x))
	\leq (1-\theta)H^u (x, \nabla \overline{W}(t,x)) + \theta C.
\]
Then, we have
\[
	\frac{\partial \overline{W}^{\theta,\beta}}{\partial t}
	+H^u(x, \nabla \overline{W}^{\theta,\beta}(t,x))
	\leq (1-\theta)\left( \frac{\partial \overline{W}}{\partial t}+H^u(x, \nabla \overline{W}(t,x)) \right)
	-(L-C)\theta-\frac{\beta}{t^2}.
\]
If we take $L>0$ such that $L>C$,
\[
\frac{\partial \overline{W}^{\theta,\beta}}{\partial t}
	+H^u(x, \nabla \overline{W}^{\theta,\beta}(t,x))
	\leq (1-\theta)\left( \frac{\partial \overline{W}}{\partial t}+H^u(x, \nabla \overline{W}(t,x)) \right)
	-\frac{\beta}{t^2}.
\]
Taking the minimum over $u \in A(x,\overline{W}(t,x))$, we obtain
\begin{align*}
	&\frac{\partial \overline{W}^{\theta,\beta}}{\partial t}
	+\min_{u \in A(x,\overline{W}(t,x))} H^u(x, \nabla \overline{W}^{\theta,\beta}(t,x))  \\
	&\leq (1-\theta)\left( \frac{\partial \overline{W}}{\partial t}
	+\mathcal{H}(x, \overline{W}(t,x), \nabla \overline{W}(t,x)) \right)
	-\frac{\beta}{t^2} \\
	&\leq -\frac{\beta}{t^2}.
\end{align*}
In the last line, we used that $\overline{W}(t,x)$ is a supersolution of \eqref{u_Isaacs}.
Since $\overline{W}(t,x) \leq \overline{W}^{\theta, \beta}(t,x)$, we have
$A(x,\overline{W}(t,x)) \subset A(x, \overline{W}^{\theta, \beta}(t,x))$. Thus, we have
\[
	\mathcal{H}(x,\overline{W}^{\theta,\beta}(t,x), \nabla \overline{W}^{\theta, \beta}(t,x))
	\leq \min_{u \in A(x,\overline{W}(t,x))} H^u(x, \nabla \overline{W}^{\theta,\beta}(t,x))
\]
Therefore we obtain \eqref{eq_approx_super}.

We shall  prove $U \leq \overline{W}^{\theta,\beta}$ in $(0,T] \times \mathbb{R}^n$.
For $\eps>0$, we define $\Phi_\eps : (0,T] \times (0,T] \times \mathbb{R}^n \times \mathbb{R}^n
\to \mathbb{R}$ by
\[
	\Phi_\eps(t,s,x,y) = U(t,x)-\overline{W}^{\theta,\beta}(s,y)
	-\frac{1}{2\eps}|t-s|^2-\frac{1}{2\eps}|x-y|^2.
\]
By using standard arguments in viscosity theory, it is not difficult to see that
for sufficiently small $\eps>0$, 
there exists $(t_\eps, s_\eps,x_\eps,y_\eps)
\in  (0,T] \times (0,T] \times \mathbb{R}^n \times \mathbb{R}^n$ such that
\begin{gather}
	\Phi_\eps (t_\eps, s_\eps, x_\eps, y_\eps) =
	\sup_{\substack{t,s \in (0,T] \\ x,y \in \mathbb{R}^n}} \Phi_\eps(t,s,x,y), 
	\label{exist_eps-max} \\
	\theta g(y_\eps) +\frac{\beta}{s_\eps} +\frac{1}{2\eps}|t_\eps-s_\eps|^2
	+\frac{1}{2\eps}|x_\eps-y_\eps|^2 \leq M, \label{est_eps-opt}
\end{gather}
where $M>0$ is a constant independent of $\eps>0$. 
Again, by standard arguments, we can show that
\begin{equation}
	\frac{1}{2\eps}|t_\eps-s_\eps|^2 +\frac{1}{2\eps}|x_\eps-y_\eps|^2 \to 0
	\ (\eps \to 0).
\label{fine_asym}
\end{equation}

Now we are ready to prove $U \leq \overline{W}^{\theta,\beta}$ in $(0,T] \times \mathbb{R}^n$, 
equivalently,
\[
	\max_{(0,T] \times \mathbb{R}^n} (U -\overline{W}^{\theta,\beta}) \leq 0.
\]
On the contrary, we suppose that
\begin{equation}
		\max_{(0,T] \times \mathbb{R}^n} (U -\overline{W}^{\theta,\beta})>0.
	\label{a_contr}
\end{equation}
Under the assumption \eqref{a_contr}, it can be seen that
\begin{equation}
	0< t_\eps, s_\eps <T \text{ for sufficiently small } \eps>0.
\label{int_max}
\end{equation}
If \eqref{int_max} does not hold, there exists a sequence $\{ \eps_n \}_{n=1}^{\infty}$
($\eps_n \downarrow 0$) such that for each $n=1,2,\cdots$,
\[
	t_{\eps_n} =T \text{ or } s_{\eps_n}=T.
\]
Since $|t_\eps - s_\eps| \to 0$ as $\eps \to 0$ by \eqref{est_eps-opt},
\[
	t_{\eps_n}, s_{\eps_n} \to T \ (n \to \infty).
\]
As for the asymptotics on $x_\eps$, $y_\eps$ through $\eps=\eps_n$,
we can see from \eqref{est_eps-opt}  that there exists $\bar{x} \in \mathbb{R}^n$ such that
\[
	x_{\eps_n}, y_{\eps_n} \to \bar{x} \ (n \to \infty)
	\text{ by taking a subsequence of } \{ \eps_n \}.
\]
We estimate $U(t_{\eps_n}, x_{\eps_n})-\overline{W}^{\theta,\beta}(s_{\eps_n}, y_{\eps_n})$ from the above.
In the case where $t_{\eps_n}=T$, 
\begin{align*}
	U(t_{\eps_n},x_{\eps_n})-\overline{W}^{\theta,\beta}(s_{\eps_n},y_{\eps_n})
	&= U(T, x_{\eps_n})-	\overline{W}^{\theta,\beta}(s_{\eps_n},y_{\eps_n}) \\
	&= U(T,x_{\eps_n}) - \overline{W}^{\theta,\beta}(T,x_{\eps_n})
		+\overline{W}^{\theta,\beta}(T,x_{\eps_n})-\overline{W}^{\theta,\beta}(s_{\eps_n},y_{\eps_n}) \\
	&\leq \overline{W}^{\theta,\beta}(T,x_{\eps_n})-\overline{W}^{\theta,\beta}(s_{\eps_n},y_{\eps_n}).
\end{align*}
In the last line, we used $U \leq \overline{W}_\ast \leq  \overline{W}^{\theta,\beta} $
on $\{ t=T \} \times \mathbb{R}^n$.
In the case where $s_{\eps_n}=T$,
\begin{align*}
	U(t_{\eps_n},x_{\eps_n})-\overline{W}^{\theta,\beta}(s_{\eps_n},y_{\eps_n})
	&= U(t_{\eps_n}, x_{\eps_n})-\overline{W}^{\theta,\beta}(T,y_{\eps_n}) \\
	&= U(t_{\eps_n},x_{\eps_n})-U(T,y_{\eps_n})
		+U(T,y_{\eps_n})-\overline{W}^{\theta,\beta}(T,y_{\eps_n}) \\
	&\leq  U(t_{\eps_n},x_{\eps_n})-U(T,y_{\eps_n}).
\end{align*}
Similarly, we used $U  \leq \overline{W}^{\theta,\beta} $ on $\{ t=T \} \times \mathbb{R}^n$.
Therefore, we have
\[
	U(t_{\eps_n}, x_{\eps_n})-\overline{W}^{\theta,\beta}(s_{\eps_n}, y_{\eps_n}) 
	\leq
	\max \{ \overline{W}^{\theta,\beta}(T,x_{\eps_n})-\overline{W}^{\theta,\beta}(s_{\eps_n},y_{\eps_n}),
		U(t_{\eps_n},x_{\eps_n})-U(T,y_{\eps_n})\} .
\]
Taking limsup as ${\eps_n \to 0}$,
\begin{align*}
	&\limsup_{\eps_n \to 0}(U(t_{\eps_n}, x_{\eps_n})-\overline{W}^{\theta,\beta}(s_{\eps_n}, y_{\eps_n}) )\\
	& \leq \limsup_{\eps_n \to 0}
	\max \{ \overline{W}^{\theta,\beta}(T,x_{\eps_n})-\overline{W}^{\theta,\beta}(s_{\eps_n},y_{\eps_n}),
		U(t_{\eps_n},x_{\eps_n})-U(T,y_{\eps_n})\}  \\
	&\leq \max\{
		\limsup_{\eps_n \to 0}
		(\overline{W}^{\theta,\beta}(T,x_{\eps_n})-\overline{W}^{\theta,\beta}(s_{\eps_n},y_{\eps_n})),
		  \limsup_{\eps_n \to 0} (U(t_{\eps_n},x_{\eps_n})-U(T,y_{\eps_n}))\}
\end{align*}
Noting that $x \mapsto \overline{W}^{\theta, \beta}(T,x)$ is continuous,
$\overline{W}^{\theta,\beta}$ is lower semi-continuous
and $U(t,x)$ is continuous, we have
\begin{equation}
	\limsup_{\eps_n \to 0}(U(t_{\eps_n}, x_{\eps_n})-\overline{W}^{\theta,\beta}(s_{\eps_n}, y_{\eps_n}) )
	\leq 0.
\label{asym_contr}
\end{equation}

On the other hand, from \eqref{a_contr},
\[
	0<\max_{(0,T] \times \mathbb{R}^n} (U-\overline{W}^{\theta,\beta})
	\leq \max_{\substack{t,s \in (0,T] \\ x, y \in \mathbb{R}^n}}\Phi_\eps(t,s,x,y)
	\leq U(t_\eps, x_\eps) -\overline{W}^{\theta, \beta} (s_\eps, y_\eps).
\]
Thus, we have
\[
	0< \max_{(0,T] \times \mathbb{R}^n} (U-\overline{W}^{\theta,\beta})
	\leq \liminf_{\eps_n \to 0} (U(t_{\eps_n}, x_{\eps_n})-\overline{W}^{\theta,\beta}(s_{\eps_n},y_{\eps_n})),
\]
which contradicts to \eqref{asym_contr}. Therefore, \eqref{int_max} has to hold.

For sufficiently small $\eps>0$, it implies from \eqref{int_max} that
$(t_\eps, s_\eps, x_\eps, y_\eps)$ is a maximum point of $\Phi_\eps(t,s,x,y)$
in $(0,T) \times (0,T) \times \mathbb{R}^n \times \mathbb{R}^n$. 
In particular, since $(t_\eps,x_\eps)$
is a maximum point of $\Phi_\eps(t,s_\eps ,x,y_\eps)$
and $U(t,x)$ is a viscosity subsolution, we have
\begin{equation}
	\frac{1}{\eps}(t_\eps-s_\eps)
	+\mathcal{H}^\ast (x_\eps, U(t_\eps, x_\eps), \frac{1}{\eps}(x_\eps-y_\eps)) \geq 0.
\label{eps-sub}
\end{equation}
In a similar way, since $(s_\eps,y_\eps)$ is a minimum point of
$-\Phi_\eps(t_\eps,s,x_\eps,y)$ and $\overline{W}^{\theta,\beta}$ is
a viscosity supersolution of \eqref{eq_approx_super}, 
\begin{equation}
	\frac{1}{\eps}(t_\eps-s_\eps)
	+\mathcal{H}_\ast(y_\eps, \overline{W}^{\theta,\beta}(s_\eps,y_\eps), \frac{1}{\eps}(x_\eps-y_\eps))
	\leq -\frac{\beta}{s_\eps^2}.
\label{eps-super}
\end{equation}
Subtracting \eqref{eps-super} from \eqref{eps-sub}, we have
\begin{equation}
	\mathcal{H}^\ast (x_\eps, U(t_\eps, x_\eps), \frac{1}{\eps}(x_\eps-y_\eps)) 
	-\mathcal{H}_\ast(y_\eps, \overline{W}^{\theta,\beta}(s_\eps,y_\eps), \frac{1}{\eps}(x_\eps-y_\eps))
	\geq \frac{\beta}{s_\eps^2}.
\label{eps-ineq}
\end{equation}

Let us take $\alpha>0$ such that
\[
	0< \alpha < \max_{(0,T] \times \mathbb{R}^n}(U-\overline{W}^{\theta,\beta}).
\]
Since $\max_{(0,T] \times \mathbb{R}^n} (U-\overline{W}^{\theta,\beta})
\leq U(t_\eps,x_\eps) -\overline{W}^{\theta,\beta}(s_\eps, y_\eps)$ for $\eps>0$,
\[
	\alpha < U(t_\eps,x_\eps) -\overline{W}^{\theta,\beta}(s_\eps, y_\eps), \ \forall \eps>0.
\]
We choose $r^\eps_1<r^\eps_2$ such that
\[
	\overline{W}^{\theta,\beta}(s_\eps, y_\eps)<r^\eps_1<r^\eps_2<U(t_\eps, x_\eps)  \text{ and }
	r_2^\eps-r_1^\eps=\alpha.
\]
By Lemma \ref{H-err} (i) and Lemma \ref{H-semi-lim} (ii),
\begin{align*}
	\mathcal{H}^\ast(x_\eps,U(t_\eps,x_\eps),\frac{1}{\eps}(x_\eps-y_\eps))
	&=\mathcal{H} (x_\eps, U(t_\eps, x_\eps)-0, \frac{1}{\eps}(x_\eps-y_\eps)) \\
	&\leq \mathcal{H} (x_\eps, r^\eps_2, \frac{1}{\eps}(x_\eps-y_\eps)).
\end{align*}
By Lemma \ref{H-err} (i) and Lemma \ref{H-semi-lim} (i),
\begin{align*}
	\mathcal{H}_\ast (y_\eps, \overline{W}^{\theta,\beta}(s_\eps,y_\eps),\frac{1}{\eps}(x_\eps-y_\eps))
	&=\mathcal{H} (y_\eps, \overline{W}^{\theta,\beta}(s_\eps,y_\eps)+0, \frac{1}{\eps}(x_\eps-y_\eps)) \\
	&\geq \mathcal{H} (y_\eps, r_1^\eps, \frac{1}{\eps}(x_\eps-y_\eps)).
\end{align*}
Thus, we have from \eqref{eps-ineq}
\begin{equation}
	\mathcal{H} (x_\eps, r^\eps_2, \frac{1}{\eps}(x_\eps -y_\eps))
	-\mathcal{H}(y_\eps, r^\eps_1, \frac{1}{\eps}(x_\eps-y_\eps)) \geq \frac{\beta}{s_\eps^2}.
\label{eps-ineq2}
\end{equation}
By Lemma \ref{H-err} (ii), there exists $\delta=\delta(|r^\eps_2-r^\eps_1|)=\delta(\alpha)$
such that
\[
	\mathcal{H}(y,r^\eps_1,p) \geq \mathcal{H} (x_\eps, r^\eps_2, p) -L(|p|+1)|p| |x_\eps -y|, \
	\forall y \in B_\delta (x_\eps), \forall p \in \mathbb{R}^n.
\]
In particular, taking $p=(x_\eps -y)/\eps$,
\[
	\mathcal{H}(y,r^\eps_1,\frac{1}{\eps}(x_\eps-y))
	\geq \mathcal{H}(x_\eps, r^\eps_2, \frac{1}{\eps}(x_\eps-y))
	-L \left( \left| \frac{x_\eps -y}{\eps} \right| +1 \right)\frac{|x_\eps -y|^2}{\eps},
	\ \forall y \in B_{\delta}(x_\eps).
\]
Since $|x_\eps-y_\eps| \to 0$ $(\eps \to 0)$ and $\delta$ is independent of $\eps$,
\begin{gather*}
	\mathcal{H}(y_\eps,r^\eps_1,\frac{1}{\eps}(x_\eps-y_\eps))
	\geq \mathcal{H} (x_\eps, r^\eps_2, \frac{1}{\eps}(x_\eps-y_\eps))
	-L  \left( \left| \frac{x_\eps -y_\eps}{\eps} \right | +1 \right)\frac{|x_\eps -y_\eps|^2}{\eps} 
\end{gather*}
for sufficiently small $\eps>0$.
Thus, we obtain from \eqref{eps-ineq2}
\[
	\mathcal{H}(y_\eps,r^\eps_1,\frac{1}{\eps}(x_\eps-y_\eps))
	+L  \left( \left| \frac{x_\eps -y_\eps}{\eps} \right| +1 \right)\frac{|x_\eps -y_\eps|^2}{\eps}
	- \mathcal{H} (y_\eps, r^\eps_1, \frac{1}{\eps}(x_\eps-y_\eps)) \geq \frac{\beta}{s_\eps^2},
\]
which implies that
\begin{equation}
	\frac{\beta}{s_\eps^2}
	\leq L  \left( \left| \frac{x_\eps -y_\eps}{\eps}  \right| +1 \right)\frac{|x_\eps -y_\eps|^2}{\eps}	
	\text{ for sufficiently small } \eps>0.
\label{max_pt_err}
\end{equation}

Recall that $x \mapsto U(t,x)$ is Lipschitz continuous uniformly on $t$.
Then, the superdifferential $D^{+}_{x} U(t,x)$ is included in $\bar{B}_{K}(0)$,
where $K>0$ is a uniform Lipschitz constant of  $x \mapsto U(t,x)$.
Thus, since $(x_\eps-y_\eps)/\eps \in D^{+}_x U(t_\eps, x_\eps)$, 
\[
	\left| \frac{x_\eps-y_\eps}{\eps} \right| \leq K,
\]
Therefore, we have from \eqref{max_pt_err}
\[
	\frac{\beta}{s_\eps^2}
	\leq L  ( K +1)\frac{|x_\eps -y_\eps|^2}{\eps}	
	\text{ for sufficiently small } \eps>0.
\]
Since $|x_\eps - y_\eps|^2/\eps \to 0$ as $\eps \to 0$ from \eqref{fine_asym}, we have
\[
	\lim_{\eps \to 0} \frac{\beta}{s_\eps^2}=0.
\]
This contradicts to $s_\eps \in (0,T]$.
Therefore, we have $\max_{(0,T] \times \mathbb{R}^n}(U-\overline{W}^{\theta,\beta})\leq 0$, \textit{i.e.},
\[
	U \leq \overline{W}^{\theta, \beta} \text{ in } (0,T] \times \mathbb{R}^n.	
\]

Finally, if we take the limit as $\theta \downarrow 0$ and $\beta \downarrow 0$, we obtain
$$
	U \leq \overline{W}_\ast \text{ in } (0,T] \times \mathbb{R}^n.
\eqno{\qed}
$$

\section{Risk-sensitive control limits}
Let $\nu=(\Omega, \{ \mathcal{F}_s \}, P, \{ W(s) \})$ be a reference probability system
with $d$-dimensional $\{ \mathcal{F}_s \}$-Brownian motion $\{ W(s) \}$.
Under $\nu$, consider the following stochastic system
modelled by the stochastic differential equation:
\begin{equation}
	\left\{
	\begin{aligned}
		dX(s) &= f(X(s),u(s))ds + \theta^{-1/2} \sigma (X(s),u(s))dW(s), \ t \leq s \leq T, \\
		X(t) &= x \in \mathbb{R}^n,
	\end{aligned}
	\right.
\label{stoch_sys}
\end{equation}
where $\{ u(s) \}$ is a 
$U$-valued $\{ \mathcal{F}_s \}$-progressively measurable control process
and $\theta>0$ is a parameter.
For this system, we introduce the risk-sensitive type criterion
\begin{equation}
	J_\theta(t,x; u(\cdot)) = E_{tx}\left[ \int_t^T e^{\theta l(X(s),u(s))} ds \right].
\label{r-s_cr}
\end{equation}
The optimal control problem for \eqref{stoch_sys} with \eqref{r-s_cr}
is to minimize \eqref{r-s_cr} over the control processes. Thus, we are interested in
the value function
\begin{equation}
	\Psi_{\theta}(t,x)
	=\inf_{u(\cdot)} J_\theta (t,x; u(\cdot) ),
\label{r-s_value}
\end{equation}
where the infimum is taken over all $U$-valued $\{ \mathcal{F}_s \}$-progressively measurable processes.
If we set $L(x,u)=e^{-l(x,u)}$, the value function becomes
\[
	\Psi_{\theta}(t,x) =\inf_{u(\cdot)} E_{tx}\left[ \int_t^T L(X(s),u(s))^{-\theta} ds \right].
\]
In an example considered in Section 7, $\theta$ can be understood as the parameter related to risk-sensitivity
in optimal investment/consumption problems of mathematical finance.

Our aim in the present section is to connect the stochastic control problem 
of risk-sensitive type
described in the above with
the max-plus stochastic control via risk-averse limits.
 More precisely, we shall study the asymptotics of $\Psi_{\theta}(t,x)$
as $\theta \to \infty$.
The arguments we use are mainly from the PDE techniques in viscosity theory.
We utilize the general stability result of discontinuous solutions,
a so-called Barles-Perthame type procedure (cf.  \cite{BaP87}, \cite{FS06}).

Let us define $V_\theta(t,x)$ by
\[
	V_\theta(t,x) = \frac{1}{\theta}\log \Psi_\theta (t,x), \
	(t,x) \in (0,T) \times \mathbb{R}^n.
\]
To derive an equation of $V_\theta(t,x)$, we first note that under (A1)--(A3),
$\Psi_\theta(t,x)$ is a (unique) bounded uniformly continuous viscosity solution
of the following DPE (cf. \cite[Chap.V.9] {FS06}): 
\begin{gather}
	\frac{\partial \Psi_\theta}{\partial t}
	+F_{\theta}(x, \nabla \Psi_\theta(t,x), D^2 \Psi_\theta (t,x))=0,
	\ (t,x) \in (0,T) \times \mathbb{R}^n, 
	\label{DPE_Psi}\\
	\Psi_\theta (T,x)=0, \ x \in \mathbb{R}^n,
	\label{Psi_T} 
\end{gather}
where $F_\theta$ is defined for $\phi: \mathbb{R}^n \to \mathbb{R}$ by
\[
	F_\theta (x, \nabla \phi(x), D^2 \phi(x))
	=\min_{u \in U} 
	\left\{ \frac{1}{2\theta} \tr (a(x,u) D^2 \phi(x)) + f(x,u) \cdot \nabla \phi(x)
		+e^{\theta l(x,u)} \right\}
\]
with $a(x,u)=\sigma(x,u) \sigma(x,u)^T$.
Thus, $V_\theta(t,x)$ is a bounded continuous viscosity solution of
\begin{equation}
	\frac{\partial V_\theta}{\partial t}
	+\mathcal{H}_{\theta}(x,V_\theta(t,x), \nabla V_\theta(t,x), D^2 V_\theta (t,x))=0,
	\ (t,x) \in (0,T) \times \mathbb{R}^n,
\label{DPE_V_theta}
\end{equation}
where
\begin{align*}
	&\mathcal{H}_\theta(x, \phi(x), \nabla \phi(x), D^2 \phi(x)) \\
	&=\min_{u \in U} 
	\bigg\{ \frac{1}{2\theta} \tr (a(x,u) D^2 \phi(x))+ H^u (x,\nabla \phi(x))
		+e^{\theta (l(x,u)-\phi(x))} \bigg\}
\end{align*}
and $H^u(x,p)$ is defined by \eqref{H^u}.
Note that the value $V_\theta(t,x)$ at $t=T$  cannot be defined since $\Psi_\theta(T,x)=0$.


If we follow Barles-Perthame type procedure, we need to verify the following two steps:
(i) stability on solutions, \textit{i.e.}, the upper (\textit{resp.\,}lower) semi-continuous limit
of $V_\theta(t,x)$ as $\theta \to \infty$ is
a viscosity subsolution (\textit{resp.\,}supersolution) of the limit equation,
(ii) the semi-continuous limits have a common continuous  terminal data.
(i) will be implied by a general argument from viscosity theory.
(ii) can be proved by obtaining good estimates for the semi-continuous limits
in the following way.

We recall the semi-continuous limits of $V_\theta(t,x)$ as $\theta \to \infty$.
For  $(t,x) \in (0,T) \times \mathbb{R}^n$, 
define $\overline{V}(t,x)$ and $\underline{V}(t,x)$ by
\[
	\overline{V}(t,x)=\limsup_{\substack{\theta \to \infty \\ (s,y) \to (t,x)}} V_{\theta}(s,y), \
	\underline{V}(t,x)=\liminf_{\substack{\theta \to \infty \\ (s,y) \to (t,x)}} V_{\theta}(s,y).
	\
\]
Since $l(x,u)$ is bounded on $\mathbb{R}^n \times U$, $\overline{V}(t,x)$ and $\underline{V}(t,x)$
are well-defined.
\begin{prop} \label{est_semi-lim}
There exists $M>0$ such that
for any $(t,x) \in (0,T) \times \mathbb{R}^n$,
\begin{equation}
	\min_{u \in U}l(x,u) \leq \underline{V}(t,x)
	\leq \overline{V}(t,x) \leq \min_{u \in U}l(x,u) + M(T-t).	
\label{est_V-bar}
\end{equation}
\end{prop}
The following corollary is immediate from the above estimates.
\begin{cor} \label{est_terminal}
$\overline{V}(t,x)$ and $\underline{V}(t,x)$
are uniquely extended to an upper and a lower semi-continuous functions
on $(0,T] \times \mathbb{R}^n$ with $\overline{V}(T,x)=\underline{V}(T,x)=\min_{u \in U}l(x,u)$,
respectively.
\end{cor}
\noindent
\textit{Proof of Proposition \textup{\ref{est_semi-lim}}.} \
We shall first show the lower estimate in \eqref{est_V-bar}.
Take $\alpha>0$. If we set $\underline{l}(x)=\min_{u \in U}l(x,u)$,
we can find a smooth bounded function $g(x)$ 
with bounded first and second derivatives such that
\begin{equation}
	\alpha +g(x) \leq \underline{l}(x) \leq g(x) + 2\alpha, \ x \in \mathbb{R}^n.
\label{g_est}
\end{equation}
This can be seen because $\underline{l}(x)$ is Lipschitz and
we can take  a smooth bounded $\tilde{g}(x)$
with bounded first and second derivatives such that
\[
	-\frac{\alpha}{2} \leq \underline{l}(x)- \tilde{g}(x) \leq \frac{\alpha}{2}.
\]
If we choose $g(x)=\tilde{g}(x) -(3/2)\alpha$, $g(x)$ satisfies \eqref{g_est}.

Let us define $\phi(t,x)$ by
\[
	\phi(t,x)=(T-t) e^{\theta g(x)}= (T-t) G(x)^\theta.
\]
If we calculate the left-hand side (LHS) of \eqref{DPE_Psi} for $\phi(t,x)$,
\begin{align*}
	& \frac{\partial \phi}{\partial t}
	+F_{\theta}(x, \nabla \phi(t,x), D^2 \phi(t,x))  \\
	&=G(x)^\theta \bigg(-1 +\min_{u \in U}
	\bigg\{ \frac{1}{2} (T-t)  G^{-1}(x) \tr (a(x,u) D^2G(x) ) \\
	&	\qquad  \qquad \qquad \qquad \qquad 
	+\frac{1}{2}(T-t)(\theta-1)G^{-2}(x)a(x,u) \nabla G(x) \cdot \nabla G(x) \\
	& \qquad  \qquad \qquad \qquad \qquad 
		+(T-t)\theta G^{-1}(x) f(x,u) \cdot \nabla G(x) + e^{\theta (l(x,u)-g(x))} \bigg\} \bigg) \\ 
	&\geq G(x)^\theta \bigg( -1 + \min_{u \in U}
		\bigg\{ \frac{1}{2} (T-t)  G^{-1}(x) \tr (a(x,u) D^2G(x) ) \\
	& \qquad  \qquad \qquad \qquad \qquad 
		+(T-t)\theta G^{-1}(x) f(x,u) \cdot \nabla G(x) + e^{\theta (l(x,u)-g(x))} \bigg\} \bigg)
\end{align*}
Here we assumed $\theta>1$ without loss of generality.
Since $g(x)$, $\nabla g(x)$ and $D^2 g(x)$ are bounded , we  can see that 
there exists a constant $C>0$ such that
\[
	 \frac{1}{2} (T-t)  G^{-1}(x) \tr (a(x,u) D^2G(x) ) 
	+(T-t)\theta G^{-1}(x) f(x,u) \cdot \nabla G(x)
	\geq -C -C\theta.
\]
By \eqref{g_est}, we have
\[
	 \frac{\partial \phi}{\partial t}
	+F_{\theta}(x, \nabla \phi(t,x), D^2 \phi(t,x)) 
	\geq G(x)^\theta (-1-C-C\theta +e^{\alpha \theta}).
\]
Therefore, there exists $\tilde{\theta}(\alpha)>0$ such that
for $\theta \geq \tilde{\theta}(\alpha)$
\[
 \frac{\partial \phi}{\partial t}
	+F_{\theta}(x, \nabla \phi(t,x), D^2 \phi(t,x)) 
	>  0, \ (t,x) \in (0,T) \times \mathbb{R}^n,
\]
that is, $\phi(t,x)$ is a strict classical subsolution.
Since $\Psi_\theta(t,x)$ is a viscosity (super) solution and $\Psi_\theta(T,0)=\phi(T,0)=0$, 
the argument of the classical comparison theorem for parabolic type equations 
works by using $\phi(t,x)$ as a test function.
Thus, we have 
\[
	\phi(t,x) \leq \Psi_\theta(t,x), \ (t,x) \in (0,T) \times \mathbb{R}^n.
\]
(Also, see \cite[Theorem V.9.1]{FS06}).
Hence, for $\theta \geq\tilde{\theta}(\alpha)$,
we obtain the lower estimate for $V_\theta(t,x)$
\[
	\frac{1}{\theta}\log (T-t) +\underline{l}(x)-2\alpha
	\leq \frac{1}{\theta}\log (T-t) +g(x)
	\leq V_\theta(t,x).
\]
Taking the liminf on $\theta$ and $(t,x)$,  we have
\[
	\underline{l}(x)-2\alpha \leq \underline{V}(t,x).
\]
Sending $\alpha \to 0$,  we have the lower bound for $\underline{V}(t,x)$
\[
	\underline{l}(x) \leq \underline{V}(t,x).
\]

For the upper estimate in \eqref{est_V-bar}, let us consider
$J_\theta (t,x; u(\cdot))$ with constant control $u(s)=u$:
\begin{equation}
	J_\theta(t,x; u) = E_{tx}\left[ \int_t^T e^{\theta l(X(s), u)} ds \right]
	= \int_t^T  E_{tx} \left[  e^{\theta l(X(s), u)} \right]	 ds .
\label{J_const}
\end{equation}
Define $W_\theta^u (r,x)$ $((r,x) \in [t,s] \times \mathbb{R}^n$) by
\[
	W_{\theta}^u(r,x)= E_{rx} \left[  e^{\theta l(X(s),u)} \right].
\]
Then, $W^u_\theta(r,x)$ satisfies the following linear PDE of parabolic type
(in viscosity sense):
\begin{gather*}
	\frac{\partial W_\theta^u}{\partial r}
	+\frac{1}{2\theta}\tr (a(x,u)D^2 W^u_\theta(r,x))+f(x,u) \cdot \nabla W^u_\theta (r,x)=0,
	\ (r,x) \in (t,s) \times \mathbb{R}^n, \\
	W^u_{\theta}(s,x)=e^{\theta l(x,u)}, x \in \mathbb{R}^n.
\end{gather*}
Taking the transformation
\[
	Z^u_\theta(r,x)=\frac{1}{\theta}\log W_{\theta}^u(r,x),
\]
we have the nonlinear PDE for $Z^u_\theta(r,x)$
\begin{gather*}
	\frac{\partial Z^u_\theta}{\partial r}+\frac{1}{2\theta}\tr(a(x,u)D^2 Z^u_\theta(r,x))
	+\frac{1}{2}a(x,u) \nabla Z^u_\theta(r,x) \cdot \nabla Z^u_\theta (r,x) \\
	\qquad \qquad \qquad \qquad \qquad \qquad \qquad \qquad 
	+f(x,u) \cdot \nabla Z^u_\theta 
	(r,x)=0, \ (r,x) \in (t,s) \times \mathbb{R}^n, \\
	Z^u_\theta(s,x)=l(x,u), \ x \in \mathbb{R}^n.
\end{gather*}

Take $\alpha>0$. Since $x \mapsto l(x,u)$ is uniformly Lipschitz on $u$,
there exists a smooth function $h^u (x)$ such that
\begin{gather}
	|h^u(x) -l(x,u)|  \leq \alpha, \ x \in \mathbb{R}^n, \label{approx_l(x,u)} \\
	|D_i h^u(x) | \leq C_1, \ |D_{ij}h^u(x)| \leq C_2(\alpha) \notag
\end{gather}
for some constants $C_1$, $C_2(\alpha)$. 
Note that we can take $C_1$ which does not depend on $\alpha$ and $u$.
Let us consider $\tilde{Z}^u_\theta(r,x)$ defined by
\begin{equation}
	\tilde{Z}^u_\theta(r,x)=Z^u_\theta (r,x)-h^u(x).
\label{def_Z^u_theta}
\end{equation}
After some calculations, we can see that $\tilde{Z}^u_\theta(r,x)$ satisfies
\begin{gather}
	\frac{\partial \tilde{Z}^u_\theta}{\partial r}(r,x)
	+\frac{1}{2\theta}\tr (a^u(x) D^2 \tilde{Z}^u_\theta(r,x))
	+\frac{1}{2}a^u(x)\nabla \tilde{Z}^u_\theta(r,x) \cdot \nabla \tilde{Z}^u_\theta(r,x) 
	\qquad \qquad 
	\notag \\
	\qquad \qquad \qquad \qquad
	+F^u(x) \cdot \nabla \tilde{Z}^u_\theta(r,x)
	+L^u_\theta (x)=0, \ (r,x) \in (t,s) \times \mathbb{R}^n, 
	\label{aug_DPE} \\
	 \tilde{Z}^u_\theta (s,x) =l(x,u)-h^u(x), \ x \in \mathbb{R}^n.
\end{gather}
where 
$a^u(x)=a(x,u)$,
$ F^u (x)=f(x,u)+a(x,u)\nabla h^u(x)$ and 
\begin{equation}
	L^u_\theta (x)=\frac{1}{2\theta}\tr(a(x,u)D^2 h^u(x))
	+\frac{1}{2}a(x,u) \nabla h^u(x) \cdot \nabla h^u(x)
	+f(x,u) \cdot \nabla h^u(x),
	\label{est-L^u_theta}
\end{equation}
Note that  we can take large $M>0$ (independent of $u$ and $\alpha$)
and $\hat{\theta}(\alpha)>0$ such that
\[
	|L^u_\theta (x)| < M, \ x \in \mathbb{R}^n, \ u \in U, \ \theta \geq \hat{\theta}(\alpha).
\]
Then, $M(s-r)+\alpha$ is a strict classical supersolution of \eqref{aug_DPE}.
Thus, for $\theta \geq \hat{\theta}(\alpha)$, we have
\[
	\tilde{Z}^u_\theta (r,x) \leq M (s-r) + \alpha.
\]
In particular, if we set $r=t$, we obtain
\begin{align*} 
	Z^u_\theta (t,x)=\frac{1}{\theta}\log E_{tx}[e^{\theta l(X(s), u)} ] 
	& \leq h^u(x) + M(s-t)+ \alpha \\
	& \leq h^u(x)+M(T-t)+\alpha.
\end{align*}
Therefore, we have the estimate
\[
	E_{tx}[e^{\theta l(X(s),u)} ] \leq e^{\theta(h^u(x) + M(T-t)+ \alpha)}.
\]
Hence the following upper bound for $J_\theta (t,x;u)$ holds:
\[
	J_\theta (t,x;u) \leq (T-t) e^{\theta(h^u(x) + M(T-t)+ \alpha)}.
\]

From the definition of $V_\theta(t,x)$ and \eqref{approx_l(x,u)},
\[
	V_\theta(t,x) \leq \frac{1}{\theta} \log J_\theta(t,x;u) 
	\leq \frac{1}{\theta}\log (T-t) +M(T-t) + l(x,u)+ 2\alpha.
\]
for $\theta \geq \hat{\theta}(\alpha)$.
Taking the limsup on $\theta$ and $(t,x)$,
\[
	\overline{V}(t,x) \leq l(x,u) + M(T-t) + 2\alpha.
\]
Then, taking the limit as $\alpha \to 0$, we have
\[
	\overline{V}(t,x) \leq l(x,u) +M(T-t).
\]
Since $M>0$ does not depend on $u$, we have the upper estimate in \eqref{est_V-bar}. $\qed$

Under Corollary \ref{est_terminal} and the comparison theorem for the limit equation,
the argument is standard to identify the semi-continuous limits.
We shall give the proof for convenience.
\begin{thm} \label{asym_r-s_value}
$V_\theta(t,x)$ converges to $V(t,x)$ as $\theta \to \infty$
uniformly on each compact set in $(0,T) \times \mathbb{R}^n$.
\end{thm}
\noindent
\textit{Proof.} \
We shall show $\overline{V}(t,x)$ is a viscosity subsolution of \eqref{u_Isaacs}.
Let $\varphi(t,x)$ be a smooth function on $(0,T) \times \mathbb{R}^n$ and
$(\hat{t},\hat{x})$ be a maximum point
of $\overline{V}(t,x)-\varphi(t,x)$ in $(0,T) \times \mathbb{R}^n$.
We may assume $(\hat{t},\hat{x})$ is a strict maximum point and
$\overline{V}(\hat{t},\hat{x})=\varphi(\hat{t},\hat{x})$.
Since $(\hat{t},\hat{x})$ is a strict maximum,
we can take some sequences $\{ \theta_n \}$, $\{ (t_n, x_n) \}$
($\theta_n \to \infty$, $(t_n,x_n) \to (\hat{t},\hat{x})$) such that
$V_{\theta_n}(t_n, x_n) \to \overline{V}(\hat{t},\hat{x})$ and
$(t_n, x_n)$ is a local maximum point of $V_\theta (t,x) - \varphi(t,x)$.
Noting that $V_\theta(t,x)$ is a viscosity (sub) solution of \eqref{DPE_V_theta},
\begin{equation}
	0 \leq \frac{\partial \varphi}{\partial t}(t_n,x_n)
	+\mathcal{H}_{\theta_n}(\hat{x}, \varphi(t_n,x_n), \nabla \varphi(t_n,x_n)),
			D^2 \varphi(t_n,x_n)).
\label{sub-V_theta}
\end{equation}
For $\alpha>0$, let $A_\alpha=A(\hat{x},\overline{V}(\hat{t},\hat{x})-\alpha)$.
From \eqref{sub-V_theta}, we have
\begin{multline}
	0 \leq \frac{\partial \varphi}{\partial t}(t_n, x_n)
	+\min_{u \in A_{\alpha}}
	\bigg\{ 
	\frac{1}{2\theta_n}\tr(a(x_n,u)D^2\varphi(t_n,x_n)) \\
	+H^u(x_n,\nabla \varphi(t_n,x_n))+ e^{\theta_n (l(x_n,u)-\varphi(t_n,x_n))} \bigg\}.
\label{sub-V_theta_app}
\end{multline}
If we take large $n$,
\begin{align*}
	l(x_n,u)-\varphi(t_n,x_n)
	&=l(x_n,u)-l(\hat{x},u) +l(\hat{x},u)-\overline{V}(\hat{t},\hat{x}) \\
	& \qquad \qquad 
   +\overline{V}(\hat{t},\hat{x})-\varphi(t_n,x_n) \\
	&\leq -\frac{\alpha}{2} \text{ for any } u \in A(\hat{x},\overline{V}(\hat{t},\hat{x})-\alpha).
\end{align*}
From \eqref{sub-V_theta_app}, 
\[
	0 \leq \frac{\partial \varphi}{\partial t}(t_n, x_n)
	+\min_{u \in A_{\alpha}}
	\bigg\{ 
	\frac{1}{2\theta_n}\tr(a(x_n,u)D^2\varphi(t_n,x_n))
	+H^u(x_n,\nabla \varphi(t_n,x_n)) \bigg\} + e^{-\theta_n \alpha /2}.
\]
Taking the limit as $n \to \infty$,  we have
\[
	0 \leq  \frac{\partial \varphi}{\partial t}(\hat{t},\hat{x})
	+\mathcal{H}(\hat{x},\overline{V}(\hat{t},\hat{x})-\alpha, \nabla \varphi(\hat{t},\hat{x})).
\]
By Lemma \ref{H-semi-lim}, if we send $\alpha \to 0$,
\[
	0 \leq  \frac{\partial \varphi}{\partial t}(\hat{t},\hat{x})
	+\mathcal{H}^\ast(\hat{x},\overline{V}(\hat{t},\hat{x}), \nabla \varphi(\hat{t},\hat{x})).
\]
Therefore $\overline{V}(t,x)$ is a viscosity subsolution of \eqref{u_Isaacs}.

Next, we shall prove that $\underline{V}(t,x)$ is a viscosity supersolution.
Let $\varphi(t,x)$ be a smooth function on $(0,T) \times \mathbb{R}^n$
and $(\hat{t},\hat{x})$ be a strict minimum point of $\underline{V}(t,x)-\varphi(t,x)$
on $(0,T) \times \mathbb{R}^n$ with $\underline{V}(\hat{t},\hat{x})=\varphi(\hat{t},\hat{x})$.
We take some sequences $\{ \theta_n \}$, $\{ (t_n, x_n) \}$
$(\theta_n \to \infty, \ (t_n,x_n) \to (\hat{t},\hat{x}))$ such that
$V_{\theta_n}(t_n,x_n) \to \underline{V}(\hat{t},\hat{x})$ and
$(t_n, x_n)$ is a  local minimum of $V_{\theta_n}(t,x) -\varphi(t,x)$.
Since $V_\theta(t,x)$ is a viscosity (super) solution of \eqref{u_Isaacs},
\begin{multline}
	\frac{\partial \varphi}{\partial t}(t_n, x_n)
	+\min_{u \in U}
	\bigg\{ 
	\frac{1}{2\theta_n}\tr(a(x_n,u)D^2\varphi(t_n,x_n)) \\
	+H^u(x_n,\nabla \varphi(t_n,x_n))+ e^{\theta_n (l(x_n,u)-\varphi(t_n,x_n))} \bigg\}
	\leq 0.
\label{super-V_theta}
\end{multline}

Take $\alpha>0$.  We note that for large $n$, the minimum of \eqref{super-V_theta} is the same
if we replace $U$ with
\[
	C^n_\alpha =\{ u \in U \,; \, l(x_n,u)-V_{\theta_n}(t_n, x_n) \leq \alpha/2 \}
\]
Let $u \in C_\alpha^n$. Then, we have
\begin{align*}
	&l(\hat{x},u)-\underline{V}(\hat{t},\hat{x}) \\
	&\leq l(\hat{x},u)-l(x_n,u)+ l(x_n,u) -V_{\theta_n}(t_n,x_n) 
	+V_{\theta_n}(t_n, x_n)-\underline{V}(\hat{t},\hat{x})
	\leq \alpha
\end{align*}
for large $n$ uniform on $u$. 
Thus, for large $n$, we see that
$ C_\alpha^n \subset B_\alpha =A(\hat{x}, \underline{V}(\hat{t},\hat{x})+\alpha)$.
Therefore we can obtain
\begin{multline*}
	\frac{\partial \varphi}{\partial t}(t_n, x_n)
	+\min_{u \in B_\alpha}
	\bigg\{ 
	\frac{1}{2\theta_n}\tr(a(x_n,u)D^2\varphi(t_n,x_n)) \\
	+H^u(x_n,\nabla \varphi(t_n,x_n))+ e^{\theta_n (l(x_n,u)-\varphi(t_n,x_n))} \bigg\}
	\leq 0.
\end{multline*}
for large $n$. Since  $e^{\theta_n (l(x_n,u)-\varphi(t_n,x_n))}>0$,
\[
\frac{\partial \varphi}{\partial t}(t_n, x_n)
	+\min_{u \in B_\alpha}
	\bigg\{ 
	\frac{1}{2\theta_n}\tr(a(x_n,u)D^2\varphi(t_n,x_n)) \\
	+H^u(x_n,\nabla \varphi(t_n,x_n)) \bigg\}
	\leq 0
\]
If we take the limit as $n \to \infty$, 
\[
	\frac{\partial \varphi}{\partial t}(\hat{t}, \hat{x})
	+\mathcal{H}(\hat{x}, \underline{V}(\hat{t},\hat{x})+\alpha, \nabla \varphi(\hat{t},\hat{x})) \leq 0.
\]
By Lemma \ref{H-semi-lim}, if we take $\alpha \to 0$, 
\[
	\frac{\partial \varphi}{\partial t}(\hat{t}, \hat{x})
	+\mathcal{H}_\ast(\hat{x}, \underline{V}(\hat{t},\hat{x}), \nabla \varphi(\hat{t},\hat{x})) \leq 0.
\]
Thus, $\underline{V}(t,x)$ is a viscosity supersolution of \eqref{u_Isaacs}.

From the above results and Corollary \ref{est_terminal}, we proved that
$\overline{V}(t,x)$ (\textit{resp.\,}$\underline{V}(t,x)$)
is a bounded upper (\textit{resp.\,}lower) semi-continuous viscosity subsolution 
(\textit{resp.\,}supersolution)
satisfying the terminal condition
\[
	\overline{V}(T,x)=\underline{V}(T,x)=\min_{u \in U} l(x,u), \ x \in \mathbb{R}^n.
\]
We recall that  $V(t,x)$ is  a bounded Lipschitz continuous viscosity solution
of \eqref{u_Isaacs} with the terminal condition $\min_{u \in U}l(x,u)$
from Theorem \ref{ch_value}.
Thus, by applying Theorem \ref{comp}, we have
\[
	\overline{V}(t,x) \leq V(t,x) \leq \underline{V}(t,x), \ (t,x) \in (0,T) \times \mathbb{R}^n.
\]
Since $\underline{V}(t,x) \leq \overline{V}(t,x)$ is trivial,
we obtain
\[
	\overline{V}(t,x)=\underline{V}(t,x)=V(t,x), \ (t,x) \in (0,T) \times \mathbb{R}^n.
\]
Hence $V_\theta(t,x)$ converges to $V(t,x)$ as $\theta \to \infty$
uniformly on compact sets. $\qed$


\section{Mathematical finance example}

In this section we consider the classical Merton optimal consumption problem in mathematical finance, on a finite time horizon. This problem has an explicit solution, and the kind of risk sensitive limit considered in Section 6 can be found by direct calculations.

In the Merton problem, let $X(s)$ denote an investor's wealth at time $s$. The wealth is divided between a risky and a riskless asset. Let $k(s)$ be the fraction of wealth in the risky asset and $C(s)$ the consumption rate. Wealth obeys the SDE
\begin{equation}
	\left\{
	\begin{aligned}
		dX(s) &= X(s) [ r(1-k(s))ds +k(s) (\mu ds +\Sigma dW(s))] -C(s) ds, \ t \leq s \leq T, \\
		X(t) &= x >0
	\end{aligned}	
	\right.
\label{stoch_w}
\end{equation}
with $r$ the riskless interest rate and  $\mu, \Sigma$ the mean return rate and volatility for the risky asset. For the Merton problem on a finite time interval, the goal is to minimize
$$
	E_{tx} \left [ \int^T_tC(s)^{-\theta}ds\right ], \ \theta > 0
$$
The ``HARA parameter'' is $-\theta$ and the discount rate is 0. Let $C(s) = c(s)X(s)$. The control in this stochastic control problem is $u(s)= (k(s),c(s))$, and $X(s)$ is the state. There is the control constraint $c(s) > 0$ and the state constraint $X(s) > 0$. We require that $k(s)$ is bounded and that $c(s)$ is bounded on $[t,T_1]$ for any $T_1 < T$.

As in \eqref{stoch_sys}, let $\Psi_\theta(t,x)$ denote the value function for this stochastic control problem. By using the dynamic programming PDE \eqref{DPE_Psi} there is the explicit solution
\cite[p.\ 161]{FR75}
\begin{gather}
	\Psi_\theta(t,x) = [h_\theta(t)]^{1+\theta}x^{-\theta},
		\label{Psi_th} \\
	h_\theta(t) = \frac{1+\theta}{\nu_\theta \theta} \left (1-e^{-\frac{\nu_\theta \theta}{1+\theta}(T-t)}\right ),
		\label{h_th} \\
	\nu_\theta = \frac{(\mu-r)^2}{2\Sigma^2(1+\theta)} + r.
		\label{nu}
\end{gather}
The optimal controls are
\begin{equation}
	k^*_\theta(s) = \frac{\mu-r}{\Sigma^2(1+\theta)}, \ c^*_\theta(s) = [h_\theta(s)]^{-1}.
\label{opt_ctr}	
\end{equation}
Note that $k^*_\theta(s)$ is constant and that $c_\theta^*(s)$ depends only on $s$ and not on $X(s)$. Since $h_\theta(T) = 0$, $c^*_\theta(s) \rightarrow \infty$ as $s\rightarrow T^-$.

The HARA parameter $-\theta$ is a measure of risk aversion, and $\theta \rightarrow \infty$ is a ``totally risk averse limit.'' We suppose that the volatility $\Sigma =\Sigma(\theta)$ is such that $\theta \Sigma^2(\theta)$ tends to $\bar{\sigma}^2$ as $\theta \rightarrow \infty$ $(\bar{\sigma} > 0)$. For the Merton problem, some of the assumptions in Section 6 are not satisfied. Instead of using Theorem \ref{r-s_value}, we obtain a solution to the limit max-plus control problem directly. 

Let $V_\theta = \theta^{-1} \log \Psi_\theta$. From \eqref{Psi_th}, \eqref{h_th}
\begin{gather}
	\lim_{\theta \rightarrow \infty} V_\theta (t,x) = V(t,x), \\
	V(t,x) = -\log x + B(t), \label{expl_value} \\
	B(t) = \log [\nu^{-1}(1-e^{-\nu(T-t)})], \label{time_part}
\end{gather}
where $\nu=(\mu-r)^2/(2\bar{\sigma}^2)+r$.
From \eqref{opt_ctr}, as $\theta \rightarrow \infty$, $k^*_\theta (s)$ and $c^*_\theta(s)$ tend to 
\begin{equation}
	k^*(s) = k^* = \frac{\mu-r}{\bar{\sigma}^2}, c^*(s) = e^{-B(s)}.
\label{lim_part}
\end{equation}
In the limiting max-plus control problem, the state $x(s) >0$ satisfies \eqref{system} with
\begin{align*}
	f(x,k,c) & =  x[r+(\mu-r)k-c]\\
	\sigma(x,k)& =  \bar{\sigma}kx.
\end{align*}
From \eqref{H^u}
\begin{equation} 
	H^u(x,V_x) = (r-c) xV_x + (\mu-r) kxV_x + \frac{1}{2} \bar{\sigma}^2k^2x^2V^2_x.
\end{equation}
For $V(t,x)$ of the form \eqref{expl_value}
\begin{equation}
	\min_k H^u(x,V_x) = c-\nu.
\label{min_k}
\end{equation}
The minimum in \eqref{min_k} occurs at the constant $k^*$ in \eqref{lim_part}. For the Merton problem
$$l(x,c) = -\log x-\log c.$$
For $V$ as in \eqref{expl_value} the dynamic programming QVI \eqref{QVI} becomes 
\begin{equation}
	0 = \min_c \{ \max [-\log c-B(t), \dot{B} (t) + c -\nu]\}.
\end{equation}
At the minimum $c=c^*(t)$,
\begin{align}
	0  = & -\log c^*(t) - B(t), \label{QVI_min1} \\
	0  = & \dot{B} (t) + c^*(t) - \nu. \label{QVI_min2}
\end{align}
This agrees with \eqref{time_part} and \eqref{lim_part}. Since $B(T) = - \infty$, $c^*(t) \rightarrow \infty$ as $t\rightarrow T^-$.

For the Merton problem, the optimal controls $k^*_\theta(s)$, $c^*_\theta(s)$ in \eqref{opt_ctr} are functions of time only. This suggests (but does not prove) that the minimum among strategies $\alpha$ in the max-plus control problem is the same as the minimum among controls $u(s) = (k(s), c(s))$ which depend only on time $s$. Let us therefore consider only $\alpha^u$ such that $\alpha^u[v] (s) = u(s)$ for all disturbances $v(\cdot)$. We assume that $u(s)$ is bounded on $[t,T]$ and right continuous with left limits, as in (S1) of Section 3. It is easy to show that
\begin{equation}
	E^+\left [ -\int^s_t\bar{\sigma}k(\rho) v(\rho) d\rho\right ]
	= \frac{\bar{\sigma}^{2}}{2} \int^s_t k(\rho)^2 d\rho.
\end{equation}

From \eqref{tower} we then have
\begin{eqnarray*}
J(t,x; \alpha^u) & = & -\log x + \tilde{J}(t,u),\\
\\
\tilde{J}(t,u) & = & \max_{s \in [t,T]}
\left\{ \int^s_t(c(\rho)-r) d\rho + \int^s_t \left [ \frac{\bar{\sigma}^2}{2} k(\rho)^2 - (\mu - r)k (\rho)\right ] d\rho - \log c (s) \right\}.
\end{eqnarray*}
The second integrand is minimized by taking $k(\rho)=k^*$, with $k^*$ as in \eqref{lim_part}. Hence, $k^*$ is the optimal fraction of wealth in the risky asset. With this choice of $k(s)$,
$$\tilde{J} (t,k^*,c) = \max_{s\in [t,T]} \left [ \int^s_t (c(\rho) - \nu)d\rho - \log c(s)\right ]$$
with $\nu$ as in \eqref{nu}. Let us show that
\begin{equation}
	\tilde{J}(t,k^*,c) \geq B(t).
\label{time_l_bd}
\end{equation}
By \eqref{QVI_min1} and \eqref{QVI_min2}
\begin{equation}
	\int^s_t (c(\rho)-\nu)d\rho - \log c(s)
	=B(t) - \log c(s) + \log c^*(s) + \int^s_t (c(\rho) - c^*(\rho)) d\rho . 
\label{comp_c}
\end{equation}
If $c(t) \leq c^*(t)$, then \eqref{time_l_bd} clearly holds. Suppose that $c^* (t)< c(t)$ and let
$$s_0 = \inf \{ s> t\colon c(s) \leq c^*(s) \}.$$
Since $c(s)$ is bounded and $c^*(s) \rightarrow \infty$ as $s\rightarrow T^-$, we have $s_0< T$. Since $c(\cdot)$ is right continuous, $c(s_0) = c^*(s_0)$. Then
$$\tilde{J} (t,k^*, c) \geq B(t) + \int^{s_0}_t (c(\rho) - c^*(\rho)) d\rho > B(t),$$
which proves \eqref{time_l_bd}. If $c(s) = c^*(s)$ for all $s$, then the right side of \eqref{comp_c} equals $B(t)$. Hence, $\tilde{J}(t,k^*c^*) = B(t)$. Since $c^*(s) \rightarrow \infty$ as $s\rightarrow T^-$, $u^* = (k^*,c^*)$ is not an admissible control. However, for $\delta > 0$ let $c_\delta(s) = 1$ if $s\in [t,t+\delta ]$ and $c_\delta (s) = c^*(s-\delta)$ if $s\in [t+ \delta, T]$. For $s\in [t+\delta,T]$
$$\int^s_{t+\delta} (c_\delta(\rho) - \nu) d\rho - \log c_\delta (s) = \int^{s-\delta}_t (c^*(\rho)-\nu) d\rho - \log c^*(s-\delta) = B(t) .$$
Hence
$$B(t) \leq \tilde{J}(t,k^*,c_\delta)\leq |1-\nu|\delta + B(t).$$
Since $\delta$ is arbitrarily small, $B(t)$ is the infimum of $\tilde{J}(t,k,c)$ among all admissible controls $u=(k,c)$.

\begin{rem}
The Verification Theorem \ref{ver_thm} cannot be applied directly to this example for two reasons.
There is a state constraint $x(s)>0$ in the Merton problem.
Also, the value function $V(t,x)$ in \eqref{expl_value} is unbounded as $t \to T-$.
However, Theorem \ref{ver_thm} can be applied to the following modified version of the Merton problem.
We take $y(s)=\log x(s)$ as the state at time $s$.
Also we require that $0<c(s) \leq C <\infty$, where $C>\nu$ is a fixed constant.
The corresponding max-plus value function $\tilde{V}(t,y)$ has the form
\[
	\tilde{V}(t,y)=-y+\tilde{B}(t).
\]
By a calculation similar to \eqref{lim_part}--\eqref{QVI_min2}, the optimal control $(k^\ast,\tilde{c}^\ast(s))$
satisfies \eqref{lim_part} with $B(s)$ replaced by $\tilde{B}(s)$.
Moreover, $\tilde{B}(s)$, $\tilde{c}^\ast(s)$ satisfy \eqref{QVI_min1}--\eqref{QVI_min2}
with $\tilde{c}^\ast(T)=C$.
Theorem \ref{ver_thm} implies that $(k^\ast,\tilde{c}^\ast(s))$ is optimal compared to any Lipschitz
Markov control policy $(\underline{k}(s,y), \underline{c}(s,y))$ such that $0 <\underline{c}(s,y) \leq C$.
\end{rem}


\section{Infinite time horizon}
In this section we consider initial time $t=0$ and final time $T$ finite but arbitrarily large. We are concerned with inequalities of the form
\begin{equation}
	E^+_{0x} \left[ \int^\oplus_{[0,T]} l (x(s), u(s)) ds \right] \leq W(x)
\label{dis_ineq}
\end{equation}
which are required to hold for every $T> 0$. Such inequalities are of interest in nonlinear $H$-infinity control, and are related to results in \cite{DMc04}, \cite{F04}.

Let $V(0,x;T)$ denote the value function, where in \eqref{value_u} we write $V(t,x;T)$ to emphasize dependence on $T$. If a strategy $\alpha$ can be chosen such that \eqref{dis_ineq} holds for all $T$, with $u(s) = \alpha [v](s)$, then \eqref{dis_ineq} implies
\begin{equation}
	V(0,x;T) \leq W(x).
\end{equation}
Since $V(0,x;T)$ is a nondecreasing function of $T$, it tends to a limit $V_\infty(x)$ as $T\rightarrow \infty$, and
\begin{equation}
	V_\infty(x) \leq W(x).
\label{dis_lim}
\end{equation}

We will show later, under the additional assumption  (A5)  that $V_\infty$ is a Lipschitz continuous viscosity solution to the steady state form of \eqref{QVI}. However, we first consider the following elementary result in which $\alpha [v](s) = \underline{u}(x(s))$ with $\underline{u}$ a stationary control policy. For this result we assume:
\renewcommand{\labelenumi}{(A\arabic{enumi})}
\renewcommand{\labelenumii}{(\roman{enumii})}
\begin{enumerate}
\setcounter{enumi}{3}
\item
	\begin{enumerate}
	\item $f(x,u)$ and $\sigma(x,u)$ are Lipschitz continuous on $\mathbb{R}^n \times U$
		 and $\sigma(x,u)$ is bounded;
	\item $l(x,u)$ is continuous on $\mathbb{R}^n \times U$.
\end{enumerate}
\end{enumerate}
\renewcommand{\labelenumi}{(\roman{enumi})}
\renewcommand{\labelenumii}{(\alph{enumii})}

\begin{prop} \label{super_opt_infty}
Assume \textup{(A4)}. Suppose that $W(x)$ is $C^1$, $\underline{u}(x)$ is Lipschitz on
$\mathbb{R}^n$ and that for any $y\in \mathbb{R}^n$
\begin{equation}
	\max \{ H^{\underline{u}(y)}(y,\nabla W(y)),l (y, \underline{u}(y)) - W(y) \} \leq 0.
\label{infty_super}
\end{equation}
Then \eqref{dis_ineq} holds with $u(s) = \underline{u}(x(s))$.
\end{prop}

\noindent
\textit{Proof}. \ 
For any $v\in L^2 [t,T]$ the ODE \eqref{system} has a unique solution $x(s)$ on $[t,T]$ with initial data $x(t) = x$. By the Fundamental Theorem of Calculus
\begin{align*}
	W(x(s)) - \frac{1}{2}\int^s_0 |v(\rho)|^2 d \rho 
	&= W(x) + \int^s_0 H^{\underline{u}(x(\rho))} (x(\rho), \nabla W(x(\rho)))d\rho \\
	&\quad \quad		-\frac{1}{2} \int^s_0
		|v(\rho)-\sigma^T(x(\rho), \underline{u} (x(\rho)))\nabla W(x(\rho))|^2 d\rho \\
	&\leq W(x).
\end{align*}
Since $l (x(s), \underline{u}(x(s))\leq W(x(s))$, we obtain \eqref{dis_ineq}. $\qed$

We note that \eqref{infty_super} implies that $W(x)$ is a viscosity supersolution of equation \eqref{infty_QVI} below, which is the steady state form of the QVI \eqref{QVI}.

Let us illustrate the use of Proposition \ref{super_opt_infty} in nonlinear $H$-infinity control theory. Suppose that a ``running cost'' $l_1(x, u)$, a ``terminal cost'' function $G(x)$ and a parameter $\mu > 0$ are given. Consider control policies $\underline{u}\colon \mathbb{R}^n \rightarrow U$. We seek a policy $\underline{u}$ such that
\begin{equation}
	\mu \left [ \int^T_0 l_1 (x(s), \underline{u} (x(s))ds + G(x(T))\right ]
	\leq W(x) + \frac{1}{2} \int^T_0 |v(s)|^2 ds
\label{dis_finite}
\end{equation}
for every initial state $x(0) = x$ and every $T> 0$. This formulation is considered in \cite{DMc04}. In much of the  $H$-infinity control literature, the terminal cost term $G(x(T))$ is omitted. Inequality \eqref{dis_finite} is often rewritten by multiplying each side by $\gamma^2 = \mu^{-1}$, where $\gamma$ is an $H$-infinity norm parameter. We may expect \eqref{dis_lim} to hold only under some further assumptions on $f$, $\sigma$, $l_1$  $G$, and for $\mu$ in some interval $[0, \mu_1]$.

To rewrite \eqref{dis_finite} in the form \eqref{dis_ineq}, we consider an augmented state $\hat{x}(s) = (x(s), x_{n+1}(s))$, where
$$
	\frac{dx_{n+1}(s)}{ds} = l_1 (x(s), \underline{u}(x(s))).
$$
Let $l (\hat{x},u) = \mu (x_{n+1} + G(x))$. Suppose that for every $s\geq 0$,
\begin{equation}
	E^+_{0\hat{x}}[l (\hat{x}(s), \underline{u}(x(s)))] \leq x_{n+1} + W(x).
\label{dis_ext}
\end{equation}
When we take $x_{n+1} = x_{n+1} (0) = 0$, this implies \eqref{dis_finite}. 

To apply Proposition \ref{super_opt_infty}, with $W(x)$ replaced by $\hat{W}(\hat{x}) = x_{n+1} + W(x)$, we require that for all $\hat{y}$
\begin{gather}
	H^{\underline{u}(y)} (y,\nabla_yW(y)) + \mu l_1(y,\underline{u}(y)) \leq 0, \label{super_feed}\\
	\mu (y_{n+1} + G(y)) \leq y_{n+1} + W(y).
\label{ext_sys}
\end{gather}
The following example illustrates the use of Proposition \ref{super_opt_infty} in nonlinear $H$-infinity control.

\begin{ex}
Assume that $f(0,0)=0$, $\underline{u} (0) = 0$ and $f(x,\underline{u}(x)) \cdot x\leq -c|x|^2$ for some $c> 0$. This implies that $x(s)$ is exponentially stable to the equilibrium point 0 if there is no disturbance $(v(s) = 0)$. Also assume that
$$0\leq l_1 (x,\underline{u}(x)) \leq C_1 |x|^2, \; 0 \leq G(x) \leq C_2 |x|^2$$
for some $C_1$, $C_2$. We choose $W(x) = K|x|^2$. An easy calculation shows that \eqref{super_feed} holds if $K\|a \| < c $ and $\mu$ is small enough $(a = \sigma \sigma^T)$. Inequality \eqref{ext_sys} holds if $\mu < 1$, $\mu C_2 < K$, $C_1\mu +2K^2 \| a \|^2 -Kc \leq 0$.
\end{ex}

To conclude this section, let us return to consider the limit function $V_\infty(x)$ in \eqref{dis_lim}. We now assume conditions (A1)--(A3) in Section 2 and also:
\renewcommand{\labelenumi}{(A\arabic{enumi})}
\renewcommand{\labelenumii}{(\roman{enumii})}
\begin{enumerate}
\setcounter{enumi}{4}
\item
\begin{enumerate}
	\item $\sigma=\sigma(u)$, 
	\item $(f(x,u)-f(y,u)) \cdot (x-y) \leq 0$ for all $x,y \in \mathbb{R}^n$, $u \in U$.
\end{enumerate}
\end{enumerate}
\renewcommand{\labelenumi}{(\roman{enumi})}
\renewcommand{\labelenumii}{(\alph{enumii})}

\begin{prop}
For every $x,y\in \mathbb{R}^n$, $T> 0$ and $\alpha \in \Gamma (0,T)$
\begin{equation}
	|J(0,x; T, \alpha) - J(0,y; T,\alpha)|
	\leq \| l_x \|  |x-y|,
\label{est_state}
\end{equation}
where $\| \cdot \|$ is the sup norm.
\end{prop}

\noindent
\textit{Proof}. \ 
Let $x(s)$, $y(s)$ be the solutions of \eqref{system} with $x(0) = x$, $y(0) = y$ and $u(s) = \alpha [v](s)$. Since $\sigma = \sigma (u)$,
$$\frac{d}{ds} |x(s) - y(s)|^2 = 2 (f(x(s),u(s)) - f(y(s), u(s))) \cdot (x(s) - y(s)) \leq 0.$$
Hence, $|x(s)-y(s)| \leq |x-y|$. Since
$$
l (x(s), u(s)) - l (y(s), u(s)) | \leq \; \| l_x \| |x(s) - y(s)|
$$
and $\alpha$ is arbitrary, this implies \eqref{est_state}. $\qed$

By \eqref{dis_lim}, $V_\infty (x)$ is Lipschitz with the same Lipschitz constant $\| l_x \|$ as in (8.9). It can be shown that $V_\infty$ is a viscosity solution of the steady state  form of \eqref{QVI}:
\begin{equation}
\min_{u\in U} \max \{H^u(x,\nabla V_\infty (x)), l (x,u) - V_\infty (x)\} = 0, \; x\in \mathbb{R}^n.
\label{infty_QVI}
\end{equation}


\appendix
\section{Proof of Lemma \ref{lem_approx_str}}

We first construct an approximation in a short time after the initial time $t$.
Set a constant $u^\eps_1 \in U$ by
\begin{equation}
	u^\eps_1 \equiv \alpha[v](t).
\label{min_const_1}
\end{equation}
Note that $u^\eps_1$ does not depend on $v \in L^2 ([t,T];\mathbb{R}^d)$
from Lemma \ref{inst_delay} (i). 
Although $u^\eps_1$ actually does not  depend on $\eps$, 
we use this notation because of the consistency.

We define $\hat{v}_1^\eps \in L^2[t,T]$ by
\[
	\hat{v}_1^\eps (s) \equiv \beta[u_1^\eps](s), \ s \in [t,T],
\]
where $u_1^\eps$ in $\beta[u_1^\eps]$ is understood as a constant curve on $[t,T]$
taking  the constant $u_1^\eps$.
Consider the first time when the error of $u_1^\eps$ and $\alpha[\hat{v}_1^\eps]$
gets larger than $\eps$, \textit{i.e.}, define $s_1 \in [t,T]$ by
\[
	s_1 \equiv
		\inf\{ s \in [t,T] \,; \, |u_1^\eps-\alpha[\hat{v}_1^\eps](s)| \geq \eps \} \wedge T.
\]
We understand that $\inf \emptyset = \infty$.
Here recall that $u_1^\eps=\alpha[\hat{v}_1^\eps](t)$
from \eqref{min_const_1} and Lemma \ref{inst_delay} (i).

By using the definition of ${\Gamma}(t,T)$,
we can see that
\begin{gather}
	|u_1^\eps -\alpha[\hat{v}_1^\eps](s)| < \eps, \ t \leq  s < s_1,
		\label{hat_bf_1} \\
	|u_1^\eps-\alpha[\hat{v}_1^\eps](s_1)| \geq \eps \text{ if } s_1<T.
	\label{hat_at_1}
\end{gather}
\eqref{hat_bf_1} is immediate from the definition of $s_1$.
In order to see \eqref{hat_at_1}, take a sequence $\{r_k \}_{k=1}^\infty$
such that
\begin{equation}
	r_k \downarrow s_1 \ (k \uparrow \infty) \text{ and }
	|u_1^\eps-\alpha[\hat{v}_1^\eps](r_k)| \geq \eps, \ k=1,2,\cdots.
\label{approx_ab_1}
\end{equation}
Since $s \mapsto \alpha[\hat{v}_1^\eps](s)$ is right-continuous, 
by taking the limit as $k \to \infty$ in \eqref{approx_ab_1}, we have \eqref{hat_at_1}.

Define $v^\eps_1 \in L^2[t,s_1)$ by
\[
	v^\eps_1 \equiv \hat{v}_1^\eps|_{[t,s_1)} =\beta[u_1^\eps]|_{[t,s_1)}.
\]
For  $v_2 \in L^2[s_1, T]$, we denote by $v^\eps_1 \cdot v_2$
the concatenation of $v^\eps_1$ and $v_2$:
\[
	v^\eps_1 \cdot v_2 (s) \equiv
	\begin{cases}
		v_1^\eps(s), & t \leq s < s_1, \\
		v_2(s), & s_1 \leq s \leq T.	
	\end{cases}
\]
From (ii) of the definition of ${\Gamma}(t,T)$,
the values of $\alpha[v_1^\eps \cdot v_2]$ on $[t,s_1]$ are uniquely determined
by $v_1^\eps$. 
More precisely, the following claim holds:
\begin{equation}
	\text{For any }v_2, \ \tilde{v}_2 \in L^2[s_1, T], \
	\alpha[v^\eps_1 \cdot v_2](s) = \alpha[v^\eps_1 \cdot \tilde{v}_2](s), \
	\forall s \in [t, s_1].
\label{indep_af_1}
\end{equation}
Taking this property into mind, we define $\alpha[v_1^\eps]: [t,s_1] \to U$ by
\[
	\alpha[v_1^\eps](s) \equiv \alpha[v_1^\eps \cdot v_2](s), \
	t \leq s \leq s_1.
\]
Then, \eqref{hat_bf_1} and \eqref{hat_at_1} are rewritten by
\begin{gather}
	|u_1^\eps -\alpha[v_1^\eps](s)| < \eps, \ t \leq  s < s_1,
		\label{bf_1} \\
	|u_1^\eps-\alpha[v_1^\eps](s_1)| \geq \eps \text{ if } s_1<T.
	\label{at_1}	
\end{gather}

If $s_1=T$, we stop the approximating procedure
and set $s_n=s_1$ for $n = 2,3,\cdots$.
If $s_1<T$,
we continue the procedure and
construct an approximation in a short time after $s_1$.
Let us define a constant $u_2^\eps$ by
\[
	u_2^\eps \equiv \alpha[v^\eps_1](s_1) =\alpha[v^\eps_1 \cdot v_2](s_1).
\]
Note that $u_2^\eps$ does not depend on $v_2 \in L^2[s_1, T]$.

Define $\hat{v}_2^\eps \in L^2[t,T]$ by
\[
	\hat{v}_2^\eps(s) \equiv \beta[u_1^\eps \cdot u_2^\eps](s), \
	t \leq s \leq T,
\]
where $u_1^\eps \cdot u_2^\eps$ is a piecewise constant curve on $[t,T]$
defined by the concatenation of constant curves taking
the constant $u_1^\eps$ on $[t,s_1)$ and the constant $u_2^\eps$ on $[s_1, T]$.
Consider the first time $s_2 \in [s_1, T]$
when the error of $u^\eps_2$ and $\alpha[\hat{v}_2^\eps]$ gets larger than $\eps$;
\begin{equation}
	s_2 \equiv \inf \{ s \in [s_1, T] \,;\, |u_2^\eps -\alpha[\hat{v}_2^\eps](s)| \geq \eps\}.
\label{1st_af_2}
\end{equation}
Since $u_1^\eps=u_1^\eps \cdot u_2^\eps $ on $[t,s_1)$,
$\beta[u_1^\eps]=\beta[u_1^\eps \cdot u_2^\eps]$ a.e.\,on $[t,s_1]$.
Thus, from the definitions of $v^\eps_1$ and $\hat{v}_2^\eps$,
\[
	v^\eps_1 =\hat{v}_2^\eps \text{ a.e.\,on } [t,s_1).
\]
Then, \eqref{1st_af_2} can be written as
\[
	s_2 = \inf \{ s \in [s_1, T] \,;\, 
	|u_2^\eps -\alpha[v^\eps_1 \cdot \hat{v}_2^\eps|_{[s_1,T]}](s)| \geq \eps\}.
\]
Note that $u_2^\eps =\alpha[v^\eps_1 \cdot \hat{v}^\eps_2|_{[s_1,T]}](s_1)$.
In a similar way to \eqref{hat_bf_1} and \eqref{hat_at_1}, we have
\begin{gather}
	|u_2^\eps-\alpha[v^\eps_1 \cdot \hat{v}_2^\eps|_{[s_1,T]}](s)|<\eps, \
	s_1 \leq  s <s_2, \label{hat_bf_2} \\
	|u_2^\eps -\alpha[v^\eps_1 \cdot \hat{v}_2^\eps|_{[s_1,T]}](s_2)| \geq \eps
	\text{ if } s_2 < T. \label{hat_at_2}
\end{gather}

Let us define $v_2^\eps \in L^2[s_1,s_2)$ by
\[
	v_2^\eps \equiv \hat{v}_2^\eps |_{[s_1,s_2)}=\beta[u^\eps_1 \cdot u^\eps_2]|_{[s_1,s_2)}.
\]
We note that from (ii) of the definition of $\Gamma(t,T)$,
$\alpha[v^{\eps}_1 \cdot v^\eps_2 \cdot v_3]|_{[t,s_2]}$ does not
depend on $v_3  \in L^2[s_2,T]$. 
So, we denote $\alpha[v^\eps_1 \cdot v^\eps_2]: [t,s_2] \to U$ by
\[
	\alpha[v^\eps_1 \cdot v^\eps_2](s) \equiv
		\alpha[v^\eps_1 \cdot v^\eps_2 \cdot v_3](s), \ t \leq s \leq s_2.
\]
From \eqref{hat_bf_2}, \eqref{hat_at_2} and the definition of $v^\eps_2$, we obtain
\begin{gather}
	|u_2^\eps-\alpha[v^\eps_1 \cdot v^\eps_2](s)|<\eps, \
	s_1 \leq  s <s_2, \label{bf_2} \\
	|u_2^\eps -\alpha[v^\eps_1 \cdot {v}_2^\eps](s_2)| \geq \eps
	\text{ if } s_2 < T. \label{at_2}
\end{gather}

If we continue this procedure, we have $\{ s_n \}_{n=1}^\infty$
$(t<s_1 \leq s_2 \leq \cdots \leq s_n \leq  \cdots \leq T)$,
$\{u_n^\eps \}_{n=1}^\infty \subset U$ and
$v_n^\eps \in L^2[s_{n-1}, s_n)$ such that
for each $n=1,2,\cdots$,
\begin{gather}
	u^\eps_n=\alpha[v^\eps_1 \cdot v^\eps_2 \cdot \cdots \cdot v^\eps_{n-1}](s_{n-1}),
	\label{def_const_n} \\
	|u^\eps_n -\alpha[v^\eps_1 \cdot v^\eps_2 \cdot \cdots \cdot v^\eps_n](s)| < \eps,
	\	s_{n-1} \leq  s < s_n, 
	\label{bf_n} \\
	|u^\eps_n -\alpha[v^\eps_1 \cdot v^\eps_2 \cdot \cdots \cdot v^\eps_n](s_n)|
	 \geq \eps \text{ if } s_n <T, 
	\label{at_n} \\
	v_n^\eps = \beta[u^\eps_1 \cdot u^\eps_2 \cdot \cdots \cdot u^\eps_n]
	 \text{ on } [s_{n-1},s_n). \label{def_v_n}
\end{gather}
Here we make a convention that $s_{0}=t$
and $u^\eps_1\equiv \alpha[v](t)$ in \eqref{def_const_n} for $n=1$.
Note that the approximating procedure could stop in a finite step if $s_n=T$ for some $n$.

We shall show  $s_n \to T \ (n \to \infty)$.
On the contrary, suppose $\tau \equiv \lim_{n \to \infty}s_n < T$.
Define ${v}^\eps \in L^2[t,T]$ by
\begin{equation}
	{v}^\eps(s) \equiv
	\begin{cases}
		v^\eps_n(s), & s_{n-1} \leq s <s_n, \ n=1,2,\cdots \\
		v^0, & \tau \leq s \leq T,
	\end{cases}
\label{v^eps_tau}
\end{equation}
where $v^0 \in \mathbb{R}^d$ is a constant.
Since ${v}^\eps=v^\eps_1\cdot v^\eps_2 \cdot \cdots \cdot v^\eps_n$
a.e.\,on $[t,s_n]$,
\begin{equation}
	\alpha[v^\eps](s)=\alpha[v^\eps_1\cdot v^\eps_2 \cdot \cdots \cdot v^\eps_n](s),
	\ s_{n-1} \leq s \leq s_n.
\label{eq_on_n}
\end{equation}
In particular,
\begin{equation}
	\alpha[v^\eps](s_n)
	=\alpha[v^\eps_1\cdot v^\eps_2 \cdot \cdots \cdot v^\eps_n](s_n).
\label{eq_at_n}
\end{equation}
From \eqref{at_n}, we have
\[
	|u_n^\eps-\alpha[v^\eps](s_n)| \geq \eps.
\]
Noting that $u_n^\eps
=\alpha[v^\eps_1\cdot v^\eps_1 \cdot \cdots \cdot v^\eps_{n-1}](s_{n-1})
=\alpha[v^{\eps}](s_{n-1})$
 by \eqref{def_const_n} and \eqref{eq_at_n},
\[
	|\alpha[v^\eps](s_{n-1})-\alpha[v^\eps](s_n)| \geq \eps.
\]
Since $s \mapsto \alpha[v^\eps](s)$ has left limits on $(t,T]$,
by taking the limit as $n \to \infty$ in the above inequality, we have
\[
	0=|\alpha[v^\eps](\tau-) -\alpha[v^\eps](\tau-)| \geq \eps.
\]
This is a contradiction. Therefore  we have
\[
	T=\lim_{n \to \infty} s_n.
\]

We define $u^\eps : [t,T] \to U$ and $v^\eps: [t,T] \to \mathbb{R}^d$ as follows:
\begin{gather}
	u^\eps(s) \equiv
	\begin{cases}
		u^\eps_n, & s_{n-1} \leq s < s_n, \ n=1,2,\cdots, \\
		u^0, & s=T,
	\end{cases} \\
	v^\eps(s) \equiv
	\begin{cases}
		v^\eps_n(s), & s_{n-1} \leq s <s_n, \ n=1,2,\cdots \\
		v^0, & s = T,
	\end{cases} \label{v^eps}
\end{gather}
where $u^0$ and $v^0$ are  constants in $U$ and $\mathbb{R}^d$, respectively.
\eqref{v^eps} is exactly the same as \eqref{v^eps_tau} when $\tau=T$.

We shall check $u^\eps$ and $v^\eps$ are actually what we want.
As already noted in \eqref{eq_on_n},
\[
	\alpha[v^\eps](s) = \alpha[v^\eps_1\cdot v^\eps_2 \cdot \cdots \cdot v^\eps_n](s),
	\ s_{n-1} \leq s  \leq s_n.
\]
So, from \eqref{bf_n} and the definition of $u^\eps$, 
\[
	|u^\eps(s)-\alpha[v^\eps](s)| < \eps, \ s_{n-1} \leq s< s_n, \ n=1,2,\cdots.
\]
Therefore  we have
\[
	|u^\eps(s)-\alpha[v^\eps](s)|< \eps, \ t \leq s <T.
\]

Since $u^\eps=u^\eps_1 \cdot u^\eps_2 \cdot \cdots \cdot u^\eps_n$
on $[t,s_n)$,
\[
	\beta[u^\eps]=\beta[u^\eps_1 \cdot u^\eps_2 \cdot \cdots \cdot u^\eps_n]
	\text{ a.e.\,on } [t,s_n].
\]
Thus, from \eqref{def_v_n} and the definition of $v^\eps$, we can see
$$
	v^\eps=\beta[u^\eps] \text{ a.e.\,on } [t,T].
$$

\section{Generators for the case of Elliott-Kalton strategies} \label{EK_gen}
We note that the same results on $\mathcal{K}(x,r,p)$ as Lemmas \ref{H-semi-lim}
and \ref{mono_lim} hold:
\begin{gather*}
	\mathcal{K}_\ast(x,r,p) =\mathcal{K}(x,r+0,p)=\mathcal{K}(x,r,p), \\
	\mathcal{K}^\ast(x,r,p)=\mathcal{K}(x,r-0,p)
	=\max_{v \in \mathbb{R}^d} \inf_{u \in A'(x,r)}
		\left\{ (f(x,u)+\sigma(x,u)v) \cdot p -\frac{1}{2}|v|^2 \right\},
\end{gather*}
where $A'(x,r)=\{ u \in U \,;\, l(x,u)<r \}$.

We first show \eqref{EK-gen_u}.
If $A'(x,\varphi(t,x))=\emptyset$, \eqref{EK-gen_u} is trivial.
We consider the case where $A'(x,\varphi(t,x))\not=\emptyset$.
From the definition of $F_{t,t+\delta}^{EK}$,
\begin{multline*}
	F_{t,t+\delta}^{EK}\varphi(t+\delta,\cdot)(x)
	=\inf_{\alpha \in \Gamma_{EK}(t,t+\delta)} \sup_{v \in L^2[t,t+\delta]}
	\bigg\{  \int_{[t,t+\delta]}^\oplus l(x(s),\alpha[v](s)) ds  \oplus \varphi(t+\delta,x(t+\delta)) \\
	-\frac{1}{2}\int_{t}^{t+\delta}|v(s)|^2 ds \bigg\}.
\end{multline*}
Since $l$ and $\varphi$ are bounded, 
we may replace the range of the supremum with $v \in L^2[t,t+\delta]$
such that
\begin{equation}
	\int_t^{t+\delta} |v(s)|^2 ds \leq M
\label{dist_bound}
\end{equation}
for some constant $M$ independent of $\delta$ and $\alpha$.

Fix $\rho>0$ and take $\eps>0$ arbitrarily.
Note that for each $v \in \mathbb{R}^d$,
\begin{multline*}
	\max_{v \in \mathbb{R}^d}
	\min_{u \in A(x,\varphi(t,x)-\rho)}
	\left\{ (f(x,u)+\sigma(x,u)v) \cdot \nabla \varphi(t,x) -\frac{1}{2}|v|^2 \right\} \\
	\geq \min_{u \in A(x,\varphi(t,x)-\rho)}
	\left\{ (f(x,u)+\sigma(x,u)v) \cdot \nabla \varphi(t,x) -\frac{1}{2}|v|^2 \right\}.
\end{multline*}
We choose a measurable mapping $\bar{u}: \mathbb{R}^d \to U$ such that
$\bar{u}(v) \in A(x,\varphi(t,x)-\rho)$  and
\begin{multline}
	\min_{u \in A(x,\varphi(t,x)-\rho)}
	\left\{ (f(x,u)+\sigma(x,u)v) \cdot \nabla \varphi(t,x) -\frac{1}{2}|v|^2 \right\} \\
	> (f(x,\bar{u}(v))+\sigma(x,\bar{u}(v))v) \cdot \nabla \varphi(t,x) -\frac{1}{2}|v|^2
	-\eps 
\label{eps_str}
\end{multline}
for each $v \in \mathbb{R}^d$.

Define $\alpha \in \Gamma_{EK}(t,t+\delta)$ by
\[
	\alpha[v](s)=\bar{u}(v(s)), \ t \leq s \leq t+\delta.
\]
Let $x(s)$ be the solution of \eqref{mp-SDE} with $\| v \|_{L^2[t,t+\delta]}^2 \leq M$.
Then, we have
\begin{align*}
	&\varphi(t+\delta, x(t+\delta))-\varphi(t,x) \\
	&=\int_{t}^{t+\delta}
		\frac{\partial \varphi}{\partial s}(s,x(s))
		+(f(x(s),\alpha[v](s))+\sigma(x(s),\alpha[v](s))v(s)) \cdot \nabla \varphi(s,x(s)) ds.
\end{align*}
 Since we consider $v$ satisfying \eqref{dist_bound}, there exists $c(M)>0$ such that
 \begin{equation}
	|x(s)-x| \leq c(M) |s-t|^{1/2}, \ t \leq s \leq t+\delta.
\label{est_err_ini}
 \end{equation}
 Thus, we have
 \begin{multline}
	 \varphi(t+\delta,x(t+\delta)) -\varphi(t,x)\\
	 =\int_{t}^{t+\delta}
		\frac{\partial \varphi}{\partial t}(t,x)
		+(f(x,\alpha[v](s))+\sigma(x,\alpha[v](s))v(s)) \cdot \nabla \varphi(t,x)ds	
		+o(\delta)
\label{ini_est1}
 \end{multline}
 where $o(\delta)$ is uniform on $\| v \|_{L^2[t,t+\delta]}^2 \leq  M$.
By noting $\alpha[v](s)=\bar{u}(v(s))$, we have from \eqref{eps_str}, \eqref{ini_est1}
\begin{align}
	& \varphi(t+\delta,x(t+\delta)) -\frac{1}{2}\int_t^{t+\delta}|v(s)|^2 ds-\varphi(t,x) \notag \\
	&=\int_{t}^{t+\delta}
		\frac{\partial \varphi}{\partial t}(t,x)
		+(f(x,\alpha[v](s))+\sigma(x,\alpha[v](s))v(s)) \cdot \nabla \varphi(t,x)
		-\frac{1}{2}|v(s)|^2 ds	
		+o(\delta)  \notag \\
	&\leq \left( \frac{\partial \varphi}{\partial t}(t,x)
		+\mathcal{K}(x,\varphi(t,x)-\rho,\nabla \varphi(t,x))  \right) \delta
		+\eps \delta + o(\delta).
	\label{F-est}
\end{align}
On the other hand, recall $\bar{u}(v) \in A(x,\varphi(t,x)-\rho)$. 
Thus we have
\[
	l(x, \alpha[v](s)) \leq \varphi(t,x) -\rho, \ t \leq s \leq t+\delta.
\]
By \eqref{est_err_ini}, we see that for small $\delta$,
\[
	l(x(s),\alpha[v](s)) < \varphi(t+\delta, x(t+\delta))-\rho/2, \  t \leq s \leq t+\delta
\]
for any $\| v \|_{L^2[t,t+\delta]}^2 \leq M$.
Then, for small $\delta$,
\begin{align}
	&F^{\alpha[v],v}_{t,t+\delta}\varphi(t+\delta,\cdot)(x)-\varphi(t,x)  \notag \\
	&=\int_{[t,t+\delta]}^\oplus l(x(s),\alpha[v](s))ds \oplus  \varphi(t+\delta,x(t+\delta)) 
		-\frac{1}{2}\int_t^{t+\delta}|v(s)|^2 ds -\varphi(t,x) \notag \\
	&= \varphi(t+\delta,x(t+\delta)) 
		-\frac{1}{2}\int_t^{t+\delta}|v(s)|^2 ds -\varphi(t,x).
	\label{F-eq}
\end{align}
Therefore, \eqref{F-est}, \eqref{F-eq} imply
\[
	F^{\alpha[v],v}_{t,t+\delta}\varphi(t+\delta,\cdot)(x)-\varphi(t,x)
	\leq \left( \frac{\partial \varphi}{\partial t}(t,x)
		+\mathcal{K}(x,\varphi(t,x)-\rho,\nabla \varphi(t,x))  \right) \delta
		+\eps \delta + o(\delta)
\]
Since $o(\delta)$ is uniform on $\| v \|_{L^2[t,t+\delta]}^2 \leq M$,
\[
	F^{EK}_{t,t+\delta}\varphi(t+\delta,\cdot)(x)- \varphi(t,x)
	\leq \left( \frac{\partial \varphi}{\partial t}(t,x)
		+\mathcal{K}(x,\varphi(t,x)-\rho,\nabla \varphi(t,x))  \right) \delta
		+\eps \delta + o(\delta).
\]
Dividing by $\delta$ and taking the limsup as $\delta \to 0+$,
\[
	\limsup_{\delta \to 0+} \frac{1}{\delta}
	\{ F^{EK}_{t,t+\delta}\varphi(t+\delta,\cdot)(x)- \varphi(t,x) \}
	\leq  \frac{\partial \varphi}{\partial t}(t,x)
		+\mathcal{K}(x,\varphi(t,x)-\rho,\nabla \varphi(t,x))  +\eps.
\]
Since $\eps>0$ is taken arbitrarily, we have
\[
\limsup_{\delta \to 0+} \frac{1}{\delta}
	\{ F^{EK}_{t,t+\delta}\varphi(t+\delta,\cdot)(x)- \varphi(t,x) \}
	\leq  \frac{\partial \varphi}{\partial t}(t,x)
		+\mathcal{K}(x,\varphi(t,x)-\rho,\nabla \varphi(t,x)).
\]
Finally, by sending $\rho \to 0$, we can prove \eqref{EK-gen_u}.

We next show \eqref{EK-gen_l}.
For $\alpha \in \Gamma_{EK}(t,t+\delta)$ and $v \in L^2[t,t+\delta]$,
\[
	F^{\alpha[v],v}_{t,t+\delta}\varphi(t+\delta,\cdot)(x)
	\geq \inf_{u \in L^\infty([t,t+\delta];U)} F^{u,v}_{t,t+\delta}\varphi(t+\delta,\cdot)(x).
\]
Then, we have
\[
	F^{EK}_{t,t+\delta}\varphi(t+\delta,\cdot)(x)
	\geq \sup_{v \in L^2[t,t+\delta]}
 \inf_{u \in L^\infty([t,t+\delta];U)} F^{u,v}_{t,t+\delta}\varphi(t+\delta,\cdot)(x).
 \]
From this inequality, 
\begin{multline}
 	\liminf_{\delta \to 0+} \frac{1}{\delta}\{ F^{EK}_{t,t+\delta}\varphi(t+\delta,\cdot)(x)
	 	-	\varphi(t,x) \} \\
	 \geq \liminf_{\delta \to 0+} \frac{1}{\delta}
	\left\{ \sup_{v \in L^2[t,t+\delta]}
	 \inf_{u \in L^\infty([t,t+\delta];U)} F^{u,v}_{t,t+\delta}\varphi(t+\delta,\cdot)(x)
 	-\varphi(t,x) \right\}.
\label{est_sup_inf}
\end{multline}
Let $\rho>0$ and fix $v \in \mathbb{R}^d$.
For $u(\cdot) \in L^\infty([t,t+\delta];U)$, 
let $x(s)$ be the solution of \eqref{system} with $u(s)$ and ${v}(s) \equiv {v}$.
We consider two cases.

For the first case, suppose $u(s) \in A(x,\varphi(t,x)+\rho)$ a.e.\,on $[t,t+\delta]$.
Then,
\begin{align}
	&F^{u,{v}}_{t,t+\delta}\varphi(t+\delta,\cdot)-\varphi(t,x)  \notag \\
	&\geq \varphi(t+\delta,x(t+\delta)) -\varphi(t,x) -\frac{1}{2}\int_t^{t+\delta}|{v}|^2ds
		\notag \\	
	&=\int_{t}^{t+\delta}
		\frac{\partial \varphi}{\partial t}(s,x(s))
		+(f(x(s),u(s))+\sigma(x(s),u(s)){v}) \cdot \nabla \varphi(s,x(s))
		 -\frac{1}{2}|{v}|^2 ds  \notag \\
	&\geq \int_t^{t+\delta}
	\frac{\partial \varphi}{\partial t}(t,x)
		+(f(x,u(s))+\sigma(x,u(s)){v}) \cdot \nabla \varphi(t,x)
		 -\frac{1}{2}|{v}|^2 ds +o(\delta)   \notag \\
	&\geq
		\left( \frac{\partial \varphi}{\partial t}(t,x)
		+\min_{u \in A(x,\varphi(t,x)+\rho)}
		\left\{ (f(x,u)+\sigma(x,u){v}) \cdot \nabla \varphi(t,x) 
		-\frac{1}{2}|{v}|^2 \right\} \right)\delta
		+o(\delta),
	\label{case1}
\end{align}
where $o(\delta)$ is uniform on $u(\cdot)$. 

For the second case, we suppose there exists $I_\delta \subset [t,t+\delta]$ with positive measure
such that $u(s) \not\in A(x,\varphi(t,x)+\rho)$ for any $s \in I_\delta$.
In a similar way to the argument in Theorem \ref{cal_gen}, we can see that
\[
	\int_{[t,t+\delta]}^\oplus l(x(s),u(s)) ds >\varphi(r,x(r))+\frac{\rho}{2}, 
	\ t \leq r \leq t+\delta
\]
for small $\delta$ uniform on $u(\cdot)$. Then, we have
\begin{multline}
	F^{u,v}_{t,t+\delta}\varphi(t+\delta,\cdot)(x)-\varphi(t,x) \\
	\geq
	\left( \frac{\partial \varphi}{\partial t}(t,x)
		+\min_{u \in A(x,\varphi(t,x)+\rho)}
		\left\{ (f(x,u)+\sigma(x,u){v}) \cdot \nabla \varphi(t,x) 
		-\frac{1}{2}|{v}|^2 \right\} \right)\delta
	\label{case2}
\end{multline}
for small $\delta$ uniform on $u(\cdot)$.

Combining \eqref{case1} and \eqref{case2},
\begin{multline*}
	\inf_{u \in L^\infty([t,t+\delta];U)} F^{u,v}_{t,t+\delta}\varphi(t+\delta,\cdot)(x)
	-\varphi(t,x) \\
	\geq 
	\left( \frac{\partial \varphi}{\partial t}(t,x)
		+\min_{u \in A(x,\varphi(t,x)+\rho)}
		\left\{ (f(x,u)+\sigma(x,u){v}) \cdot \nabla \varphi(t,x) 
		-\frac{1}{2}|{v}|^2 \right\} \right)\delta+o(\delta).
\end{multline*}
Thus, by \eqref{est_sup_inf}, we have
\begin{align*}
	&\liminf_{\delta \to 0+} \frac{1}{\delta}\{ F^{EK}_{t,t+\delta}\varphi(t+\delta,\cdot)(x)
	 	-	\varphi(t,x) \} \\
	&\geq \liminf_{\delta \to 0+} \frac{1}{\delta}
	\left\{ \sup_{v \in L^2[t,t+\delta]}
	 \inf_{u \in L^\infty([t,t+\delta];U)} F^{u,v}_{t,t+\delta}\varphi(t+\delta,\cdot)(x)
 	-\varphi(t,x) \right\} \\
	& \geq \liminf_{\delta \to 0+} \frac{1}{\delta}
	\left\{
	 \inf_{u \in L^\infty([t,t+\delta];U)} F^{u,v}_{t,t+\delta}\varphi(t+\delta,\cdot)(x)
 	-\varphi(t,x) \right\} \\
	&\geq 
	\frac{\partial \varphi}{\partial t}(t,x)
		+\min_{u \in A(x,\varphi(t,x)+\rho)}
		\left\{ (f(x,u)+\sigma(x,u){v}) \cdot \nabla \varphi(t,x) 
		-\frac{1}{2}|{v}|^2 \right\}.
\end{align*}
Since $v \in \mathbb{R}^d$ is taken arbitrarily, 
\[
	\liminf_{\delta \to 0+} \frac{1}{\delta}\{ F^{EK}_{t,t+\delta}\varphi(t+\delta,\cdot)(x)
	 	-	\varphi(t,x) \} 
	\geq \frac{\partial \varphi}{\partial t}(t,x)
		+\mathcal{K}(x,\varphi(t,x)+\rho, \nabla \varphi(t,x)).
\]
Taking $\rho \to 0$, we obtain \eqref{EK-gen_l}.


\end{document}